\documentclass{article}
\usepackage{
amsmath,
amsfonts,
latexsym, 
amssymb
}

\usepackage{bm}

\usepackage{tocloft}

\usepackage{cite}

\usepackage{pstricks}

\usepackage{subcaption}

\usepackage{comment}

\usepackage[hidelinks]{hyperref}   

\usepackage{tikz}
\usepackage{tikz,fullpage}
\usetikzlibrary{arrows,%
                petri,%
                topaths, automata}%
\usepackage{floatrow}

\usepackage{graphicx}
\usepackage{enumerate}
\usepackage{mathrsfs}
\usepackage{wrapfig}

\usepackage{color}

\usepackage{color}

\newcommand{\Id}{{\mathrm{Id}}}

\numberwithin{equation}{section}

\newcommand{\eps}{\varepsilon}

\newcommand{\rd}{{\rm d}}

\newcommand{\D}{\mathbb{D}}

\newcommand{\bw}{{\bf{w}}}

\newcommand{\bz} {{\bf {z}}}

\newcommand{\be}{\begin{equation}}
\newcommand{\ee}{\end{equation}}

\newcommand{\e}{{\varepsilon}}

\usepackage{amsmath} 
\usepackage{amssymb}
\usepackage{amsthm}

\newcommand{\wt}{\widetilde}

\newcommand{\ii}{\mathrm{i}} 
\newcommand{\dd}{\mathrm{d}}

\renewcommand{\epsilon}{\varepsilon}
\renewcommand{\leq}{\leqslant}
\renewcommand{\geq}{\geqslant}

\newcommand{\E}{\mathbb{E}}
\newcommand{\R}{\mathbb{R}}
\newcommand{\C}{\mathbb{C}}
\newcommand{\N}{\mathbb{N}}

\newcommand{\p}[1]{({#1})}

\DeclareMathOperator{\tr}{Tr}

\DeclareMathOperator{\supp}{supp}

\DeclareMathOperator{\re}{Re}
\DeclareMathOperator{\im}{Im}

\DeclareMathOperator{\OO}{O}
\DeclareMathOperator{\oo}{o}

\DeclareMathOperator{\cpcty}{cap}

\DeclareMathOperator{\dist}{dist}

\theoremstyle{plain} 
\newtheorem{theorem}{Theorem}[section]
\newtheorem*{theorem*}{Theorem}
\newtheorem{lemma}[theorem]{Lemma}
\newtheorem*{lemma*}{Lemma}
\newtheorem{corollary}[theorem]{Corollary}
\newtheorem*{corollary*}{Corollary}
\newtheorem{proposition}[theorem]{Proposition}
\newtheorem*{proposition*}{Proposition}

\newtheorem*{definition*}{Definition}

\newtheorem*{example*}{Example}
\newtheorem{remark}[theorem]{Remark}

\newtheorem*{remark*}{Remark}
\newtheorem*{remarks*}{Remarks}

\makeatletter
\renewcommand{\subsection}{\@startsection
{subsection}
{2}
{0mm}
{-\baselineskip}
{0 \baselineskip}
{\normalfont\bf\itshape}} 
\makeatother

\usepackage{dsfont}
\usepackage{stmaryrd}

\def\@empty{}

\def\author#1{\par
    {\centering{\authorfont#1}\par\vspace*{0.05in}}
}

\def\titlefont{\fontsize{13}{15}\bfseries\boldmath\selectfont\centering{}}
\def\authorfont{\fontsize{13}{15}}
\def\abstractfont{\fontsize{8}{10}}

\let\affiliationfont\rhfont

\def\address#1{\par
    {\centering{\affiliationfont#1\par}}\par\vspace*{11pt}
}

\def\body{
\setcounter{footnote}{0}
\def\thefootnote{\alph{footnote}}
\def\@makefnmark{{$^{\rm \@thefnmark}$}}
}

\def\title#1{
    \thispagestyle{plain}
    \vspace*{-14pt}
    \vskip 79pt
    {\centering{\titlefont #1\par}}%
    \vskip 1em
}

\renewenvironment{abstract}{\par%
    \vspace*{6pt}\noindent 
    \abstractfont
    \noindent\leftskip10pt\rightskip10pt
}{%
  \par}

\usepackage{tocloft}

\makeatletter
\renewcommand{\section}{\@startsection
{section}
{1}
{0mm}
{-2\baselineskip}
{1\baselineskip}
{\normalfont\large\scshape\centering}} 
\makeatother

\usepackage[non-sorted-cites]{amsrefs}
\usepackage{wrapfig}
\usepackage[none]{hyphenat}
\usepackage{caption}
\usepackage{titlesec}

\titleformat{\subsection}[runin]
  {\normalfont\bfseries\itshape}   
  {\thesubsection}                 
  {0.5em}                          
  {}
  []

\titlespacing*{\subsection}
  {0pt}                             
  {3.25ex plus 1ex minus .2ex}      
  {0.5em}                           

\titleformat{\subsubsection}[runin]
  {\normalfont\itshape}            
  {\thesubsubsection}
  {0.5em}
  {}
  []

\titlespacing*{\subsubsection}
  {0pt}
  {3.25ex plus 1ex minus .2ex}
  {0.5em}

\begin{document}

~\hspace{-0.3cm}
\title{
Fisher-Hartwig asymptotics for non-Hermitian random matrices}

\vspace{1cm}

\noindent\begin{minipage}[b]{0.19\textwidth}
 \author{Paul Bourgade}
\address{\scriptsize Courant Institute,\\ New York University}
 \end{minipage}
\noindent\begin{minipage}[b]{0.32\textwidth}
 \author{Guillaume Dubach}
\address{\scriptsize Centre~de~Mathématiques~Laurent~Schwartz,\\  École Polytechnique}
 \end{minipage}
\noindent\begin{minipage}[b]{0.24\textwidth}
 \author{Lisa Hartung}
\address{\scriptsize Institut für Mathematik,\\ Universität Mainz}
 \end{minipage}
\noindent\begin{minipage}[b]{0.23\textwidth}
 \author{Ahmet Keles}
\address{\scriptsize Courant Institute,\\ New York University}
 \end{minipage}

\begin{abstract}

\noindent We prove the two-dimensional analogue of the asymptotics for Toeplitz determinants with Fisher-Hartwig singularities, for general real symbols.
This formula has applications to random normal matrices with complex spectra:
(i) the characteristic polynomial converges to a Gaussian multiplicative chaos random measure on the limiting droplet, in the  subcritical phase;
(ii) the electric potential converges pointwise to a logarithmically correlated field; 
(iii) the measure of its level sets (i.e.  thick points) is identified;
(iv) the associated free energy undergoes a freezing transition.

This establishes emergence of the Liouville quantum gravity measure from free fermions in 2d,   and universality with respect to the external potential. 

\end{abstract}

\tableofcontents

\noindent

\section{Introduction}

Determinants of the form
\[
D_N(\mu)=\det\left(\int z^i\bar z^j\rd \mu(z)\right)_{0\leq i,j\leq N-1}   \ \ \ (\mu\ \mbox{a measure on}\  \mathbb{C})
\]
arise naturally in probability,   approximation theory,  and complex geometry.
In the primary settings,    $\mu$ is supported on either the unit circle (Toeplitz determinants) or the real line (Hankel determinants).  The case in which it has a smooth density except at a finite number of singularities is of particular interest in statistical physics,  as such singularities can manifest phase transitions.  These singularities can be either a  polynomial vanishing (root singularity) or a discontinuity (jump singularity), and the 
large $N$ evaluation of  $D_N$ in these settings -- commonly referred to as Fisher–Hartwig asymptotics -- is well understood (see Subsection \ref{subsec:FH}).

For general $\mu$ supported on $\mathbb{C}$ presenting singularities, much less is known.  These determinants are then  related to free fermions in 2d and random matrices with complex spectrum,  through Andreiev's identity:
\[
D_N(\mu)=\frac{1}{N!}\int_{\mathbb{C}^N} \prod_{j<k}|z_j-z_k|^2\prod_{i=1}^N\rd\mu(z_i).
\]
This paper gives asymptotics of $D_N$ for $\mu$ supported on the complex plane in the presence of root singularities.
A particular case of our main results gives the limiting fractional moments of the characteristic polynomial for the  Ginibre ensemble, 
defined as  an $N\times N$ random  matrix $G$ with independent complex Gaussian entries: 
$G_{ij}$ has distribution
\[
\frac{N}{\pi}\,e^{-N|z|^2}\,\rd m(z),\hspace{2cm} (\mbox{$\rd m$ is the Lebesgue measure on $\C$}).
\]
With this normalization of the matrix entries,  the empirical measure of eigenvalues converges to the uniform probability measure on the unit disc $\D$.
As a consequence of our asymptotics for fractional moments,  the characteristic polynomial converges to 
the Gaussian multiplicative chaos on $\mathbb{D}$ (see Subsection \ref{subsec:GMC}). These results are informally stated as follows.

\begin{theorem*}
For any fixed $m$, distinct $\zeta_j$ in the interior of $\mathbb{D}$ and positive exponents 
$\gamma_j$, $1\leq j\leq m$,
we have
\[
\E\Big[\prod_{j=1}^m|\det(G-\zeta_j)|^{\gamma_j}\Big]\underset{N\to\infty}{\sim}
\ \prod_{j=1}^m e^{\frac{\gamma_j}{2}N(|\zeta_j|^2-1)}N^{\frac{\gamma_j^2}{8}}
\frac{(2\pi)^{\frac{\gamma_j}{4}}}{{\rm G}\left(1+\frac{\gamma_j}{2}\right)}
\prod_{1\leq j<k\leq m}|\zeta_j-\zeta_k|^{-\frac{\gamma_j\gamma_k}{2}},
\]
where ${\rm G}$ is the Barnes ${\rm G}$-function.  
Such non-Hermitian Fisher-Hartwig asymptotics -- in the full form stated in Subsection \ref{subsec:FH} -- imply the following convergence on $\mathbb{D}$,  in distribution with respect to the weak topology and in the subcritical phase (i.e.  $\gamma<2\sqrt{2}$):
\[
\frac{|\det(G-z)|^{\gamma}}{\mathbb{E}[|\det(G-z)|^{\gamma}]}\,\rd m(z)  \to e^{\frac{\gamma}{\sqrt{2}} h(z)}\, \rd m(z),
\]
where $h$ is the Gaussian field with covariance $\E[h(w)h(z)]=-\log|z-w|$ on $\mathbb{D}^2$,  and $e^{\frac{\gamma}{\sqrt{2}} h(z)}\, \rd m(z)$ denotes the Gaussian multiplicative chaos random measure associated to $h$, with parameter $\frac{\gamma}{\sqrt{2}}$.
\end{theorem*}
The main Theorem \ref{Theorem.2singularities} below generalizes the above moment formula to the exponential generating function of the electric potential,  coupled with mesoscopic linear statistics,  for general confining external potentials. 
A remarkable aspect of this general formula is the relevance of the harmonic measure on the boundary of the droplet, with the following consequence:
the Gibbs measure associated to the electric potential converges to a Gaussian multiplicative chaos measure, which is
universal up to a global Gaussian shift inherited from  the capacity of the droplet (see Theorem \ref{GMCforGinibre}).

\begin{center}
\vspace{-0.3cm}
\includegraphics[width=0.3\linewidth]{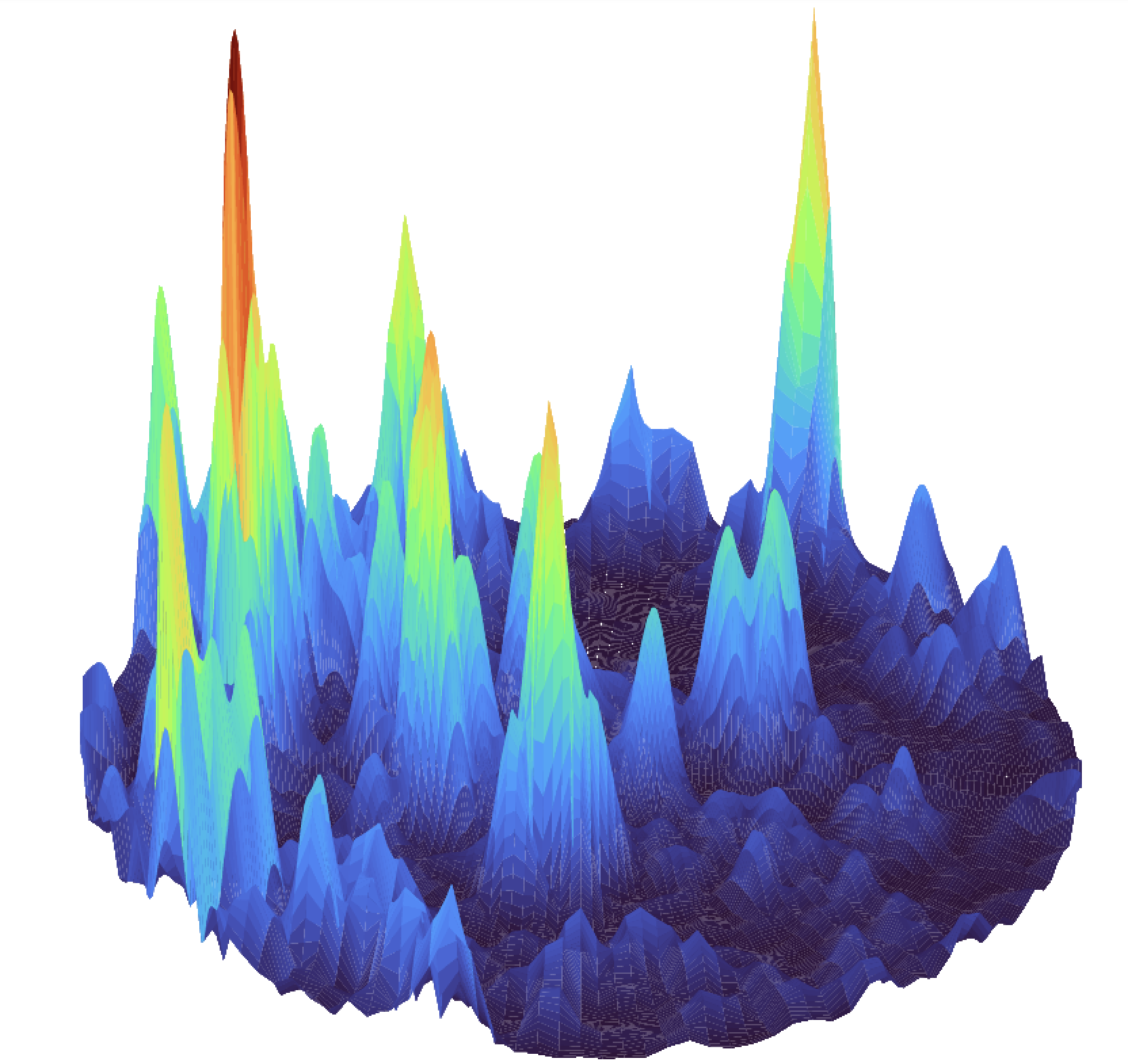}
\vspace{-0.2cm}
\captionof{figure}{The random field $\frac{|\det(G-z)|}{\E{|\det(G-z)|}}$ for $|z|<1$ and a Ginibre matrix $G$ of size $N=1000$.}
\end{center}

\subsection{Fisher-Hartwig asymptotics.}\ \label{subsec:FH} 
The main results of this paper are explicit large $N$ asymptotics of $D_N(\mu)$ in the 2d context, involving special functions that  parallel the Toeplitz and Hankel cases.
More precisely,  for Toeplitz determinants a simple version of the strong Szeg{\H o} theorem states that if $\rd\mu=e^V\frac{\rd\theta}{2\pi}$ with $V$ real-valued and smooth enough, 
\begin{equation}\label{eqn:SZ}
D_N(\mu)\sim \exp(N\hat V_0+\frac{1}{2}\|V\|^2_{{\rm H}^{1/2}})
\end{equation}
for large $N$,
where the Fourier transform is normalized as $\hat f_k=\int_0^{2\pi} f(e^{\ii\theta})e^{-\ii k\theta}\frac{\rd\theta}{2\pi}$ and  $\|f\|_{{\rm H}^{1/2}}^2=\sum_{n\in\mathbb{Z}} |n|\,|\hat f_n|^2$.
For $\mu$ with general  1d support,  similar formulas were obtained which connect $D_N(\mu)$ to
the Loewner energy of curves and Weil-Petersson quasicircles \cite{Joh1988,Joh2022,JohVik2023}.

In presence of singularities, Fisher and Hartwig \cite{FisHar} made a seminal general conjecture about the asymptotic form of $D_N(\mu)$ in the Toeplitz case, which has been corrected by Basor and Tracy \cite{BasTra} and is settled in full generality by Deift, Its and Krasovsky,  using Riemann-Hilbert methods  \cite{DeiItsKra}, after multiple important contributions, e.g.\cite{Wid,Bas,Ehr}.
For example, in the special case where $\rd\mu(z)=e^{V(z)}\prod_{j=1}^m|z-z_j|^{2\alpha_j}\frac{\rd\theta}{2\pi}$ with $m\geq 1$ fixed singularities $z_j$ on the unit circle, $\alpha_j>-1/2$,  and smooth centered real $V$,
the Fisher-Hartwig asymptotics  states that
\begin{equation}\label{eqn:Widom}
D_N(\mu)=e^{\frac{1}{2}\|V\|^2_{{\rm H}^{1/2}}-\sum_{j=1}^m\alpha_jV(z_j)}\, N^{\sum_{j=1}^m \alpha_j^2}
\prod_{1\leq j<k\leq m}|z_j-z_k|^{-2\alpha_j\alpha_k}\prod_{j=1}^m\frac{{\rm G}(1+\alpha_j)^2}{{\rm G}(1+2\alpha_j)}\ (1+\oo(1)),
\end{equation}
where the Barnes function ${\rm G}$ is defined in Subsection \ref{subsec:NotConv}.
Applications of general Fisher-Hartwig asymptotics for Toeplitz determinants (including the case of jump singularities) are multiple, we refer  for example to \cite{DeiItsKraII} for the relevance to phase transitions in the Ising model.  
The main contribution of our paper is the  analogue of (\ref{eqn:Widom}) -- a formula due to Harold Widom in 1973 -- for non-Hermitian random matrices, as explained below.\\

We consider the ensemble of $N\times N$ random normal matrices, defined on the manifold $\{M\in\C^{N\times N}:MM^*=M^*M\}$, endowed with the probability measure,
\begin{align*} 
\mathbb{P}^{V}(\rd M)\propto e^{-N\tr V(M)}\rd M
\end{align*}
where $\rd M$ is the measure induced by the Lebesgue measure on all $N\times N$ complex matrices and $V:\C\to\R\cup\{+\infty\}$ be a confining potential.  This induces the following probability measure on the eigenvalues of $M$:
\[
\frac{1}{Z_{N,V}} \prod_{1\leq j<k\leq N} |z_j - z_k|^2
\prod_{i=1}^N\mu(\rd z_i),\quad \mu(\dd z)= e^{-NV(z)} \dd m(z)
\]
where  $Z_{N,V}$ is the normalizing constant. In the following we will mostly suppress the dependence on $V$ and $N$, writing $\mathbb{P}=\mathbb{P}^{V}$ and $\mathbb{E}=\mathbb{E}^{V}$ for the corresponding expectation, and only make it explicit when necessary.

 Under suitable regularity and growth assumptions on $V$ (see Subsection \ref{subsec:potential}), the normalized empirical distribution of the particles $N^{-1}\sum\delta_{z_i}$ converges almost surely, as $N\to\infty$, to a compactly supported probability measure $\mu_V$, called the equilibrium measure. Its support has non-empty interior, denoted by $S$, which is referred to as the droplet. The equilibrium measure is then given by
\begin{align*} 
\rd\mu_{V}(z)=\rho_V(z)\,\rd m(z),\quad \rho_V(z)=\frac{\Delta V(z)}{4\pi}\mathds{1}_{S}(z),\ \ \ \mbox{($\Delta$ is the Laplacian in $\R^2$)}.
\end{align*}
The fluctuations around this deterministic limit were first fully understood by Rider and Vir\'{a}g \cite{RidVir2007} when the potential is quadratic, $V(z)=|z|^2$; then, the eigenvalues density coincides with the complex Ginibre ensemble
previously mentioned and $S=\mathbb{D}$.  Their central limit theorem is the 2d analogue of the strong Szeg{\H o} theorem \eqref{eqn:SZ}: for smooth enough $f$,
\begin{equation}\label{eqn:RiderVirag}
\E\Big[e^{{\rm Tr}f(G)-\E[{\rm Tr}f(G)]}\Big]\sim
e^{\frac{1}{8\pi}\int_{\mathbb{D}}|\nabla f|^2\rd m
+\frac{1}{4}\|f\|_{{\rm H}^{1/2}}^2
}.
\end{equation}
Ameur, Hedenmalm, and Makarov proved analogous results for general $V$ and rigorously introduced the method of Ward identities (also called loop equations) in this bidimensional setting \cite{AmeHedMak2011,AmeHedMak2015}. 
The main Theorem \ref{Theorem.2singularities} below is  an extension of the results in \cite{AmeHedMak2015} to include root singularities.\\

Before stating it, we introduce some notation and recall some facts from potential theory (see e.g. \cite[Chapters II-III]{GarMar2005}).
For a function $g$ with domain including $\bar{S}$, let $g^{S}$ denote the unique bounded harmonic extension of $g$ to $S^{\rm c}$, with $g^{S}=g$ on $\bar{S}$. Note that for this harmonic extension, $g^{S}(\infty):=\lim_{|z|\to\infty}g^{S}(z)$ exists and
\begin{align*} 
|g^{S}(z)-g^{S}(\infty)|=\OO\big(|z|^{-1}\big),\quad |\nabla^n g^{S}(z)|=\OO\big(|z|^{-n-1}\big)\quad \textnormal{for all}\ n\geq 1\ \textnormal{as}\ |z|\to\infty.
\end{align*}
Omitting the dependence on $S$ from the notation, we define the Neumann jump operator by $\mathcal{N}(g)=\frac{\partial g|_{S}}{\partial n}-\frac{\partial g^{S}|_{S^{\rm c}}}{\partial n}$ on $\partial S$ where $n$ is the exterior normal vector (this is the negative of the definition used in \cite{AmeHedMak2015}).

Denote $\log^{\zeta}=\log|\zeta-\cdot\,|$. The Green's function with a pole at $\infty$ on the exterior domain $S^{\rm c}$ is defined by (see e.g. \cite[Chapter II.2]{GarMar2005})
\begin{align*} 
G(z,\infty)=\log^{\zeta}z-(\log^{\zeta})^S(z), \quad \zeta\in S
\end{align*}
where the right hand side is independent of the choice of $\zeta\in S$ (see Appendix \ref{app:neumann}) and we omitted the dependence on $S$ from the notation. The capacity of the set $\bar{S}$ is defined by (see e.g. \cite[Chapter III.1]{GarMar2005}) 
\begin{align}\label{eqn:capacity_defn} 
\mathfrak{s}:=\log\cpcty(\bar{S})=(\log^{\zeta})^S(\infty)
\end{align}
where the right hand side is again independent of the choice of $\zeta\in S$ (see Appendix \ref{app:neumann}). 
Equivalently, $-\log\cpcty(\bar{S})$ is the infimum of the energy integral
\begin{align}\label{eqn:energy_functional} 
I(\mu)=\iint\log\frac{1}{|z-w|}\rd\mu(z)\rd\mu(w)
\end{align}
over all probability measures $\mu$ supported on $\bar{S}$, \cite[Thm. III.4.1]{GarMar2005}. The infimum is uniquely attained by the harmonic measure seen from infinity which is supported on the boundary $\partial S$ and given by
\begin{align}\label{eqn:Neumann_log}
\rd\omega^{\infty}(z)=\frac{1}{2\pi}\frac{\partial G(z,\infty)}{\partial n}\rd s=\frac{1}{2\pi}\mathcal{N}(\log^{\zeta})(z)\rd s, \quad z\in\partial S
\end{align}
where $\rd s$ denotes arc-length measure on $\partial S$, \cite[Thm. II.2.5]{GarMar2005}. See Appendix \ref{app:neumann} for a self-contained proof of the independence of the right-hand side from the choice of $\zeta\in S$ and the fact that $\omega^{\infty}$ is a probability measure. Equivalently, $\omega^{\infty}$ is characterized probabilistically as follows. Let $B^{z}$ be planar Brownian motion started from $z\in\C$, and $\tau$ denote the hitting time of $\partial S$. Then, $\omega^{\infty}$ coincides with the weak limit of the law of $B^{z}_\tau$ as $|z|\to\infty$. Finally, for any $\zeta\in S$, the Poisson-Jensen formula applied to $(\log^{\zeta})^S$ yields
\begin{align} \label{eqn:lop_capacity}
\mathfrak{s}=(\log^{\zeta})^S(\infty)=\int_{\partial S}\log^{\zeta} (z)\,\rd\omega^{\infty}(z),\quad \textnormal{for all}\ \zeta\in S.
\end{align}


We  also need the notations  $L=\log\frac{\Delta V}{4}$, and refer to Subsection \ref{subsec:NotConv} for the definition of the functional space $\mathscr{S}_{n,C,\kappa}$,  to state the following 
main theorem.  It gives the asymptotics for the joint Laplace transform of the electrostatic potential,  or equivalently the joint moments of the characteristic polynomial of random normal matrices. 

\begin{theorem}[General potential]\label{Theorem.2singularities}
Let $C>1$, $m\in\N$, $0<\kappa\leq\frac{1}{100C}$ be fixed constants. Then, uniformly in $f\in\mathscr{S}_{5,C,\kappa}$, $\zeta_1,\dots,\zeta_m\in \{\zeta\in S:\dist(\zeta,\partial S)>N^{-1/2+\kappa}\}$ satisfying $\min_{1\leq j<k\leq m}|\zeta_j-\zeta_k|>N^{-1/2+\kappa}$ and $\gamma_1,\dots,\gamma_m\in[0,C]$ we have
(${\rm G}$ stands for the Barnes ${\rm G}$-function, see Subsection \ref{subsec:NotConv})
\begin{align*}
\E\Big[e^{\sum_{i=1}^Nf(z_i)}\prod_{j=1}^m\prod_{i=1}^{N}|z_i-\zeta_j|^{\gamma_j}\Big]=
e^{N\int f\rd\mu_{V}
+\frac{1}{8\pi}\int_{\C}(|\nabla f^S|^2
+\mathds{1}_{S}\Delta f
+\Delta f L^S)\rd m
}
\prod_{j=1}^{m}e^{\frac{\gamma_j}{2}\big(\int_{\partial S}f\rd\omega^{\infty}-f(\zeta_j)\big)} 
\\
\times\prod_{j=1}^m e^{\gamma_jN\int\log^{\zeta_j}\rd\mu_{V}}N^{\frac{\gamma_j^2}{8}}
\frac{(2\pi)^{\frac{\gamma_j}{4}}}{{\rm G}\left(1+\frac{\gamma_j}{2}\right)}
e^{\frac{\gamma_j}{4}(L(\zeta_j)-L^{S}(\infty))+
\frac{\gamma_j^2}{8}(L(\zeta_j)+2\mathfrak{s})}
\\
\times\prod_{1\leq j< k\leq m}|\zeta_j-\zeta_k|^{-\frac{\gamma_j\gamma_k}{2}}
e^{\frac{\gamma_j\gamma_k}{2}\mathfrak{s}}
\left(1+\OO(N^{-\kappa/10})\right).
\end{align*}
\end{theorem}

The Ginibre ensemble is the paradigmatic non-Hermitian random matrix model; see \cite{ByuFor} for its many facets.
For this ensemble, our asymptotics take a simpler form. We have $\int_{\D} \log^{\zeta}\frac{\rd m}{\pi}=\frac{|\zeta|^2-1}{2}$ and choosing $\zeta=0$ gives $\mathcal{N}(\log^{\zeta})=1$ and $\int_{\partial\D}\log^{\zeta}\rd s=0$,  so $\log\cpcty(\bar{\D})=0$. Hence the main theorem takes the following simplified expression,  conjectured by Webb and Wong in \cite{WebWong2019}.

\begin{theorem}[Ginibre ensemble]\label{Theorem.2singularitiesGin}
Let $G$ be a complex Ginibre matrix with eigenvalues $z_i$. Fix constants $C>1$, $m\in\N$, $0<\kappa\leq\frac{1}{100C}$. Then, uniformly in $f\in\mathscr{S}_{5,C,\kappa}$, $\zeta_1,\dots,\zeta_m\in{\rm B}(0,1-N^{-1/2+\kappa})$ satisfying $\min_{1\leq j<k\leq m}|\zeta_j-\zeta_k|>N^{-1/2+\kappa}$ and $\gamma_1,\dots,\gamma_m\in[0,C]$ we have
\begin{multline*}
\E\Big[e^{{\rm Tr}f(G)}\prod_{j=1}^m|\det(G-\zeta_j)|^{\gamma_j}\Big]=
e^{\frac{N}{\pi}\int_{\mathbb{D}}f\rd m
+\frac{1}{8\pi}\int_{\mathbb{D}}|\nabla f|^2\rd m
+\frac{1}{4}\|f\|_{{\rm H}^{1/2}}^2
+\frac{1}{8\pi}\int_{\mathbb{D}}\Delta f\rd m
}
\prod_{j=1}^{m}e^{\frac{\gamma_j}{2}(\hat{f}_0-f(\zeta_j))} 
\\
\times\prod_{j=1}^m e^{\frac{\gamma_j}{2}N(|\zeta_j|^2-1)}N^{\frac{\gamma_j^2}{8}}
\frac{(2\pi)^{\frac{\gamma_j}{4}}}{{\rm G}\left(1+\frac{\gamma_j}{2}\right)}
\prod_{1\leq j<k\leq m}|\zeta_j-\zeta_k|^{-\frac{\gamma_j\gamma_k}{2}}
\left(1+\OO(N^{-\kappa/10})\right).
\end{multline*}
\end{theorem}

\begin{remark}[Complex exponents]\label{rem:Extensions}
The case $m=1,  V(z)=|z|^2, f=0$ was proved in \cite{WebWong2019} with a different method which analyses a related Riemann-Hilbert problem \cite{BalBerLeeMcL},  and allows  
any complex exponent ${\rm Re}\,\gamma>-2$.  The current  proof of Theorem \ref{Theorem.2singularities} extends to complex $\gamma_i$'s in a neighbourhood of $[0,C]$, 
but the statement is expected to hold uniformly in compact subsets of ${\rm Re}\, \gamma_k> -2$, $1\leq k\leq m$.  We do not address
this  here but note potential applications to fluctuations of the measure of thick points of the field, see Remark \ref{rem:FluctField}.
\end{remark}

\begin{remark}[Jump singularities]
In the Toeplitz case,  the most general version of Fisher-Hartwig asymptotics \cite{DeiItsKra} allows jump singularities with general complex exponents,  with asymptotics involving a subtle variational problem.  
As a natural 2d analogue, one may consider the  asymptotics of $\E\Big[e^{\sum_{j=1}^m\gamma_j\#\{z_i\in\Omega_j\}}\Big]$.
No conjecture is known for general $\Omega_j\subset\mathbb{C}$ and complex $\gamma_j$.  Results were obtained for rotationally invariant ensembles, concentric rings
$\Omega_j$'s,  and exponents in a neighbourhood of $\mathbb{R}$, see  \cite{Cha2022} and references therein.
\end{remark}

\begin{remark}[Painlev\'e transcendents]
For merging singularities on microscopic scale (namely $\zeta_j-\zeta_k\asymp N^{-1/2}$),   the joint moments of characteristic polynomials are expected to give rise to Painlev\'e transcendents,  mirroring a similar phenomenon identified for Toeplitz determinants in \cite{ClaItsKra2015}.  In dimension 2,  this was recognized by Dea\~{n}o and Simm,  who proved it by algebraic methods for the Ginibre ensemble,  integer moments, with no linear statistics $f$. Then
  Painlev\'e V appears for two merging singularities in
the bulk (see \cite[Theorem 1.5]{DeaSim2019}), and  Painlev\'e IV on the boundary of the disk (see  \cite[Theorem 1.3]{DeaSim2019}).
In a forthcoming work \cite{LeeYang2025},  these asymptotics are extended to fractional moments in the bulk, 
for $V(z)=|z|^2$ and $f=0$.
The methods developed for Theorem \ref{Theorem.2singularities} have potential to prove the emergence of these Painlev\'e transcendents for general $V,f$,  by comparison with \cite{DeaSim2019,LeeYang2025},
but we do not pursue this goal here.
\end{remark}

A straightforward corollary of Theorem \ref{Theorem.2singularities} is the following central limit theorem for the log-characteristic polynomial. Gaussian fluctuations of the determinant of random matrices have been studied extensively, both for  1d \cite{KeaSna2000, TaoVu2014, BouMod2019, BMP2022,Mod2023} and 2d spectra \cite{Gir1981, TaoVu2014,SerSim2023}.  Our main theorem proves such fluctuations for random normal matrices, generalizing \cite[Corollary 1.2]{WebWong2019} to the multidimensional setting and general $V$.

\begin{corollary}[Log-correlations pointwise]\label{cor:CLT}
Fix  $m\in\N$ and $\kappa\in(0,1/2)$. Assume that the (possibly $N$-dependent) points $\zeta_1,\dots,\zeta_m\in \{\zeta\in S:\dist(\zeta,\partial S)>N^{-1/2+\kappa}\}$ satisfy $\min_{1\leq j<k\leq m}|\zeta_j-\zeta_k|>N^{-1/2+\kappa}$ and the limit
$
c_{j,k}=\lim_{N\to\infty}\frac{-\log|\zeta_j-\zeta_k|}{2\log N}
$
exists for every $j\neq k$. Let $\Sigma=(c_{i,j})_{i,j=1}^{m}$ where $c_{i,i}=\frac{1}{4}$ for each $i=1,\dots,m$. Then
\begin{align*} 
\left(\frac{\sum_{i=1}^{N}\log|\zeta_j-z_i|-N\int \log |\zeta_j-z|\rd\mu_V(z)}{\sqrt{\log N}}\right)_{j=1}^m\xrightarrow[N\to\infty]{\rm (d)}\mathcal{N}(0,\Sigma).
\end{align*}
\end{corollary}

In the case $m=1$, \cite[Theorem 1.2]{AmeKanSeo}  shows that such a CLT also holds at points where the limiting density $\Delta V$ may vanish or diverge,  with a variance depending on the order of the singularity.\\

To conclude this subsection, we note that moments of random characteristic polynomials of wide classes of random matrices have been a topic of major interest, see e.g.
 \cite{AkeVer,FyoSup,BorStr,ForRai2009,BaiKea,Afa2019,SerSim2025}  in the case of integer exponents by algebraic and supersymmetric  methods, and \cite{Kra2007,ItsKra,Cha,BerWebWon2017,WebWong2019,ChaFahWebWon,CFLW2021,ByuSeoYan,DeaMclMolSim2025} for fractional exponents by Riemann-Hilbert methods.  
Motivations for the study of random determinants in the non-Hermitian setting include for example the stability of large complex systems \cite{BenFyoKho2021}.
Theorems \ref{Theorem.2singularities} and \ref{Theorem.2singularitiesGin} introduce joint fractional moments for non-Hermitian ensembles, motivated by two-dimensional random geometry from free fermions, as explained below.

\subsection{Gaussian multiplicative chaos.}\  \label{subsec:GMC}
The Gaussian multiplicative chaos (GMC) appeared in attempts to model intermittency in fluid mechanics, with its premises in studies by Kolmogorov \cite{Kolmogorov1941,Kolmogorov1962},   Obukhov \cite{Obukhov1962} and Mandelbrot \cite{Mandelbrot1,Mandelbrot1974}.
Its general,  rigorous foundations were led by Kahane in \cite{Kah1985}: in his theory,  the Gaussian multiplicative chaos is the fractal measure on $\mathbb{R}^d$ defined 
by the density
\[
e^{\gamma X(z)}   = \lim_{\eps \to 0} e^{\gamma X_{\eps}(z) - \frac{\gamma^2}{2} \E(X_\eps(z)^2)},
\]
with respect to the Lebesgue measure,
where $X_{\eps}$ is a mollification of a $\log$-correlated Gaussian field $X$ (meaning $\E[X(v)X(w)]=-\log|v-w|+{\rm O(1)}$ as $v\to w$) on a domain $D \subset \mathbb{R}^d$. The regularization and renormalization are necessary because of the negative Sobolev regularity of the field. The convergence holds in probability with respect to the topology of weak convergence,  the limit does not depend on the mollification and it is non-trivial  for any $\gamma \in (0, \sqrt{2d})$, see for example \cite{RhodesVargas}
for these properties. The specific case where $X$ is a two dimensional Gaussian free field  (a Gaussian field whose covariance function is  the inverse of the Laplacian)  has proved to be connected with many different domains in mathematical physics. To name a few, it is the volume form in Liouville quantum gravity, a metric measure space corresponding to the formal Riemannian metric tensor $``e^{\gamma X} (\rd x^2 +\rd y^2)"$  \cite{DDDF2020, GwyMil2021}, it appears in the scaling limit of random planar maps \cite{LeGall2013, MilShe2020, HolSun2019}, and it plays a central role in the rigorous approach to Liouville Conformal Field Theory \cite{DavKupRhoVar2016, DOZZ}. We refer to the surveys \cite{She2023,RhodesVargas2025} for more references on these connections.\\

In \cite{Web2015}, Webb connected the GMC and random matrix theory, showing that the characteristic polynomial of random unitary matrices evaluated on the unit circle converges to a one-dimensional GMC.  He conjectured that similar results also hold for the Gaussian ensembles,   one-dimensional $\beta$-ensembles, and more generally for random matrix models presenting log-correlations, including in dimension two.  His proof and the ones of the following works \cite{BerWebWon2017, NikSakWeb2018,CFLW2021} relied on Fisher-Hartwig asymptotics based on the Riemann-Hilbert method. 
Another approach based on sparse matrix models showed that the limit of  the spectral measure of circular $\beta$-ensembles,  and the characteristic polynomial,  also converge to a GMC,  still for $d=1$  \cite{ChaNaj,LamNaj2024}. 

In two-dimensional space-time settings when dynamics are involved, the works \cite{BouFal2025, Keles2025} obtained convergence of the eigenvalues counting function and the electric potential to a 2d GMC, 
for the geometry of the 2d torus (unitary Brownian motion) and the 2d strip (non-intersection Brownian motions), respectively.
These works rely on a general surgery method replacing the Riemann-Hilbert approach when there is no known equivalence to asymptotics of orthogonal polynomials.  

In the 2d ambient space the Liouville quantum gravity measure also appears as stated in Theorem \ref{GMCforGinibre} below.  The key ingredient is the precise asymptotics of the moments of the characteristic polynomial, Theorem \ref{Theorem.2singularities} above.
To state our result, note that from Corollary \ref{cor:CLT}, the log-characteristic polynomial of a random normal matrix is $\frac{1}{2}$-log-correlated, so the corresponding critical value becomes $\gamma_{\rm c}=2\sqrt{2}$. The regime $\gamma<2$ is often referred to as the $L^{2}$-phase, while the full subcritical range $\gamma<2\sqrt{2}$ is called the $L^{1}$-phase. All arguments below are carried out in the $L^{1}$-phase.

We consider the field $X_N$ associated to the characteristic polynomial,
\begin{align} \label{eqn:defn_X_N}
X_{N}(z)=\log|\det(M-z)|-\E\big[\log|\det(M-z)|\big],
\end{align}
and note that the fine asymptotics of the centering term (as follows from the proof of Theorem \ref{Theorem.2singularities}):  for some $\epsilon>0$
\[
\E\big[\log|\det(M-z)|\big]=N\int\log^{z}\rd\mu_{V}+\frac{1}{4}+\frac{1}{4}(L(z)-L^S(\infty))+\OO(N^{-\epsilon}).
\]
 At the level of fluctuations,  $\sqrt{2}\,X_{N}$ converges to a Gaussian free field $X$ (see  \cite{AmeHedMak2011}) with covariance kernel $T_{\mathfrak{s}}$ on $S^2$ defined by 
\begin{equation*}
T_{\mathfrak{s}}(z,w) = \frac{1}{2\pi}\int_{\C}\nabla(\log^{z})^{S}\cdot\nabla(\log^{w})^{S}\rd m=\log\frac{1}{|z-w|}+\mathfrak{s},
\end{equation*}
where the second equality is  shown in Appendix \ref{app:neumann}. Note that $\mathfrak{s}=\log\cpcty(\bar{S})$ is the smallest real number for which the kernel $T$ is positive semi-definite, since $-\log\cpcty(\bar{S})$ is the infimum of the energy integral given in \eqref{eqn:energy_functional}.

Let $\nu_{\gamma,(\mathfrak{s})}$ be the GMC measure associated to $X$, formally given by
\[
\nu_{\gamma,(\mathfrak{s})}(A)=\int_{A}e^{\gamma X(z)-\frac{\gamma^2}{2}\E[X(z)^2]}\rd m(z).
\]
More precisely,  $\nu_{\gamma,(\mathfrak{s})}$ the limit,  as $\varepsilon\to 0$, of $e^{\gamma X_\varepsilon(z)-\frac{\gamma^2}{2}\E[X_\varepsilon(z)^2]}\rd m(z)$, where $X_\e$ is the convolution of $X$ with a bump function on scale $\e$.  As already mentioned,  this limit is  non-trivial and does not depend on the convolution kernel,  see e.g.  \cite[Theorem 3.1]{RhodesVargas}.
Then  the random Borel measure
\begin{equation}\label{eqn:GMCN}
\nu_{\gamma,N} (A) = \int_A \frac{e^{\gamma X_N(z)}}{\E[e^{\gamma X_N(z)}]}\rd m(z)
\end{equation}
converges to a GMC measure in the full $L^1$-phase,  as conjectured in \cite{Web2015,WebWong2019,Lam2020}.

\begin{theorem}[Convergence to a GMC]\label{GMCforGinibre}
For any $0<\gamma<2\sqrt{2}$,  as $N\to \infty$ we have the following convergence in distribution with respect to the weak topology on $S$:
\[
\nu_{\gamma,N} \to \nu_{\gamma/\sqrt{2},(\mathfrak{s})}.
\]
\end{theorem}
\noindent More explicitly, the above convergence means that
$
\int_{S}\varphi\,\rd\nu_{\gamma,N}\xrightarrow[N\to\infty]{{\rm (d)}}\int_{S}\varphi\,\rd\nu_{\gamma/\sqrt{2},(\mathfrak{s})}
$
for any continuous bounded $\varphi$ supported in $S$.
Moreover, we also note that the denominator in \eqref{eqn:GMCN} is explicit thanks to Theorem \ref{Theorem.2singularities}.

\begin{remark}[Outside the droplet] The field $X_N$ behaves very differently on $S^{\rm c}$ and the fractal nature of the limit disappears: in the case of the Ginibre ensemble,  $\det(\Id-z^{-1}G)$ converges in law of the exponential of a random analytic function for $|z|<1$, as proved in \cite{BorChaGar2022} in the general setting of random matrices with i.i.d. entries.
\end{remark}

\begin{remark}[The Coulomb gas]
A natural extension of the model considered is the two-dimensional Coulomb gas confined by the external field $V$,  at general inverse temperature $\beta$, i.e. the Gibbs measure
\[e^{-\beta H}\ \ \mbox{where}\ \ H(\bz)=-\sum_{1\leq j<k\leq N} \log|z_j - z_k|
+\frac{N}{2}\sum_{i=1}^N V(z_i).\]
In this setting, distributional convergence of  $X_N(z)=\sum\log|z-z_i|-\E\sum\log|z-z_i|$ to the Gaussian free field was proved in  \cite{LebSer2018,BauBouNikYau2019},
suggesting convergence of $e^{\gamma X_N}$ to a GMC.  More precisely, in view of Theorem \ref{GMCforGinibre},  one expects
\[
\frac{e^{\gamma X_N(z)}}{\E[e^{\gamma X_N(z)}]}\to \nu_{\gamma/\sqrt{\beta},(\mathfrak{s})}
\]
in the ${\rm L}^1$ phase, which corresponds to $\gamma\leq 2\sqrt{\beta}$.
It is an interesting problem to prove the above convergence without relying on an analogue of Theorem \ref{Theorem.2singularities}, which seems out of reach with current methods.
\end{remark}

\begin{remark}[Dynamics]
An extension of \eqref{eqn:RiderVirag} was recently obtained in \cite{BouCipHua2024},  for the Ornstein-Uhlenbeck process on the space of non-Hermitian matrices. 
In this setting the field $X_{N}(z,t)=\log|\det(G_t-z)|-\E\big[\log|\det(G_t-z)|\big]$ converges to a Gaussian field $X$ on $\mathbb{D}\times \mathbb{R}$ which is log-correlated for the parabolic distance:
\[
\E[X_N(z,t)X_N(w,s)]\to \E[X(z,t)X(w,s)]=-\frac{1}{2}\log(|z-w|+\sqrt{t-s})+{\rm O}(1),
\]
raising the question of the convergence of $|\det(G_t-z)|^{\gamma}\rd m(z)\rd t$ to a related $3d$ GMC.
\end{remark}

Another corollary of our main theorems concerns the thick points of the field, where $X_N$ attains unusually large values, which play a central role in the geometry of extremes for log-correlated fields. It follows directly from \cite[Proposition 3.8]{CFLW2021}. It is the 2d analogue of  \cite[Theorem 1.3]{ArgBelBou2017} which considers high points of the circular unitary ensemble.

\begin{corollary}[Thick points]
For every compact $K\subset S$ with non-empty interior,  $\gamma\in[0,\frac{1}{\sqrt{2}})$ and $\epsilon>0$, we have
\begin{align*} 
\lim_{N\to\infty}\mathbb{P}\Big(N^{-2\gamma^2-\epsilon}\leq\big|\big\{z\in K:\,X_N(z)\geq\gamma\log N\big\}\big|\leq N^{-2\gamma^2+\epsilon}\Big)=1.
\end{align*}
\end{corollary}

\begin{remark}[Fluctuations of thick points]\label{rem:FluctField} The fluctuations of $N^{2\gamma}\big|\big\{z\in K:\,X_N(z)\geq\gamma\log N\big\}\big|$ are supposedly related to the mass of $\nu_{\gamma/\sqrt{2},(\mathfrak{s})}$.  Indeed  \cite{JunLamWeb2024} 
gives a general criterion for this correspondence,  and applies it to identify the fluctuations of the measure of high points of the characteristic polynomial of the circular unitary ensemble,  which had been conjectured in \cite{FyoKea2014}.
In our setting,   the criterion from \cite{JunLamWeb2024} requires an
extension of  Theorem \ref{Theorem.2singularities} to complex exponents $\zeta_i$,  see Remark \ref{rem:Extensions}.  However, while the mass of the Gaussian multiplicative chaos measures on the unit circle has an explicit distribution \cite{FyoBou,Rem},
to the best of our knowledge there is no conjectural explicit density in the 2d case.
\end{remark}

From the above corollary,  by following the same steps as in \cite[Proof of Corollary 1.4]{ArgBelBou2017}, we obtain the following consequence for the free energy. This extends the identity \cite[(1.15)]{Lam2020} to the non-Hermitian setting with a general potential\footnote{As a minor remark, \cite{Lam2020} appears to have omitted the factor $N$ in the numerator inside the logarithm.}. 

\begin{corollary}[Freezing]
For every compact $K\subset S$ with non-empty interior and $\gamma>0$, the free energy converges   in probability as $N\to\infty$:
\begin{align*} 
\frac{\log\big(N\int_{K}e^{\gamma X_N(z)}\rd m(z)\big)}{\gamma\log N}\to \begin{cases}
\frac{1}{\gamma}+\frac{\gamma}{8},&\mbox {for}\ \gamma\leq2\sqrt{2},\\
\frac{1}{\sqrt{2}},&\mbox{for}\ \gamma\geq 2\sqrt{2}.
\end{cases}
\end{align*}
\end{corollary}

\noindent The above transition of the free energy appears in statistical physics for models exhibiting the {\it freezing scenario}, which is instrumental in the analysis of equilibrium Gibbs measures with logarithmic spatial correlations \cite{CarPen2001}.

\begin{remark}[Extreme values]\label{rmk:max_log_char}
Our main theorems also allow to recover the leading order asymptotics for the electric field for  general potential:
For $K$ any  compact subset of $S$ with non-empty interior and any $\epsilon>0$,
\begin{equation}\label{eqn:Max} 
\lim_{N\to\infty}\mathbb{P}\Big(\max_{z\in K}\frac{X_N(z)}{\log N}\in[\frac{1}{\sqrt{2}}-\epsilon,\frac{1}{\sqrt{2}}+\epsilon]\Big)=1.
\end{equation}
The leading-order asymptotics for $X_N$ has already been established for the Ginibre ensemble in \cite{Lam2020}, for general $\beta$ with quadratic potential in \cite{LamLebZei2024}, for general $\beta$ and potential in \cite{Pei2025},
and for i.i.d. matrices in \cite{CipLan2024}. 

With Theorem \ref{GMCforGinibre}, 
the lower bound in the above equation follows directly from the support of the GMC on so-called $\gamma$-thick points (see  \cite[Theorem 3.4]{CFLW2021}),  while the upper bound follows directly from  Theorem \ref{Theorem.2singularities},  as explained in \cite[Proof of Proposition 3.1]{Lam2020}. 
A natural problem is refining the method for Theorem \ref{Theorem.2singularities} (closer singularities,  complex exponents) for finer orders in \eqref{eqn:Max}.
\end{remark}

\subsection{Outline.}
The proof of the main Theorem \ref{Theorem.2singularities} is based on the integration of two techniques,  the Ward identities as introduced in \cite{AmeHedMak2015} for non-Hermitian matrices,  and
the removal of singularities,  first implemented in space-time in \cite{BouFal2025}.  The surgery from  \cite{BouFal2025} can be applied to a variety of models involving local, special factors and long-range interactions,
but  implementation in different settings presents specific difficulties.  This method proceeds in several steps:
\begin{enumerate}[(1)]
\item Cut the long range nonsingular part of the determinants,  and prove a decoupling of the resulting product of localized root singularities.
\vspace{-0.15cm}
\item Establish a general ``gluing operation" for nonsingular terms. 
\vspace{-0.15cm}
\item Evaluate asymptotics of \textit{one} localized singularity, by gluing the opposite of the associated long range nonsingular part to the determinant,  evaluated for one specific model where integrability holds.
\vspace{-0.15cm}
\item Glue back the nonsingular parts of the determinant,  and the additional smooth function $f$ to the localized singularities.
\end{enumerate}

Step (1) requires some decay of correlations,  obtained here from  a general multiplicative comparison of Fredholm determinants introduced in \cite{BouFal2025} and kernel estimates
from \cite{AmeHedMak2010}, see Lemma \ref{lem:decoupling}.  At the technical level,  this step requires an a priori submicroscopic smoothing of the log-singularity (Proposition  \ref{prop:submic_log_reg})
and goes beyond the usual decoupling methods for determinantal point processes which give additive error terms instead of multiplicative.  

Steps (2) and (4) rely on
the method of Ward identities,  or loop equations,  which is common  in field theories.  For Hermitian random matrices,  the seminal work \cite{Joh1998} relied on loop equations
to reduce Laplace transform estimates to a shift of the density of states.   In our non-Hermitian setting, Ward identities are also instrumental as in \cite{AmeHedMak2015},  but the control of the error  due to potential anisotropy (i.e. local directional contribution from the interactions) and the boundary of the droplet,  is  delicate for measures biased by determinant powers, i.e.   in the presence of root singularities.  
Regarding the proof of local isotropy -- a problem naturally absent in the settings of \cite{BouFal2025,Keles2025} --  Lemma \ref{lem:anis} first treats nonsingular potentials,  by a moments method specific to determinantal point processes but flexible enough to cover mesoscopic observables up to the boundary of the droplet; this estimate is then transferred to singular setting,  as it holds with overwhelming probability (here  $1-e^{-(\log N)^{D}}$ for any $D>0$), while in  \cite{AmeHedMak2015} angular cancellations hold in expectation for the nonsingular potentials.
About the effect of the boundary, one novelty of our analysis is a careful control of the norms of the harmonic extensions of functions on mesoscopic scales, and an analysis of the effect of bulk singularities on the boundary field. This analysis is significant due to the harmonic measure on the boundary of the droplet,  which generates an additional 
 Gaussian shift when compared to the free field 
with free boundary conditions,  harmonically extended outside of the  droplet.

For step (3), while the emergence of the Barnes special function was rooted in the Selberg integrals in the space-time settings \cite{BouFal2025,Keles2025}, in this article it relies on 
a result by Kostlan stating that the radial parts of the eigenvalues of a Ginibre matrix are independent,  see Lemma \ref{lemma:stab_reg_log}.  In \cite{BouFal2025} this  step 
dealt with only one matrix ensemble, the CUE, due to rotational invariance; in \cite{Keles2025} this step required comparison between two different ensembles (GUE at arbitrary energy level and CUE); in the present work, this comparison is pushed further,  matching the Ginibre ensemble to any random normal matrix model,  as needed towards the proved universality of the limiting GMC.

These steps are implemented in the next Sections \ref{eqn:SecNon} and  \ref{sec:sing_pot} to prove Theorem \ref{Theorem.2singularities}.  Section \ref{sec:GMC} then proves Theorem \ref{GMCforGinibre},  the convergence of
the Gibbs measure associated to the electric potential to a GMC,  based on a general criterion from \cite{CFLW2021}.

\subsection{Assumptions on the potential.}\label{subsec:potential}\ We impose the following assumptions (A1)-(A4) on the external potential $V$.

\begin{itemize}
\item[(A1)] Growth: $V\in\mathscr{C}^2(\C)$ and there exists an $\epsilon>0$ such that $\liminf_{|z|\to\infty}\big(V(z)-(2+\epsilon)\log|z|\big)\geq 0$ for some $\epsilon>0$. For convenience we also assume that there exists a constant $C>0$ such that $|V(z)|+|\nabla V(z)|\leq C(1+|z|)^C$ for all $z\in\C$; noting that this additional condition can be relaxed easily.
\end{itemize}

\noindent Under (A1), standard results in logarithmic potential theory (see e.g. \cite{SaffTotik2024LogPot}) guarantee that the energy functional
\begin{align*} 
\mathcal{I}_V(\mu)=\iint \log\frac{1}{|z-w|}\rd\mu(z)\rd\mu(w)+\int V(z)\,\rd\mu(z)
\end{align*}
admits a unique minimizer, called the equilibrium measure $\mu_V$, among all probability measures on $\C$. Its support is compact and on the support its density is $\frac{\Delta V}{4\pi}$. Moreover, $\mu_{V}$ is the limiting distribution of the normalized empirical spectral distribution for random normal matrices with confining potential $V$.
\begin{itemize}
\item[(A2)] Boundary regularity: The boundary of $\supp(\mu_{V})$ is a real-analytic Jordan curve and the interior of $\supp(\mu_{V})$, denoted by $S$, is simply connected.

\item[(A3)] Regularity: $V$ is real-analytic in a neighbourhood of $\supp(\mu_{V})$.
\end{itemize}

\noindent We define,
\begin{align*} 
\check{V}(z)=\sup\{f(z):f\in\mathscr{C}^2(\C)\ {\rm is\ subharmonic},\ f(z)\leq2\log|z|+\OO_{|z|\to\infty}(1),\ f\leq V\}.
\end{align*}
It is known (e.g. \cite{AmeHedMak2015}) that $\check{V}\in\mathcal{C}^{1,1}(\C)$, satisfies $\check{V}=V$ on $S$, is subharmonic on $\C$, harmonic on $\C\setminus S$ and
\begin{align} 
\check{V}(z)&=2\log |z|+\OO_{|z|\to\infty}(1),\nonumber
\\
\partial \check{V}(z)&=\int \frac{\rd\mu_V(w)}{z-w},\quad \textnormal{for all}\ z\in\C,\label{eqn:EulLag}
\end{align}
where the second equation is obtained by the Euler-Lagrange equation if $z\in S$ and uniqueness of the harmonic extension if $z\in S^{\rm c}$. 

\begin{itemize}
\item[(A4)] Non-degeneracy: $\Delta V(z)>0$ on the coincidence set $\{z\in\C:\check{V}(z)=V(z)\}$ and the coincidence set coincides with $\supp(\mu_{V})$. 
\end{itemize}

\subsection{Notations and conventions.} \label{subsec:NotConv}
\ We collect here the notations used throughout the paper.

\subsubsection{Scales and geometry.} We use the symbols $\alpha$ and $\kappa$ to denote small fixed constants in the exponents. We define the mesoscopic scale $\delta$, the submicroscopic scale $\Delta$, and the slightly super-microscopic scale $\delta_N$ as follows:
\begin{align*} 
\delta=N^{-1/2+\kappa}, \quad
\Delta=N^{-1/2-\alpha}, \quad
\delta_N=(\log N)^2N^{-1/2}.
\end{align*}
We let $\tau$ vary from mesoscopic to macroscopic scales, i.e., $\tau\in[\delta,1]$. For any parameter $a\geq0$, we define the inner, outer, and edge neighbourhoods of the droplet $S$:
\begin{gather*}
S_{\rm in}(a)=\{z\in S:\dist(z,\partial S)>a\},
\\
S_{\rm out}(a)=\{z\in S^{\rm c}:\dist(z,\partial S)>a\},
\\
S_{\rm edge}(a)=\{z\in\C:\dist(z,\partial S)< a\}.
\end{gather*}

\subsubsection{Function classes and differential operators.} Denote by $\mathscr{C}^n(A)$ the class of real-valued functions on $A\subseteq\C$ whose partial derivatives up to order $n$ are continuous. We utilize the complex differential operators $\partial=\frac{1}{2}(\partial_{x}-i\partial_{y})$, $\bar\partial=\frac{1}{2}(\partial_{x}+i\partial_{y})$, and the Laplacian $\Delta=4\partial\bar{\partial}$. For $f\in\mathscr{C}^n(\C)$ we define
 \begin{align*}
|\nabla^n f(z)|=\max_{n_1+n_2=n}|\partial_x^{n_1}\partial_y^{n_2} f(z)|,\quad \|\nabla^n f\|_{\infty}=\sup_{z\in\C}|\nabla^{n} f(z)|. 
 \end{align*}
For any $\epsilon>0$ and $n\in\mathbb{N}$, we define the scale-dependent norm
\begin{equation}\label{eqn:norm}
\|f\|_{n,\epsilon}=\sum_{k=0}^n \frac{1}{\epsilon^k}\|\nabla^k f\|_\infty, \ \ \textnormal{where}\ \|\nabla^n f\|_{\infty}=\max_{n_1+n_2=n}\sup_{z\in\C}|\partial_x^{n_1}\partial_y^{n_2} f(z)|.
\end{equation}
We denote by $\mathscr{A}_{n,\epsilon}$ the collection of smooth functions that are supported in a ball of radius $\epsilon$ and are smooth up to the $n^{\rm th}$ derivative on that scale:
\begin{equation*} 
\mathscr{A}_{n,\epsilon}=\Big\{g\in\mathscr{C}^{n}(\C)\,\Big|\,\textnormal{$g$ is supported in ${\rm B}(z_0,\epsilon)$ for some $z_0\in\C$, and $\|g\|_{n,\epsilon}\leq 1$}\Big\}.
\end{equation*}
Additionally, for constants $C>0$ and $\kappa>0$, we define a collection of (sequence of) functions $\mathscr{S}_{n,C,\kappa}$ by:
\begin{equation*}
\mathscr{S}_{n,C,\kappa}=\Bigg\{f=f(N)\,\Big|\,f=\sum_{i=1}^{\lfloor C\log N\rfloor}g_i\ \textnormal{where}\ g_i\in\mathscr{A}_{n,\epsilon_i},\ \epsilon_i\in[N^{-1/2+\kappa},1]\Bigg\}.
\end{equation*}

\subsubsection{Cutoff functions and log-regularization.} Unless otherwise stated, the symbol $\chi$ denotes a fixed smooth radial cutoff function with $\chi=1$ on ${\rm B}(0,1/2)$, $\chi=0$ outside ${\rm B}(0,1)$ and let $\chi_\epsilon(z)=\chi(z/\epsilon)$. For any $\epsilon>0$, we define log-regularization on scale $\epsilon$ by 
\begin{align*} 
\log_{\epsilon}=\log|\cdot|\ast\frac{\chi_{\epsilon}}{\|\chi_{\epsilon}\|_{L^1}}.
\end{align*}
For shifted functions, we use the shorthand
\begin{align*} 
\chi_{\epsilon}^{\zeta}=\chi_{\epsilon}(\cdot-\zeta),\quad \log_{\epsilon}^{\zeta}=\log_{\epsilon}(\cdot-\zeta).
\end{align*}

\subsubsection{Biased measures.} We denote the particles by $z_i$ and the singularities by $\zeta_j$. The linear statistics and centered linear statistics are given by:
\begin{align*} 
\tr(f)=\sum_{i=1}^{N}f(z_i),\quad X_f=\sum_{i=1}^{N}f(z_i)-N\int_{S}f\rho_V\,\rd m
\end{align*}
where $\rho_{V}=\frac{\Delta V}{4\pi}$ on $S$. Finally, we define the biased probability measure $\mathbb{P}_f$ and the corresponding expectation $\mathbb{E}_f$ by
\begin{align}\label{eqn:bias_defn} 
\mathbb{P}_{f}(A)=\E\big[\mathds{1}_{A}\frac{e^{\tr(f)}}{\E [e^{\tr(f)}]}\big].
\end{align}

\subsubsection{Barnes $\rm{G}$-function.} The Barnes G-function is defined as the Weierstrass product
\begin{align*}
{\rm G}(z+1)=(2\pi)^{z/2}e^{-\frac{z+z^2(1+\gamma)}{2}}\prod_{k=1}^\infty \left(1+\frac{z}{k}\right)^ke^{\frac{z^2}{2k}-z}
\end{align*}
where $\gamma$ is the Euler-Mascheroni constant.
It satisfies the functional equation
$
G(z+1)=\Gamma(z)G(z)
$
where $\Gamma$ is the Gamma function.

\subsubsection*{Acknowledgements.}  
We wish to thank Fredrik Viklund for enlightening comments about harmonic measures.
P. B. and A. K. were supported by the NSF standard grant DMS-2348202. G. D. gratefully acknowledges support from the Fondation de l'\'{E}cole Polytechnique, as well as from the Agence Nationale de la Recherche (ANR-25-CE40-5672 and ANR-25-CE40-1380).
The research of L. H. was supported by the Deutsche Forschungsgemeinschaft (DFG, German Research
Foundation) in the framework to the collaborative research center "Multiscale Simulation Methods for Soft-Matter Systems" (TRR 146) under Project No. 233630050 and Project-No. 443891315 within SPP 2265.

\section{Preliminaries on nonsingular potentials}\label{eqn:SecNon}

In this section, we gather preliminary estimates for nonsingular potentials, in preparation for the next section treating root singularities.   Many relevant methods 
appeared in \cite{AmeHedMak2015},  to obtain concentration bounds (Proposition \ref{prop:exp_tail}) and fluctuations (Proposition \ref{prop:Laplace_loop}) for linear statistics of the $z_i$'s. 
We follow a similar strategy based on the Ward identities,  with additional ideas from \cite{BauBouNikYau2017} to cover linear statistics on any mesoscopic scale,  including at the edge of the droplet.

These improvements to mesoscopic scales are instrumental in the next section,   to treat smoothed versions of $\log$ by multiscale decompositions,  and to eventually allow the main Theorem \ref{Theorem.2singularities} to cover the electric potential up to the mesoscopic distance $N^{-1/2+\kappa}$ from the boundary.

\subsection{Ward identity.}\label{subsec:ward}\ We denote
\begin{align*}
\hat \mu_N=\frac{1}{N}\sum \delta_{z_k},\quad \mu_{V}=\frac{\Delta V}{4\pi}\mathds{1}_{S}\,\rd  m,\quad \tilde\mu_N=
\hat \mu_N-\mu_{V}
\end{align*}
for any $V\in\mathscr{C}^2(\C)$ satisfying 
\begin{align}
&\liminf_{z\to\infty}\frac{V(z)}{2\log|z|}>1,\label{eqn:V1}\\
&|V(z)|+|\nabla V(z)|\leq C(1+|z|)^C\label{eqn:V2}
\end{align}
for some fixed constant $C>0$. In this subsection, we do not impose the non-degeneracy and boundary regularity assumptions from Subsection \ref{subsec:potential}; these will be assumed in all other sections.

Moreover, assume that
$h:\mathbb{C}\to\mathbb{C}$ is a  function such that
\begin{align}
&h\ {\rm is\ continuous\ on}\ \mathbb{C},\label{eqn:h3}\\
&h\mid_{\mathbb{D}}\ {\rm  and}\  h\mid_{\mathbb{D}^{\rm c}}\ {\rm  are\ } \mathscr{C}^1,\label{eqn:h2}\\
&\|h\|_\infty+\|(\nabla h)\mid_{\mathbb{D}}\|_\infty+\|(\nabla h)\mid_{\mathbb{D}^{\rm c}}\|_\infty<C\label{eqn:h1}
\end{align}
for some fixed $C$, possibly depending on $N$.
For any
$\bold z\in(\mathbb{C}\backslash\partial S)^N$,  we define  
\begin{equation}\label{e:wvplus}
W_V^{h}({\bf z})= \sum_{j < k} \frac{h(z_j)-h(z_k)}{z_j-z_k}+\sum_j \partial h(z_j) - N \sum_j h(z_j) \partial V(z_j).
\end{equation}
The following lemma, sometimes called the Ward identity or loop equation, is analogous to \cite[Proposition 2.1]{AmeHedMak2015}, where it first appeared in the context of 2d Coulomb gases and was used to study Gaussian fluctuations of $\hat\mu_N$. It was also instrumental in   
\cite{BauBouNikYau2019} to prove such fluctuations at any inverse temperature. Its relation to Conformal Field Theory is discussed in \cite{KanMak2013}.

In \cite{AmeHedMak2015}, it was stated for Lipschitz and compactly supported test function. We give another version, with assumptions closer to
\cite[Lemma 8.3]{BauBouNikYau2017}, and essentially reproduce the proof from \cite{BauBouNikYau2017} by integration by parts.

\begin{lemma}\label{lem:expectation0}
Assume $V$ and $h$ are as described above. Then we have
\begin{align*}
\E^{V}\left[W_V^{h}\right]=0.
\end{align*}
\end{lemma}

\begin{proof}
Let $H(\bz)=-2\sum_{j<k}\log|z_i-z_j|+N\sum_{i=1}^NV(z_i)$. 
Let $\delta_\e$ be a smooth Dirac approximation on scale $\e$ and $h_\e=h*\delta_\e$.
By integration by parts, for any $j\in\llbracket1,N\rrbracket$, we have
\begin{equation*}
\E^{V}\left[\partial h_\e(z_j)\right]
=
\E^{V}\left[h_\e(z_j)\partial_{z_j} H(\bold z)\right].
\end{equation*}
 Summing this over $j\in\llbracket 1,N\rrbracket$ gives
\begin{multline*}
\sum_{j=1}^N\E^{V}[\partial h_\e(z_j)]
=\E\Big [\sum_{j=1}^Nh_\e(z_j)\Big ( N\partial V(z_j)-2\sum_{k\neq j}\partial{z_j}\log|z_j-z_k|\Big )\Big ]\\
=\E^{V}\left[N \sum_j h_\e(z_j) \partial V(z_j)-\sum_{j\neq k}\frac{h_\e(z_j)-h_\e(z_k)}{z_j-z_k}\right].
\end{multline*}
We now take $\e\to 0$ and conclude by dominated convergence. Indeed, first $\partial h_\e$ converges almost everywhere to $\partial h$ and $\|\nabla h_\e\|_\infty<C$ by \eqref{eqn:h1}.
Second, $|h_\e\partial V|<C|\partial V|$ which is integrable by \eqref{eqn:V1} and \eqref{eqn:V2}. Finally
$
\frac{h_\e(z_j)-h_\e(z_k)}{z_j-z_k}
$
is uniformly bounded thanks to \eqref{eqn:h3} and \eqref{eqn:h1}.
This concludes the proof.
\end{proof}

The following classical decomposition will also be useful in the next key proposition.

\begin{lemma}\label{lem:decompgeneral} For any $h$, $V$ as above, and $\bold z\in(\mathbb{C}\backslash\partial S)^N$, we have
\begin{multline}\label{Wdef}
\sum_j \int \frac{h(w)}{z_j-w}\mu_V(\rd w) - N\iint\frac{h(w)}{z-w}\mu_V(\rd w)\mu_V(\rd z)
=
- \frac{1}{N} W_V^h ({\bf z})
+ \frac{1}{N}  \sum_j \partial h(z_j)
+\sum h(z_j)(\partial \check{V}(z_j)-\partial V(z_j))\\
+ \frac{N}{2} \iint_{z\neq w} \frac{h(z)-h(w)}{z-w} \, \tilde \mu_N(\rd z) \, \tilde \mu_N(\rd w).
\end{multline}
\end{lemma}

\begin{proof}
The proof only relies on the definition of $W_V^h$ and \eqref{eqn:EulLag}. With \eqref{eqn:EulLag}, we can write the left-hand side of \eqref{Wdef} as 
\begin{multline*}
N\iint \frac{h(w)-h(z)}{z-w}\hat \mu_N(\rd z)\mu_V(\rd w)+\sum_j h(z_j)\partial \check{Q}(z_j)-\frac{N}{2}\iint\frac{h(w)-h(z)}{z-w}\mu_V(\rd w)\mu_V(\rd z)\\
=-\frac{1}{N}\sum_{j<k}\frac{h(z_j)-h(z_k)}{z_j-z_k}+\sum_j h(z_j)\partial \check{V}(z_j)+\frac{N}{2} \iint_{ z \not = w } \frac{h(z)-h(w)}{z-w}  \tilde \mu_N(\rd z) \tilde \mu_N(\rd w),
\end{multline*}
which is equivalent to \eqref{Wdef} for $\bold z\in(\mathbb{C}\backslash\partial S)^N$, where $W_V^h$ is well defined.
\end{proof}

For any given function $f$ of class $\mathscr{C}^{5}(\C)$ with at most logarithmic growth as $|z|\to\infty$ and $|\nabla f(z)|\leq C(1+|z|)^C$ (for a possibly $N$-dependent $C$), we denote the probability measure
\[
\rd\mu_{f}(\boldsymbol{z})=\frac{1}{Z_{N}^{f}}\prod_{i<j}|z_i-z_j|^2
e^{-N\sum_{i=1}^N (V-\frac{f}{N})(z_i)}\prod_{i=1}^N\rd m(z_i),
\]
and write $\mathbb{P}_{f}=\mathbb{P}^{V-\frac{f}{N}}$, $\E_{f}=\E^{V-\frac{f}{N}}$ for probability and expectation with this biased measure. Note that this coincides with the biased measure definition in \eqref{eqn:bias_defn}.

By definition \eqref{e:wvplus}, we have
\begin{align*} 
W_{V-\frac{f}{N}}^h=W_{V}^h+\sum_{j}h(z_j)\partial f(z_j).
\end{align*}
Thus, Lemma \ref{lem:expectation0} implies:
\begin{align*} 
\E_{f}[W_{V}^h]=-\E_{f}\Big[\sum_{j}h(z_j)\partial f(z_j)\Big].
\end{align*}
Therefore, taking the expectation of \eqref{Wdef} with respect to the biased measure $\E_{f}$ we obtain:
\begin{multline} \label{eqn:ward_main}
\E_{f}\Big[\sum_j \int \frac{h(w)}{z_j-w}\mu_V(\rd w)\Big]+\E_{f}\left[\sum h(z_j)(\partial V(z_j)-\partial \check{V}(z_j))\right]-N\iint\frac{h(w)}{z-w}\mu_V(\rd w)\mu_V(\rd z)
\\
= \frac{1}{N}\E_{f}\Big[\sum_{j}h(z_j)\partial f(z_j)\Big]+\frac{1}{N}\E_{f}\Big[\sum_j \partial h(z_j)\Big]
+ \E_{f}\Big[\frac{N}{2} \iint_{z\neq w} \frac{h(z)-h(w)}{z-w} \, \tilde \mu_N(\rd z) \, \tilde \mu_N(\rd w)\Big].
\end{multline}

\subsection{Harmonic extension and decomposition.}\ 
We define an analogous version of the $\|\cdot\|_{n,\epsilon}$ norm in \eqref{eqn:norm}  on a curve: For any smooth Jordan curve $\gamma$ parametrized by $\gamma:[0,2\pi]\to\C$, and any smooth function $g$ on the curve $\gamma$ we define,
\begin{align*} 
\|g\|_{\gamma,n,\epsilon}=\sum_{k=0}^{n}\frac{1}{\epsilon^k}\big\|\frac{\rd^k}{\rd\theta^k}(g\circ\gamma)\big\|_\infty.
\end{align*}
Although the definition depends on the choice of parametrization, we omit this dependence from the notation. We fix an order $1$ parametrization of $\partial S$ and will not need to refer to it again.

Let $\kappa>0$, $\tau\in[N^{-1/2+\kappa},1]$, $n\geq 2$ and $g\in\mathscr{A}_{n,\tau}$ be fixed. We denote the unique bounded harmonic extension of $g|_{\partial S}$ into $S^{\rm c}$ by $g_{\rm out}$. As $S$ is a simply connected open set there exists a biholomorphic map $\varphi:\bar{S}^{\rm c}\to\bar{\D}^{\rm c}$ satisfying $\varphi(\infty)=\infty$. Since $\partial S$ is real-analytic, the map $\varphi$ can be extended to a $\mathscr{C}^{\infty}$-diffeomorphism from $S^{\rm c}$ onto $\D^{\rm c}$ by the Kellogg–Warschawski Theorem (see, e.g. \cite[Theorem 5.2.4]{krantz2007}) and further to a $\mathscr{C}^{\infty}$ function on the whole complex plane $\mathbb{C}$. 
We define,
\begin{align*} 
G=g\circ\varphi^{-1},\quad \textnormal{on}\ \partial\D.
\end{align*}
Observe that since $g\in\mathscr{A}_{n,\tau}$, we have $\|g\|_{\partial S,n,\tau}<C$ for some constant $C$ depending only on $S$. By $\mathscr{C}^{\infty}$-diffeomorphism of $\varphi$, it follows that $\|G\|_{\partial \D,n,\tau}<C'$, where $C'$ is a constant depending only on $\varphi$. In addition, we may assume that $G$ is supported on an arc of length $\tau$ and $\|G\|_{\partial\D,n,\tau}\leq 1$, since $G$ can be written as a sum of finitely many such functions, the number of which depends only on $S$ and $\varphi$.

 If we denote the Fourier coefficients of $G$ by $\hat{G}_k=\frac{1}{2\pi}\int_{0}^{2\pi}G(e^{\ii\theta})e^{-\ii k\theta}\rd\theta$, the bounded harmonic extension of $G$ on $\D^{\rm c}$, denoted by $G_{\rm out}$, becomes
\begin{align*} 
G_{\rm out}(z)= \cdots+\frac{1}{z^2}\hat{G}_{-2}+\frac{1}{z}\hat{G}_{-1}+\hat{G}_0+\frac{1}{\bar{z}}\hat{G}_{1}+\frac{1}{\bar{z}^2}\hat{G}_{2}+\cdots ,\quad {\rm if}\ |z|\geq 1.
\end{align*}
Note that $g_{\rm out}=G_{\rm out}\circ\varphi$ on $S^{\rm c}$ by the uniqueness of harmonic extension. Decompose $G_{out}$ into two functions,
\begin{align}\label{eqn:pm_decomp}
G_{+}(z)=\sum_{k=0}^{\infty}\frac{1}{z^k}\hat{G}_{-k},\quad G_{-}(z)=\sum_{k=1}^{\infty}\frac{1}{\bar{z}^k}\hat{G}_k,\quad \textnormal{for}\ z\in\D^{\rm c}.
\end{align}
The following lemma proves that the regularity of $G$ on $\partial \D$ is preserved for the functions $G_{+}$ and $G_{-}$.

\begin{lemma}\label{lem:HarmReg}
Let $0<\epsilon<1$ and assume $f:\partial\D\to\C$ is smooth function supported on a curve of length $\epsilon$ with $\|f\|_{\partial\D,n,\epsilon}<1$ for an $n\geq 2$. Let $f_{+}$ and $f_{-}$ be defined as in \eqref{eqn:pm_decomp} on $\D^{\rm c}$. Then, for every $\ell\in\{0,1,\dots,n-1\}$ and $1<|z|<2$ we have
\begin{align} \label{eqn:HarmReg}
|\nabla^\ell f_{+}(z)|+|\nabla^\ell f_{-}(z)|\leq C_{\ell}\frac{\epsilon}{\max\big(\epsilon,\dist(z,\supp(f)\big)^{\ell+1}}
\end{align}
where $C_{\ell}$ is an absolute constant depending only on $\ell$. Moreover, for any $|z|>2$,
\begin{align} \label{eqn:HarmReg2}
|\nabla^\ell f_{+}(z)|+|\nabla^\ell f_{-}(z)|\leq C_{\ell}\,\epsilon.
\end{align}
\end{lemma}

\begin{proof} For any $|z|>1$,
\begin{align}\label{eqn:Poisson} 
f_{+}(z)=\int_{0}^{2\pi}f(e^{\ii\theta})\frac{z}{z-e^{\ii\theta}}\frac{\rd\theta}{2\pi},\quad f_{-}(z)=\int_{0}^{2\pi}f(e^{\ii\theta})\frac{e^{-\ii\theta}}{\bar{z}-e^{-\ii\theta}}\frac{\rd\theta}{2\pi},
\end{align}
where two integrands add up to the Poisson kernel on the unit circle. We discuss the proof for $f_{+}$ here in detail, for $f_{-}$ the proof is identical.

Fix an arbitrary point $z_0$ in the support $f$. If $|z-z_0|>2\epsilon$ and $|z|<2$,
\begin{align*} 
\partial^{\ell} f_{+}(z)=\int_{0}^{2\pi}f(e^{\ii\theta})\partial^{\ell}\big(\frac{z}{z-e^{\ii\theta}}\big)\frac{\rd\theta}{2\pi}\leq C_{\ell}\frac{\epsilon}{|z-z_0|^{\ell+1}}.
\end{align*}

We now move on to the case $|z-z_0|<2\epsilon$. Note that there exists universal constants $a_{k,\ell},1\leq \ell\leq k$ such that for any $|z|>1$ we have
\begin{equation}\label{eqn:partial}
\partial^k\frac{z}{z-e^{\ii\theta}}=\frac{1}{z^k}\sum_{\ell=1}^ka_{k,\ell}\partial_\theta^\ell \frac{z}{z-e^{\ii\theta}}.
\end{equation}
From \eqref{eqn:Poisson}, \eqref{eqn:partial} and integration by parts, we obtain
\begin{equation*}
\partial^k f_{+}(z)=
\frac{1}{z^k}\sum_{\ell=1}^k (-1)^\ell a_{k,\ell}\int_{0}^{2\pi}\frac{\partial_\theta^\ell f(e^{\ii\theta})}{1-\frac{1}{z} e^{\ii\theta}}\frac{\rd \theta}{2\pi}.
\end{equation*}
This implies that
\begin{equation*}
|\partial^k f_{+}(z)|\leq C_k \sup_{1\leq \ell\leq k}\left|\int_{0}^{2\pi}\frac{\partial_\theta^\ell f(e^{\ii\theta})}{1-\frac{1}{z} e^{\ii\theta}}\frac{\rd \theta}{2\pi}\right|.
\end{equation*}

Note that the Cauchy transform
$\int_{0}^{2\pi}\frac{\partial_\theta^\ell f(e^{\ii\theta})}{1-\frac{1}{z} e^{\ii\theta}}\frac{\rd \theta}{2\pi}$ has boundary 
\begin{align*}
\lim_{r\to 1^-}\int_{0}^{2\pi}\frac{\partial_\theta^\ell f(e^{\ii\theta})}{1- re^{\ii(\theta-\varphi)}}\frac{\rd \theta}{2\pi}=\frac{1}{2}\partial_\theta^\ell f(e^{\ii\theta})
+\lim_{\delta\to 0}\int_{|\theta-\varphi|>\delta}\frac{\partial_\theta^\ell f(e^{\ii\theta})}{1- e^{\ii(\theta-\varphi)}}\frac{\rd \theta}{2\pi}.
\end{align*}
We bound the above right-hand side by writing, for any $0<\delta<\epsilon$,
\begin{multline*}
\left|\int_{|\theta-\varphi|>\delta}\frac{\partial_\theta^\ell f(e^{\ii\theta})}{1- e^{\ii(\theta-\varphi)}}\rd \theta\right|\leq \int_{\delta <|\theta-\varphi|<\epsilon}\frac{|\partial_\theta^\ell f(\theta)-\partial_\theta^\ell f(\varphi)|}{|e^{\ii\theta}-e^{\ii\varphi}|}\rd \theta+\|\partial_\theta^\ell f\|_\infty\int_{\substack{|\theta-\varphi|>\epsilon\\ \theta\in\supp(f)}}\frac{\rd\theta}{|e^{\ii\theta}-e^{\ii\varphi}|}
\\
 \leq 
\epsilon\|\partial_\theta^{\ell+1} f\|_\infty+c\|\partial_\theta^\ell f\|_\infty =O(1/\epsilon^\ell).
\end{multline*}
This completes the proof of \eqref{eqn:HarmReg} as $\bar\partial f_{+}=0$ on $\D^{\rm c}$. The second estimation \eqref{eqn:HarmReg2} follows directly from the the maximum principle for harmonic functions.
\end{proof}

Let $\eta$ be a smooth cutoff function that is equal to $1$ on ${\rm B}(0,2)$ and vanishes on ${\rm B}(0,3)^{\rm c}$. Pick an arbitrary point $z_0\in\supp(G)$. Let $\sum_{k=1}^{\infty}\chi_k$ be a partition of unity such that $\chi_1$ is supported in ${\rm B}(z_0,2\tau)$ and $\chi_k$ is supported in ${\rm B}(z_0,2^{k+1}\tau)\setminus{\rm B}(z_0,2^{k-1}\tau)$ with $\|\chi_k\|_{n,2^k\tau}<C$ for some universal constant $C$. By Lemma \ref{lem:HarmReg}, the decompositions
\begin{align*}
G_{+}\eta=\sum_{j}G_{+}\eta\chi_j,\quad G_{-}\chi=\sum_{j}G_{-}\eta\chi_j
\end{align*}
contain $\OO(\log N)$ many (non-zero) terms satisfying: $G_{+}\eta\chi_j$ supported in a ball of radius $2^{j+1}\tau$ and smooth up to the $(n-1)^{\rm st}$ derivative on that scale (similarly for $G_{-}\eta\chi_j$).

Now, assume $f:\D^{c}\to\C$ is supported in a ball of radius $\epsilon$ that is smooth on that scale up to the $k^{\rm th}$ derivative. Then we can extend the domain of $f$ to include $S$ preserving its regularity as follows:
\begin{align*} 
f(re^{\ii\theta})=\left(\sum_{j=0}^{k}\partial_{r}^{j}f(re^{\ii\theta})\big|_{r=1^{+}}\frac{(r-1)^j}{j!}\right)\chi_{\epsilon}(r-1).
\end{align*}
This construction extends $f$ to a function on $\C$ supported in a ball of radius $C\epsilon$, satisfying  $\|\nabla^{j} f\|_{\infty}\leq C/\epsilon^{j}$ for every $j=0,\dots,k$ for some universal constant $C$. Applying this procedure to $G_{+}\eta\chi_j$ and $G_{-}\eta\chi_j$, and then gluing the resulting extensions together, we obtain smooth extensions of $G_{+}$ and $G_{-}$ into $\D$ while preserving their regularity. More precisely, if $g\in\mathscr{A}_{n,\tau}$ for a $\tau\in[N^{-1/2+\kappa},1]$, then $G_{+}\eta$ and $G_{-}\eta$ can be decomposed into sum of $\OO(\log N)$ many functions each of which is supported in a ball of radius $\epsilon_j$ and smooth up to $(n-1)^{\rm st}$ derivative on that scale for an $\epsilon_j\in[\tau,1]$.

We define 
\begin{align*} 
g_{+}=G_{+}\circ\varphi,\quad g_{-}=G_{-}\circ\varphi\quad \textnormal{on}\ S^{\rm c}.
\end{align*}
We now take a biholomorphic map $\psi$ from $S$ to $\D$. As $\partial S$ is real-analytic, $\psi$ extends holomorphically to an open neighbourhood of $\bar{S}$. Using this map, we construct $g_{+,j}=(G_{+}\eta\chi_j)\circ\psi$ and $g_{-,j}=(G_{-}\eta\chi_j)\circ\psi$ on that open neighbourhood of $\bar{S}$. Gluing  $g_{+,j}$'s together (and $g_{-,j}$'s) we obtain an extension of $g_{+}$ (and $g_{-}$) to the whole complex plane such that $g_{+}\chi$ and $g_{-}\chi$ can be decomposed as a sum of $\OO(\log N)$ many functions each of which is supported in a ball of radius $\epsilon_j$ and smooth up to the $(n-1)^{\rm st}$ derivative on that scale for an $\epsilon_j\in[\tau,1]$; where $\chi$ is an order $1$ bump function that is $1$ on $S\cup S_{\rm edge}(1/3)$ and 0 on $S_{\rm out}(1/2)$.

Note that, by construction, $\bar{\partial}g_{+}=\partial g_{-}=0$ and $g_{\rm out}=g^{S}=g_{+}+g_{-}$ on $S^{\rm c}$. We extend the definition of $g_{\rm out}$ to all $\C$ by setting $g_{\rm out}=g_{+}+g_{-}$ everywhere. Recall that, in contrast, $g^S$ was defined to coincide with $g$ on $S$. We summarize the main properties in below lemma.

\begin{lemma}\label{lem:harm_ext_reg}
Let $\kappa>0$, $\tau\in[N^{-1/2+\kappa},1]$, $\|f\|_{n,\tau}\leq 1$ with $\supp(f)\subset {\rm B}(z_0,\tau)$ for some $z_0\in\partial S$. Let $\chi$ be a smooth bump function such that $\chi=1$ on ${\rm B}(z_0,1)$ and $\chi=0$ outside ${\rm B}(z_0,2)$. Then there exist a constant $c$ depending only on $S$ and two functions $f_{+},f_{-}\in\mathscr{C}^{n-1}(\C)$ satisfying the following properties:
\begin{itemize}
\item[\textnormal{(i)}] $(f_{+}+f_{-})|_{S^{\rm c}}$ is the bounded harmonic extension of $f|_{\partial S}$ on $S^{\rm c}$.
\item[\textnormal{(ii)}] $\bar{\partial}f_{+}=\partial f_{-}=0$ on $S^{\rm c}$.
\item[\textnormal{(iii)}] $f_{+}\chi$ can be decomposed as $\sum_{k=1}^{\lfloor c\log N\rfloor}f_{+,k}$ such that $\supp(f_{+,k})\subset{\rm B}(z_0,\tau_k)$ for some $\tau_k\in[\tau,1]$ and $\|f_{+,k}\|_{n-1,\tau_k}\leq 1$. Similarly for $f_{-}$.
\item[\textnormal{(iv)}] $f_{+}(1-\chi)$ has uniformly bounded first and second derivatives on $\C$. Similarly for $f_{-}(1-\chi)$.
\end{itemize}
\end{lemma}

\subsection{Application of the Ward identity.}\
Let $\kappa>0$ and $g\in\mathscr{A}_{3,\tau}$ for some $\tau\in[N^{-1/2+\kappa},1]$. Let $\iota>0$ be a constant so that $\partial V-\partial\check{V}\neq0$ and $\Delta V>0$ on $S_{\rm edge}(2\iota)\setminus S$ (see e.g. \cite[Proof of Lemma 5.2]{AmeHedMak2015} for the existence of such $\iota$). Denote a bump function that is $1$ on $S\cup S_{\rm edge}(\iota)$ and $0$ on $S_{\rm out}(2\iota)$ by $\chi$. We

We apply \eqref{eqn:ward_main} for three different choices of $h$:
\begin{align}\label{eqn:h_1...}
h_1=\frac{\bar\partial g_{+}}{\Delta V/4},\quad h_2=\frac{\overline{\partial g_{-}}}{\Delta V/4},\quad h_3=\frac{\bar{\partial}g_0}{\Delta V/4}+r,\quad r=\Big(\frac{g_0}{\partial V-\partial \check{V}}-\frac{\bar{\partial}g_0}{\Delta V/4}\Big)\mathds{1}_{\bar{S}^{\rm c}}
\end{align}
where $g_0=(g-g_{\rm out})\chi$. As shown in \cite[Lemma 5.2]{AmeHedMak2011}, the function $h_3$ satisfies the conditions \eqref{eqn:h3}-\eqref{eqn:h1}; the continuity of $h_3$ relies on the fact that $g_0=0$ on $\partial S$ and $\Delta\check{V}=0$ on $S^{\rm c}$. The subtle point is that, although $r$ is continuous along $\partial S$, its derivatives are not. Nevertheless, the key estimates for $r$ remain valid as discussed below. 

The following lemma collects the decomposition of $h_1$, $h_2$, $h_3-r$ and $r$ together with the corresponding estimates. The part concerning $h_1$, $h_2$ and $h_3-r$ follows directly from Lemma \ref{lem:harm_ext_reg} (iii-iv), while the decomposition and estimates for $r$ are obtained by using Lemma \ref{lem:harm_ext_reg} (iii-iv) and following the argument in the proof of \cite[Lemma 5.2 (i)]{AmeHedMak2011}. We present it here for convenience, as it will be referred repeatedly, and omit the proof since it follows exactly as above.

\begin{lemma}\label{lem:dec}
Let $\kappa>0$, $n\geq 3$ and $g\in\mathscr{A}_{n,\tau}$ for some $\tau\in[N^{-1/2+\kappa},1]$. There exists a constant $C$ (depending only on $S$ and $n$) such that $g_{+}$ (and similarly $g_{-}$ and $g_{0}$) can be decomposed as $p+q$ where $p\in\mathscr{S}_{n-1,C,\kappa}$ and $q$ is supported in $S_{\rm out}(1/2)$ such that $\|\nabla q\|_{\infty}+\|\nabla^2 q\|_{\infty}<C$. 

Moreover, if $g$ is supported in $S$, then $r=0$. On the other hand, if there exists a $z_0\in\supp(g)\cap\partial S$, then $r$ can be decomposed similarly as $\sum_{k=1}^{\lfloor C\log N\rfloor}r_k+q$ where $q$ is again supported in $S_{\rm out}(1/2)$ such that $\|\nabla q\|_{\infty}+\|\nabla^2 q\|_{\infty}<C$; each $r_k$ is continuous, equal to $0$ on $\bar{S}$, with $\|\nabla^{j}r_k\|_{L^{\infty}(\bar{S}^{\rm c})}<\frac{C}{\tau_k^{j+1}}$ for $j=0,\dots,n-1$ and $\supp(r_k)\subset{\rm B}(z_0,\tau_k)$ for some $\tau_k\in[\tau,C]$.
\end{lemma}

We now proceed to apply the Ward identity to the functions $h_1$, $h_2$ and $h_3$. Recall that we denote $\lim_{|z|\to\infty}g_{\rm out}(z)$ by $g^{S}(\infty)$. Using $\bar{\partial}g_{+}=\partial g_{-}=0$ on $S^{\rm c}$, by Green's theorem we obtain:
\begin{align*} 
\int \frac{h_1(w)}{z-w}\rd\mu_{V}=\frac{1}{\pi}\int_{\C}\frac{\bar\partial g_{+}}{z-w}=g_{+}(z)-g^{S}(\infty),\quad \int \frac{h_2(w)}{z-w}\rd\mu_{V}=\frac{1}{\pi}\int_{\C}\frac{\overline{\partial g_{-}}}{z-w}=\overline{g_{-}(z)}
\end{align*}
for every $z\in\C\setminus\partial S$. On the other hand, because $g_0=0$ on the boundary
\begin{align*} 
\int \frac{h_3(w)}{z-w}\rd\mu_{V}=\frac{1}{\pi}\int_{S}\frac{\bar\partial g_{0}}{z-w}=\begin{cases}g_{0}(z), &z\in S\\ 0, & z\in \bar{S}^{\rm c}\end{cases},\quad \textnormal{and} \quad  h_3(z)(\partial V(z)-\partial \check{V}(z))=\begin{cases}0, &z\in S\\ g_{0}(z), & z\in \bar{S}^{\rm c}\end{cases}.
\end{align*}

Recall that $\tilde\mu_N=\frac{1}{N}\sum \delta_{z_k}-\mu_{V}$
Adding the three equations (after taking the conjugate of the second one) we obtain that
\begin{multline}
\E_f\Big[\sum_{k=1}^N g(z_k)\Big]-N\int g\,\rd\mu_V=\int \big(h_1\partial f +\overline{h_2 \partial f}+h_3\partial f\big) \rd\mu_{V} + \frac{1}{2}\int \big(\partial h_1+\overline{\partial h_2}+\partial{h_3}\big)\rd\mu_{V} + \mathscr{E}_N(f,g)
\\
+\E_f\left[\sum_{k=1}^N (g(z_k)-g_{\rm out}(z_k))(1-\chi(z_k))\right],\label{eqn:firstloop}
\end{multline}
where 
\begin{align*} 
\mathscr{E}_N(f,g)=&\ \E_{f}\Big[\int h_1\partial f +\overline{h_2\partial f}+h_3\partial f\,\rd\tilde{\mu}_{N}\Big]
\\
&+\E_{f}\Big[\int \partial h_1+\overline{\partial h_2}+\partial h_3 \,\rd\tilde{\mu}_{N}\Big]
\\
&+\E_{f}\Big[\frac{N}{2} \iint_{z\neq w} \frac{h_1(z)-h_1(w)}{z-w} \, \tilde \mu_N(\rd z) \, \tilde \mu_N(\rd w)\Big]+\frac{1}{2}\int \partial h_1\rd\mu_{V}
\\
&+\E_{f}\Big[\frac{N}{2} \iint_{z\neq w} \frac{\overline{h_2(z)}-\overline{h_2(w)}}{\overline{z-w}} \, \tilde \mu_N(\rd z) \, \tilde \mu_N(\rd w)\Big]+\frac{1}{2}\int \overline{\partial h_2}\rd\mu_{V}
\\
&+\E_{f}\Big[\frac{N}{2} \iint_{z\neq w} \frac{h_3(z)-h_3(w)}{z-w} \, \tilde \mu_N(\rd z) \, \tilde \mu_N(\rd w)\Big]+\frac{1}{2}\int \partial h_3\rd\mu_{V}.
\end{align*}
The expectation of the trace of $(g-g_{\rm out})(1-\chi)$ in the right-hand side of \eqref{eqn:firstloop} is already sub-polynomial, since its support lies outside the droplet and Lemma \ref{lem:1point}.(ii) gives exponential decay of the one-point function for the biased measure. We denote this contribution by $\OO(N^{-D})$ for arbitrarily large $D>0$. Furthermore, the main terms in equation \eqref{eqn:firstloop} simplify as follows. Since the expression on the left-hand side of the integral below depends only on the values of $f$ within $S$, we may, without loss of generality, replace $f$ by a smooth, compactly supported function that agrees with it in a neighbourhood of $S$. This modification is made solely for convenience in the forthcoming computations:
\begin{align*} 
\int \big(h_1\partial f +\overline{h_2 \partial f}+h_3\partial f\big) \rd\mu_{V}&=\frac{1}{\pi}\int_{S} (\bar\partial g_{+}\partial f+\partial g_{-}\bar\partial f+\bar\partial g_{0}\partial f) \rd m
\\
&=-\frac{1}{\pi}\int_{\C} (g_{+}+g_{-}) \bar{\partial}\partial f \rd m-\frac{1}{\pi}\int_{S} g_0 \bar{\partial}\partial f \rd m
\\
&=\frac{1}{4\pi}\int_{\C}\nabla (g_{+}+g_{-})\cdot\nabla f\rd m+\frac{1}{4\pi}\int_{S}\nabla g_0\cdot\nabla f\rd m
\\
&=\frac{1}{4\pi}\int_{S^{\rm c}}\nabla g^S\cdot\nabla f\rd m+\frac{1}{4\pi}\int_{S}\nabla g\cdot\nabla f\rd m
\\
&=\frac{1}{4\pi}\int_{\C}\nabla g^S\cdot\nabla f^S\rd m.
\end{align*}
Similarly, let $L=\log\frac{\Delta V}{4}$ and replace $L$ by a smooth, compactly supported function that agrees with it in a neighbourhood of $S$. The analogous calculations lead to,
{\allowdisplaybreaks
\begin{align*} 
\frac{1}{2}\int \big(\partial h_1+\overline{\partial h_2}+\partial{h_3}\big)\rd\mu_{V} &= 
\frac{1}{8\pi}\int_{S}\Delta g\,\rd m-\frac{1}{2\pi}\int_{S}\big(\bar{\partial}g_{+}\partial L+\partial g_{-}\bar\partial L+\bar{\partial}g_{0}\partial L\big)\, \rd m
\\
&=\frac{1}{8\pi}\int_{S}\Delta g\,\rd m-\frac{1}{2\pi}\int_{\C}\big(\bar{\partial}g_{+}\partial L+\partial g_{-}\bar\partial L\big)\,\rd m+\frac{1}{8\pi}\int_{S}g_0\Delta L\,\rd m
\\
&=\frac{1}{8\pi}\int_{S}\Delta g\,\rd m+\frac{1}{8\pi}\int_{\C}g_{\rm out}\Delta L\,\rd m+\frac{1}{8\pi}\int_{S}g_0\Delta L\,\rd m
\\
&=\frac{1}{8\pi}\int_{S}\Delta g\,\rd m+\frac{1}{8\pi}\int_{\C}g^{S}\Delta L\,\rd m
\end{align*}
}
and 
{\allowdisplaybreaks
\begin{align*} 
\frac{1}{8\pi}\int_{\C}g^{S}\Delta L\,\rd m &=\frac{1}{8\pi}\int_{S}g\Delta L\,\rd m+\frac{1}{8\pi}\int_{S^{\rm c}}\big(g^{S}\Delta L-L\Delta g^{S}\big)\,\rd m
\\
&=\frac{1}{8\pi}\int_{S}L\Delta g\,\rd m+\frac{1}{8\pi}\int_{\partial S}\Big(g\frac{\partial L}{\partial n}-L\frac{\partial g}{\partial n}\Big)+ \Big(g\frac{\partial L}{\partial(-n)}-L\frac{\partial g^{S}|_{S^{\rm c}}}{\partial(-n)}\Big)\,\rd s
\\
&=\frac{1}{8\pi}\int_{S}L\Delta g\,\rd m+\frac{1}{8\pi}\int_{\partial S} \Big(L^{S}\frac{\partial g}{\partial (-n)}-g\frac{\partial L^{S}|_{S^{\rm c}}}{\partial(-n)}\Big)\,\rd s 
\\
&=\frac{1}{8\pi}\int_{\C}L^{S}\Delta g\,\rd m
\end{align*}
}
where we used $0=\int_{S^{\rm c}}g^{S}\Delta L^{S}-L^{S}\Delta g^{S}=\int_{\partial S}g\frac{\partial L^{S}|_{S^{\rm c}}}{\partial(-n)}-L\frac{\partial g^{S}|_{S^{\rm c}}}{\partial(-n)}$ in the third equality.

Substituting these into \eqref{eqn:firstloop}, we obtain 
\begin{equation}\label{eqn:finalloop}
\E_f\Big[\sum_{k=1}^N g(z_k)\Big]-N\int g\,\rd\mu_V=\frac{1}{4\pi}\int_{\C}\nabla g^S\cdot\nabla f^S\rd m + \frac{1}{8\pi}\int_{S}\Delta g\,\rd m+\frac{1}{8\pi}\int_{\C}L^{S}\Delta g \,\rd m + \mathscr{E}_N(f,g)+\OO(N^{-D}).
\end{equation}
Note that every term in this equation is linear in $g$, but not in $f$. Moreover, an alternative representation involving a Neumann jump term follows from Green’s identity: 
\begin{multline} \label{eqn:green_to_neumann}
\int_{\C}L^{S}\Delta g \,\rd m=\int_{S}\Delta g L\,\rd m+\int_{S^{\rm c}}\big(\Delta g L^{S}-g\Delta L^{S}\big)\rd m
\\
=\int_{S} g\Delta L\,\rd m+\int_{\partial S}\Big(L\frac{\partial g}{\partial n}-g\frac{\partial L}{\partial n}\Big)\rd s+\int_{\partial S}\Big(L\frac{\partial g}{\partial (-n)}-g\frac{\partial L^{S}|_{S^{\rm c}}}{\partial (-n)}\Big)\rd s
=\int_{S} g\Delta L\,\rd m-\int_{\partial S}g\mathcal{N}(L)\,\rd s
\end{multline}
which holds given $g$ smooth, compactly supported.

Next, we apply the estimation \eqref{eqn:finalloop} to the Laplace transform of the linear statistics.

\begin{proposition}\label{prop:Laplace_loop}
Let $C,\kappa>0$ be fixed constants. Let $g=g_1+g_2$ where $g_1\in\mathscr{S}_{3,C,\kappa}$, $g_2\in\mathscr{C}^{3}(\C)$ supported in $S_{\rm out}(1/2)$ with at most logarithmic growth as $|z|\to\infty$ and $\|\nabla g_2\|_{\infty}=\OO(N^{\log\log N})$. Let $f$ be of class $\mathscr{C}^{5}(\C)$ with at most logarithmic growth as $|z|\to\infty$ and $|\nabla f(z)|\leq c(1+|z|)^c$ for an $N$-dependent constant $c$. Then, for any $D>0$,
\begin{align*} 
\log\E [e^{\sum_{i}(f+g)(z_i)}]=&\,\log\E [e^{\sum_{i}f(z_i)}]+\frac{1}{4\pi}\int_{\C}\nabla g^S\cdot\nabla f^S\rd m +\frac{1}{8\pi}\int_{\C}|\nabla g^S|^2\rd m + \frac{1}{8\pi}\int_{S}\Delta g\,\rd m
\\
&+\frac{1}{8\pi}\int_{S} g\Delta L\,\rd m-\frac{1}{8\pi}\int_{\partial S}g\mathcal{N}(L)\,\rd s
+N\int g\,\rd\mu_{V}
+\int_{0}^{1}\mathscr{E}_N(f+t g,g_1)\rd t
+\OO(N^{-D}).
\end{align*}
The error term is uniform in $g_1$, $g_2$ and $f$.
\end{proposition}

Observe that if $g_2$ is compactly supported $\int_{S} g\Delta L\,\rd m-\int_{\partial S}g\mathcal{N}(L)\,\rd s$ can be replaced by $\int_{\C}L^{S}\Delta g \,\rd m$.

\begin{proof}
We define $F_{N}(h)=\log \mathbb{E}\left[e^{X_{f+h}}\right],\ X_h=\sum_{i}h(z_i)-N\int h\,\rd \mu_{V}$ for any function $h$. Equation \eqref{eqn:finalloop} together with \eqref{eqn:green_to_neumann} applied to $f+tg$ and $g$, substituted into $f$ and $g$ respectively, gives
\begin{align*}
\frac{\rd}{\rd t} F_N(t g)=\E_{f+tg}\left[X_{g}\right]
=&\,
\frac{1}{4\pi}\int_{\C}\nabla g_1^S\cdot\nabla f^S\rd m+\frac{t}{4\pi}\int_{\C}|\nabla g_1^S|^2\rd m + \frac{1}{8\pi}\int_{S}\Delta g_1\,\rd m
\nonumber
\\
&+\frac{1}{8\pi}\int_{S} g_1\Delta L\,\rd m-\frac{1}{8\pi}\int_{\partial S}g_1\mathcal{N}(L)\,\rd s+\E_{f+tg}\left[X_{g_2}\right]+\mathscr{E}_N(f+tg,g_1)+\OO(N^{-D}).
\end{align*}
In the integrals on the right-hand side, each $g_1$ may be replaced by $g$, as the support of $g-g_1=g_2$ lies entirely outside of a neighbourhood of the set $\bar{S}$. Thus, we only need to prove that
\begin{align*} 
\E_{f+tg}\left[X_{g_2}\right]=\OO(N^{-D}).
\end{align*}
This follows easily by the exponential decay of the one-point correlation function outside the droplet $S$, i.e., Lemma \ref{lem:1point}.(ii).
\end{proof}

\subsection{Fluctuations.}\label{subsec:flucs}\ 
The following central limit theorem is a quantitative extension to mesoscopic scales of the main result from \cite{AmeHedMak2015}.
The proof follows the same path, with (i) Johansson's idea to reduce Laplace transform estimates to a shift of the density of states \cite{Joh1998}, (ii) the Ward identity, (iii) estimates of the kernel of the determinantal point process as an input. 

The main differences with \cite{AmeHedMak2015} are the following: for the error estimations in step (ii) we need an isotropy result for general potential, i.e. Lemma \ref{lem:anis}, and for step (iii) we need kernel asymptotics when the external potential is perturbed on a mesoscopic scale, a technicality proved in Appendix \ref{sec:kernel} similarly to the macroscopic case \cite{Ber2009, AmeHedMak2015}.

\begin{proposition}\label{proposition.nosingularities}
Let $0<\kappa\leq 1/2$, $C>0$ be fixed constants. 
Then, uniformly in $f\in\mathscr{S}_{5,C,\kappa}$, we have
\begin{align*}
\E\left[e^{\sum_{i=1}^Nf(z_i)}\right]=
e^{N\int f\,\rd\mu_{V}
+\frac{1}{8\pi}\int_{\C}\big(|\nabla f^S|^2
+\mathds{1}_{S}\Delta f
+L^{S}\Delta f\big)\rd m
}
\left(1+\OO(N^{-\kappa/10})\right).
\end{align*}
\end{proposition}

For the proof of this proposition, we apply Proposition \ref{prop:Laplace_loop} for $g$ and $f$  substituted by $f$ and $0$. This yields the desired result  provided that
\begin{equation}\label{eqn:error}
\mathscr{E}_N(tf,f)=\OO(N^{-\kappa/10})
\end{equation}
uniformly in $0\leq t\leq 1$. We now prove the above estimate, only for $t=1$ to simplify notations  (uniformity in $t$ of the following estimates is easy to check). Note that by linearity of $\mathscr{E}_N$ in the second coordinate, it suffices to prove the error estimation \eqref{eqn:error} for $\mathscr{E}_N(f,g)$ uniformly in $g$ or type $\mathscr{A}_{5,\tau}$ for a $\tau\in[N^{-1/2+\kappa},1]$.

We now prove $\mathscr{E}_N(f,g)=\OO(N^{-\kappa/10})$ holds in the following sequence of Lemmas \ref{lem1} to \ref{lem3}, which conclude the proof. We only state the lemmas for $t=1$, naturally the proof extends to any  $0\leq t\leq 1$ and the error estimates are uniform in this parameter range.

\begin{lemma}\label{lem1}
Let $C,\kappa>0$, $f\in\mathscr{S}_{3,C,\kappa}$ and $g\in\mathscr{A}_{3,\tau}$ for some $\tau\in[N^{-1/2+\kappa}, 1]$. Let $h_i$ be as defined in \eqref{eqn:h_1...}. For any $\epsilon>0$ and $i=1,2,3$ we have
\begin{align}
&\E_{f}\Big[\int h_i\partial f \,\rd\tilde{\mu}_{N}\Big]=\OO(N^{-\kappa+\epsilon}) \label{eqn:center1}
\end{align}
where the implicit constant is independent of $f$, $\tau$ and $g$.
\end{lemma}
\begin{proof} Recall the definition of function $r$ from \eqref{eqn:h_1...}. By Lemma \ref{lem:dec} $h_1$, $h_2$ and $(h_3-r)$ can be decomposed as $\frac{\bar{\partial} p}{\Delta V}+q$ where $p\in\mathscr{S}_{2,C,\kappa}$ and $q$ is supported in $S_{\rm out}(1/2)$ and uniformly Lipschitz with Lipschitz constant depending only on $S$. By Lemma \ref{lem:1point}.(ii), the contribution of $q$ is easily subpolynomial. Thus, it suffices to prove the estimation for $\frac{\bar{\partial} p}{\Delta V}$ when $p\in\mathscr{A}_{\tau}$ for some $\tau\in[N^{-1/2+\kappa},1]$.
\begin{align*}
\E_f\Big[\int\frac{\bar{\partial} p}{\Delta V}\partial f\rd\tilde\mu_N\Big]=\int \frac{\bar{\partial} p}{\Delta V}\partial f \left( N^{-1}\wt {\bf K}_N(z,z)-\rho_{V}(z)\right)\rd m(z),
\end{align*}
where we used the notation from Theorem \ref{prop:kernel}. From this theorem, for $z\in S_{\rm in}(\delta_N)$, we have $|(N^{-1}\wt {\bf K}_N(z,z)-\rho_{V}(z)|\leq N^{-1}\sup_{{\rm B}(z,2 \delta_N)}|\nabla^2 f|\leq N^{-2\kappa}\log N
$, because of $f\in\mathscr{S}_{3,C,\kappa}$. We therefore have
\begin{align}\label{eqn:needchange}
\left|\E_f\Big[\int_{S_{\rm in}(\delta_N)}\frac{\bar{\partial} p}{\Delta V}\partial f\rd\tilde\mu_N\Big]\right|\leq N^{-2\kappa+\epsilon}\int_{S}|\frac{\bar{\partial} p}{\Delta V}\partial f|\,\rd m\leq N^{-2\kappa+\epsilon}\int_{S}\big(\frac{|\nabla p|^2}{|\Delta V|^2}+|\nabla f|^2\big)\rd m\leq N^{-2\kappa+\epsilon}
\end{align}
where we used the easy estimate $\int_{S}|\nabla f|^2\leq \pi C\log N$ for any $f\in\mathscr{S}_{3,C,\kappa}$.
The contribution from $z\in S_{\rm out}(\delta_N)$ can be bounded easily with Lemma \ref{lem:1point}.(ii), which gives 
\begin{equation*}
\wt{\bf K}_N(z,z)\leq \exp(-c(\log N)^2)
\end{equation*}
for all such $z$.
Finally, for the domain $z\in S_{\rm edge}(\delta_N)$, note that by Lemma \ref{lem:1point}.(i) there exists a universal $C>0$ such that 
\begin{equation*}
\frac{\wt {\bf K}_N(z,z)}{N}\leq C.
\end{equation*}
This gives, using a shorthand notation $\wt {\bf K}_N$ for the function $\wt {\bf K}_N(z,z)$,
\begin{equation}\label{eqn:ring}
\int_{S_{\rm edge}(\delta_N)} |\frac{\bar{\partial} p}{\Delta V}\partial f| \left|\frac{\wt {\bf K}_N}{N}-\rho_{V}\right|\rd m\lesssim \int_{S_{\rm edge}(\delta_N)} |\frac{\bar{\partial} p}{\Delta V}\partial f|\,\rd m\leq N^{-\kappa+\epsilon}
\end{equation}
because $\frac{\bar{\partial} p}{\Delta V}$ is $\OO(1/\tau)$ on a domain of area $\OO(\tau\delta_N)$, and $\partial f$ is $\OO(N^{1/2-\kappa}\log N)$. 

Thus, it only remains to show $\E_{f}[\int r\partial f \,\rd\tilde{\mu}_{N}]=\OO(N^{-\kappa+\epsilon})$ the following to complete the proof of \eqref{eqn:center1}. This follows the same argument: On $S_{\rm out}(\delta_N)$, it's again subpolynomial and at the edge, the calculation is identical to \eqref{eqn:ring}.
\end{proof}

\begin{lemma}\label{lem2} Let $C,\kappa>0$, $f\in\mathscr{S}_{3,C,\kappa}$ and $g\in\mathscr{A}_{3,\tau}$ for some $\tau\in[N^{-1/2+\kappa}, 1]$. Let $h_i$ be as defined in \eqref{eqn:h_1...}. For any $\epsilon>0$ we have
\begin{align*} 
\E_{f}\Big[\int \partial h_1+\overline{\partial h_2}+\partial h_3 \,\rd\tilde{\mu}_{N}\Big]=\OO(N^{-\kappa+\epsilon})
\end{align*}
where the implicit constant is independent of $f$, $\tau$ and $g$.
\end{lemma}
\begin{proof} Recall the definition of function $r$ from \eqref{eqn:h_1...}. We have
\begin{align*} 
\partial h_1+\overline{\partial h_2}+\partial h_3=\frac{\Delta g}{\Delta V}-\Big(\frac{\bar{\partial}g_{+}\partial L}{\Delta V/4}+\frac{\partial g_{-}\bar{\partial} L}{\Delta V/4}+\frac{\bar{\partial}g_{0}\partial L}{\Delta V/4}\Big)+\partial r
\end{align*}
on $\C\setminus \partial S$ where $L=\log\frac{\Delta V}{4}$. We prove that each term in this sum yield to a negligible integral. All terms except $\partial r$ can be treated in the same way as in the previous lemma. Starting with the first term, we write
\begin{align}\label{eqn:changeinlem2} 
\E_{f}\big[\int \frac{\Delta g}{\Delta V} \,\rd\tilde{\mu}_{N}\big]=\int \frac{\Delta g}{\Delta V} \left( N^{-1}\wt {\bf K}_N-\rho_{V}\right)\rd m.
\end{align}
On $z\in S_{\rm in}(\delta_N)$, the integral is $\OO(N^{-2\kappa+\epsilon})$ by Theorem \ref{prop:kernel}. On $z\in S_{\rm out}(\delta_N)$, the integral is subpolynomial due to Lemma \ref{lem:1point}.(ii). On $S_{\rm edge}(\delta_N)$, the integrand is $\OO(1/\tau^2)$ and the integration area is $\tau\delta_N$ which gives $\OO(N^{-\kappa+\epsilon})$ bound. 

By Lemma \ref{lem:dec}, each $g_{+},g_{-},g_{0}$ can be written as $p+q$ where $p\in\mathscr{S}_{2,C,\kappa}$ and Lipschitz $q$ with $\supp(q)\subset S_{\rm out}(1/2)$. The contribution of $q$ is easily subpolynomial. Thus, it suffices to consider the following expectation for $p\in\mathscr{A}_{2,\tau}$ with $\tau\in[N^{-1/2+\kappa},1]$:
\begin{align}\label{eqn:change2inlem2} 
\E_{f}\big[\int \frac{\bar{\partial}p\partial L}{\Delta V/4} \,\rd\tilde{\mu}_{N}\big]=\int \frac{\bar{\partial}p\partial L}{\Delta V/4} \left( N^{-1}\wt {\bf K}_N-\rho_{V}\right)\rd m(z).
\end{align}
By the same way, it follows that this term is $\OO(N^{-\kappa})$.

Lastly, $r$ admits a similar decomposition by Lemma \ref{lem:dec}. For the term $\E_{f}[\int \partial r\,\rd\tilde\mu_N]$, the integral over $S_{\rm out}(\delta_N)$ is already subpolynomial, while the contribution from $S_{\rm edge}(\delta_N)$ can be bounded in the same way by $\OO(\frac{\delta_{N}\log N}{\tau})$.
\end{proof}

\begin{lemma}\label{lem3} Let $C,\kappa>0$, $f\in\mathscr{S}_{3,C,\kappa}$ and $g\in\mathscr{A}_{5,\tau}$ for some $\tau\in[N^{-1/2+\kappa}, 1]$. Let $h_i$ be as defined in \eqref{eqn:h_1...}. For any $\epsilon>0$ and $i=1,2,3$ we have
\begin{align}\label{eqn:lem3}
\E_{f}\Big[\frac{N}{2} \iint_{z\neq w} \frac{h_i(z)-h_i(w)}{z-w} \, \tilde \mu_N(\rd z) \, \tilde \mu_N(\rd w)\Big]+\frac{1}{2}\int \partial h_i\,\rd\mu_{V}=\OO(N^{-\kappa/10})
\end{align}
where the implicit constant is independent of $f$, $\tau$ and $g$.
\end{lemma}
\begin{proof} We begin with proving that $\E_{f}[N \iint_{z\neq w} \frac{r(z)-r(w)}{z-w}  \rd\tilde \mu_N  \rd\tilde \mu_N]=I(r;\C,\C)$ is negligible where we define
\begin{align*} 
I(r;A,B)=N\int_{A}\int_{B} \frac{r(z)-r(w)}{z-w}\left(\Big(\frac{\wt{\bf{K}}_N(z,z)}{N}-\rho_V(z)\Big)\Big(\frac{\wt{\bf{K}}_N(w.w)}{N}-\rho_V(w)\Big)-\frac{|\wt{\bf{K}}_N(z,w)|^2}{N^2}\right)\rd m(w)\rd m(z).
\end{align*}
 Recall the $r=q+\sum r_k$ decomposition from Lemma \ref{lem:dec}. By Lemma \ref{lem:1point}.(ii) $I(q;\C,\C)$ is subpolynomial. For $I(r_k;\C,\C)$, again by Lemma \ref{lem:1point}.(ii) when either of $z$ and $w$ is in $S_{\rm out}(\delta_N)$, the contribution to the integral is subpolynomial. Then there are two remaining cases: $z,w\in S_{\rm edge}(\delta_N)$ and $z\in S_{\rm edge}(\delta_N)$, $w\in S_{\rm in}(\delta_N)$. Recall that by Lemma \ref{lem:dec}, $r_k$ is a continuous function that is $0$ on $\bar{S}$, supported on a ball with radius $\tau_k$ and $\sup_{\C\setminus\partial S}|\nabla r|\leq \frac{1}{\tau_k}$; which implies $|r_k|\leq\frac{\delta_N}{\tau_k^2}$ on $S_{\rm edge}(\delta_N)$.
 
 If $z,w\in S_{\rm edge}(\delta_N)$, we rely only on the bound $\wt{\bf{K}}_N(z,z)=\OO(N)$ from Lemma \ref{lem:1point}.(i). Therefore,
\begin{align*} 
|I(r_k;S_{\rm edge}(\delta_N),S_{\rm edge}(\delta_N))|=N \iint_{S_{\rm edge}(\delta_N)^2}\frac{|r_k(z)|}{|z-w|} \OO(1)\rd m(z)\rd m(w)=\OO(N\delta_N\log N)\int_{S_{\rm edge}(\delta_N)}|r_k(z)|\rd m(z)
\\
=\OO(\frac{N\delta_N^3\log N}{\tau_k})=\OO(N^{-\kappa+\varepsilon}).
\end{align*}

If $z\in S_{\rm in}(\delta_N),w\in S_{\rm edge}(\delta_N)$, using Theorem \ref{prop:kernel} for $z$, we obtain
\begin{align*}
|I(r_k;S_{\rm in}(\delta_N), S_{\rm edge}(\delta_N))|= \OO(N) \iint_{\substack{z\in S_{\rm in}(\delta_N),\\w\in S_{\rm edge}(\delta_N)}}\frac{|r_k(w)|}{|z-w|}\OO\left(\frac{\sup_{{\rm B}(z,2\delta_N)}|\nabla^2 f|}{N}+N^{-3/2}+\frac{|\wt{\bf{K}}_N(z,w)|^2}{N^2}\right)\rd m(z)\rd m(w).
\end{align*}
The contribution from $\frac{|\wt{\bf{K}}_N(z,w)|^2}{N^2}$  can be bounded by the reproducing property \eqref{eqn:reprod} as follows: 
\begin{align*} 
N \iint_{\substack{z\in S_{\rm in}(\delta_N),\\w\in S_{\rm edge}(\delta_N)}}\frac{|r_k(w)|}{|z-w|}\frac{|\wt{\bf{K}}_N(z,w)|^2}{N^2}\rd m(z)\rd m(w)=\OO(\frac{1}{N\tau_k^2})\iint_{\substack{z\in\C,\\w\in S_{\rm edge}(\delta_N)\cap\supp(r)}}|\wt{\bf{K}}_N(z,w)|^2\rd m(z)\rd m(w)
\\
=\OO(\frac{1}{N\tau_k^2})\int_{S_{\rm edge}(\delta_N)\cap\supp(r)}\wt{\bf{K}}_N(w,w)\,\rd m(w)=\OO(\frac{\delta_N}{\tau_k})=\OO(N^{-\kappa+\epsilon}).
\end{align*}

Moreover, as $f\in\mathscr{S}_{3,C,\kappa}$, it can be written as $\sum_{\ell<C\log N}f_\ell$ where $f_\ell$ is supported in a ball centered at $z_{\ell}'$ with radius $\tau_\ell'$ and smooth on that scale for some $\tau_\ell'\in[N^{-1/2+\kappa},1]$, i.e. $\|\nabla^{3}f_\ell\|\leq (\tau_\ell')^{-3}$. Note that for each $\ell$ we have $\sup_{{\rm B}(z,2\delta_N)}|\nabla^2 f_{\ell}|=\OO(\frac{1}{(\tau_{\ell}')^2}\mathds{1}_{{\rm B}(z_{\ell}',2\tau_{\ell}')}(z))$.  So, the contribution from $\frac{\sup_{{\rm B}(z,2\delta_N)}|\nabla^2 f|}{N}$ will be bounded by
\begin{multline}\label{eqn:another_change} 
\sum_{\ell}\frac{1}{(\tau_{\ell}')^2} \iint_{\substack{z\in S_{\rm in}(\delta_N)\cap {\rm B}(z_{\ell}',2\tau_{\ell}'),\\ w\in S_{\rm edge}(\delta_N)}}\frac{|r_k(w)|}{|z-w|}\rd m(z)\rd m(w)
\\
=\sum_{\ell}\OO(\frac{\delta_N\log N}{\tau^2(\tau_{\ell}')^2}) \iint_{\substack{z\in S_{\rm in}(\delta_N)\cap {\rm B}(z_{\ell}',2\tau_{\ell}'),\\ w\in S_{\rm edge}(\delta_N)}}\frac{1}{|z-w|}\rd m(z)\rd m(w)
=O(\frac{\delta_N^2(\log N)^2}{\tau^2})=\OO(N^{-\kappa}).
\end{multline}
Lastly, the contribution from $N^{-3/2}$ is easily seen to be $\OO(N^{-1/2}\delta_N/\tau^2)=\OO(N^{-\kappa})$. This completes the proof of negligiblity of $\E_{f}[N \iint_{z\neq w} \frac{r(z)-r(w)}{z-w}  \rd\tilde \mu_N  \rd\tilde \mu_N]$. \\

We now proceed to proving \eqref{eqn:lem3} for $h_1$, $h_2$ and $(h_3-r)$. By Lemma \ref{lem:dec}, it suffices to consider $\E_f[\frac{N}{2}\iint \frac{h(z) - h(w)}{z-w}\rd\tilde\mu_N\rd\tilde\mu_N]$ when $\tau\cdot h\in \mathscr{A}_{3,\tau}$ and $\supp(h)\subset B(z_0,\tau)$ for some $\tau\in[N^{-1/2+\kappa},1]$, $z_0\in \{z\in\C:\dist(z,S)\leq 2\}$.

We rely on a multiscale decomposition idea from \cite{FefdeL1986} which was used for two dimensional Coulomb gases in \cite{BauBouNikYau2017,BauBouNikYau2019}. 
	Let $\varphi(z)=e^{-|z|^2}$. It's easy to verify that
	\begin{gather}
\frac{2}{\pi}\int_{\C} \varphi\left(\frac{z-\xi}{t}\right) \varphi\left(\frac{w-\xi}{t}\right)\rd m(\xi)=t^2e^{-\frac{|z-w|^2}{2 t^2}}, \label{eqn:exp_decomp1}
\\
 \int_\eta^\infty   \left(\frac{2}{\pi}\int_{\C} \varphi\left(\frac{z-\xi}{t}\right) \varphi\left(\frac{w-\xi}{t}\right)\rd m(\xi)\right) \frac{\rd t }{t^5}=\frac{1-e^{-\frac{|z-w|^2}{2\eta^2}}}{|z-w|^2}. \label{eqn:exp_decomp2}
\end{gather}  
Let $\eta<\tau$ be some small fixed parameter to be chosen. From \eqref{eqn:exp_decomp1}-\eqref{eqn:exp_decomp2}, for any function $h$ we have
	\begin{align*}
\E_f\left[\frac{N}{2}\iint \frac{h(z) - h(w)}{z-w}\rd\tilde\mu_N(z)\rd\tilde\mu_N(w)\right]=E_1(\eta,h)+E_2(\eta,h)
	\end{align*}
where
\begin{align} 
&E_1(\eta,h)=\frac{N}{2}\E_f\left[\iint_{z\neq w}\frac{h(z)-h(w)}{z-w}e^{-\frac{|z-w|^2}{2\eta^2}}\rd\tilde\mu_N(z)\rd\tilde\mu_N(w)\right], \label{eqn:E_1}
	\\
	&E_2(\eta,h)=\frac{N}{\pi}\int_\eta^{\infty}\int_{\mathbb{C}}\E_f\left[\iint_{z\neq w} \varphi\left(\frac{z-\xi}{t}\right)\varphi\left(\frac{w-\xi}{t}\right)(h(z)-h(w))(\overline{z-w})\rd\tilde\mu_N(z)\rd\tilde\mu_N(w)\right]m(\rd \xi)\frac{\rd t}{t^5}.\nonumber 
\end{align}
To bound $E_2$ we use the measure cancellation due to $\tilde{\mu}_N$, i.e., the explicit uniform bounds in Corollary \ref{prop:Laplace_S_kappa}. 
\begin{align*} 
E_2(\eta,h)=&\frac{2N}{\pi}\sum_{i=1}^{2}\int_\eta^{\infty}\int_{\mathbb{C}}\E_f\left[\int u_i(z)\rd\tilde{\mu}_N(z)\int v_i(w)\rd\tilde{\mu}_N(w)\right]m(\rd \xi)\frac{\rd t}{t^5}
\end{align*}	
where 
\begin{gather*} 
u_1(z)=\varphi\left(\frac{z-\xi}{t}\right)(h(z)-h(\xi))(\overline{z-\xi}),\quad  v_1(w)=\varphi\left(\frac{w-\xi}{t}\right),
\\
 u_2(z)=\varphi\left(\frac{z-\xi}{t}\right)(h(z)-h(\xi)),\quad  v_2(w)=\varphi\left(\frac{w-\xi}{t}\right)(\overline{\xi-w}).
\end{gather*}
With a negligible error, we can assume that the function $\varphi$ is truncated so that it's supported in ${\rm B}(0,(\log N)^{\log\log N})$. Thus, the integration region for $\xi$ can be restricted to ${\rm B}(z_0,\tau+t(\log N)^{\log\log N})$. Moreover, note that the contribution from $t>N^{10}$ is negligible, because $|\int u_1\rd\tilde{\mu}_N|\lesssim \|u_1\|_{\infty}\lesssim t$ and similarly $|\int u_2\rd\tilde{\mu}_N|\lesssim 1$, $|\int v_1\rd\tilde{\mu}_N| \lesssim 1$, $|\int v_2\rd\tilde{\mu}_N| \lesssim t$. This leads to a negligible contribution of order $N\int_{N^{10}}^{\infty} t^{-2}\rd t$. 

On the other hand, by Corollary \ref{prop:Laplace_S_kappa} we obtain
\begin{align*} 
|\int u_1\rd \tilde{\mu}_N|\leq N^{-1+\epsilon}\frac{t^2}{\tau^2}, \ |\int v_1\rd \tilde{\mu}_N| \leq N^{-1+\epsilon}, \ |\int u_2\rd \tilde{\mu}_N| \leq N^{-1+\epsilon}\frac{t}{\tau^2}, \ |\int v_2\rd \tilde{\mu}_N| \leq N^{-1+\epsilon}t
\end{align*}
with probability greater than $1-e^{-N^{\epsilon/2}}$ with respect to the unbiased measure for every $t$ and $\xi$ ($\epsilon$ is an arbitrarily small positive constant, that may change line-by-line). Uniformity of these inequalities in $t\in[\eta,N^{10}]$ and $\xi\in {\rm B}(z_0,N^{10})$ can be argued easily by a grid argument.
Moreover, by Corollary \ref{cor:stability}, the same inequalities hold with overwhelming probability with respect to the biased measure. Thus, the contribution to $E_2$ becomes:
\begin{align*} 
N^{-1+\epsilon}\int_{\eta}^{N^{10}}(t\vee\tau)^2\frac{t^2}{\tau^2}\rd m(\xi) \frac{\rd t}{t^5}\leq N^{-2\kappa+\epsilon}+\frac{N^{-1+\epsilon}}{\eta^2}.
\end{align*}
We now choose $\epsilon<\kappa/100$ and $\eta=N^{-1/2+\kappa/15}$.

 To bound $E_1$, we do not only use the measure cancellation due to $\tilde{\mu}_N$ and rely also on some phase cancellation of the kernel first, which follows from Theorem \ref{prop:kernel}. We write the Taylor expansion for $h$ first.	
	\begin{equation*}
  \frac{h(z)-h(w)}{z-w} = \partial h(z) + \bar\partial h(z) \frac{\bar z - \bar w}{z-w} + \OO(|z-w|).
\end{equation*}
The negligibility of the contribution from $\OO(|z-w|)$ to $E_1$ follows by
\begin{align} 
N\,\E_f\bigg[\iint& |z-w| e^{-\frac{|z-w|^2}{2\eta^2}}\big(\rd\hat{\mu}_N(z)+\rho_{V}(z)\,\rd m(z)\big)\big(\rd\hat{\mu}_N(w)+\rho_{V}(w)\,\rd m(w)\big)\bigg]\nonumber
\\
\lesssim &\, N\iint |z-w|e^{-\frac{|z-w|^2}{2\eta^2}}\det\begin{vmatrix} \wt {\bf K}_N(z,z)/N & \wt {\bf K}_N(z,w)/N \\ \wt {\bf K}_N(w,z)/N & \wt {\bf K}_N(w,w)/N\end{vmatrix}\rd m(z)\rd m(w)\nonumber
\\
&+N\int_{S}\int |z-w|e^{-\frac{|z-w|^2}{2\eta^2}}\frac{\wt {\bf K}_N(z,z)}{N}\rd m(z)\rd m(w)+N\int_{S}\int_{S} |z-w|e^{-\frac{|z-w|^2}{2\eta^2}}\rd m(z)\rd m(w).\label{eqn:E_1case1}
\end{align}
Using $\wt{\bf K}_N(z,z)/N\leq C$ by Lemma \ref{lem:1point}.(i) and the decay of $\wt{\bf K}_N(z,z)$ as $|z|\gg 1$ by Lemma \ref{lem:1point}.(ii), we obtain that the expression above is bounded by $N\eta^3$. 

On the other hand, the contribution from $\partial h(z)$ to $E_1$ can be bounded similarly to $E_2$ using \eqref{eqn:exp_decomp1}, except the fact that we get a deterministic extra term $\frac{1}{2}\int\partial h\rd\mu_{V}$ due to the diagonal terms. More explicitly, substituting \eqref{eqn:exp_decomp1} into $\frac{N}{2}\E_f[\iint_{z\neq w}\partial h(z)e^{-\frac{|z-w|^2}{2\eta^2}}\rd\tilde\mu_N\rd\tilde\mu_N]$ we get that it's equal to
\begin{align} \label{eqn:E_1case2}
\frac{N}{\pi\eta^2}\int_{\C}\E_f\left[\int\partial h(z) \varphi\left(\frac{z-\xi}{\eta}\right) \rd\tilde\mu_N(z)\int \varphi\left(\frac{w-\xi}{\eta}\right)\rd\tilde\mu_N(w)\right]\rd m(\xi)-\frac{1}{2N}\mathbb{E}_{f}\Big[\sum_{i=1}^{N}\partial h(z_i)\Big].
\end{align}
The first term can be bounded similarly with $E_2(\eta,h)$, by $\frac{N}{\eta^2}\tau^2\frac{N^{-1+\epsilon}}{\tau^2}N^{-1+\epsilon}=N^{-2\kappa/15+2\epsilon}$ where $\tau^2$ factor stands for the area of the integration region for $\xi$, $\frac{N^{-1+\epsilon}}{\tau^2}$ and $N^{-1+\epsilon}$ factors stand for the normalized centered linear statistics terms respectively. The second term can be evaluated easily as $-\frac{1}{2}\int\partial h\rd\mu_{V}$ up to a negligible error, by Theorem \ref{prop:kernel} and Lemma \ref{lem:1point}:
\begin{multline*} 
\frac{1}{N}\mathbb{E}_{f}\big[\sum_{i=1}^{N}\partial h(z_i)\big]-\int\partial h\,\rd\mu_{V}=\int\partial h(z)(\frac{\wt{\bf K}_N}{N}-\rho_{V})\rd m(z)
\\
= \oo(N^{-1})+\int_{S_{\rm in}(\delta_N)\cap\supp(h)}\frac{1}{\tau^2}\OO(N^{-2\kappa+\epsilon})\rd m+\int_{S_{\rm edge}(\delta_N)\cap\supp(h)}\frac{1}{\tau^2}\OO(1)\rd m=\OO(N^{-\kappa+\epsilon}).
\end{multline*}

Lastly, we finish the proof by discussing the contribution of the term $\bar\partial h(z) \frac{\bar z - \bar w}{z-w}$ to $E_1(\eta,h)$. By the isotropy lemma, i.e. Lemma \ref{lem:anis}, the combination from $\rd\hat\mu_N(z)\rd\hat\mu_N(w)$ can be bounded easily:
\begin{align} \label{eqn:anis_app}
\frac{N}{2}\E_f\left[\iint_{z\neq w}\bar\partial h(z)\frac{\bar z-\bar w}{z-w}e^{-\frac{|z-w|^2}{2\eta^2}}\rd\hat\mu_N(z)\rd\hat\mu_N(w)\right]=\OO(N^{-\kappa})
\end{align}
as $\|\tau^2\partial h\|_{\infty}\lesssim 1$. Secondly, the $\rd \tilde{\mu}_N(z)\rd\mu_{V}(w)$ integral can be bounded as follows,
\begin{multline} 
\frac{N}{2}\E_f\left[\int\bar\partial h(z)\int_{S} \frac{\bar z - \bar w}{z-w}e^{-\frac{|z-w|^2}{2\eta^2}}\rd\mu_{V}(w)\rd\tilde\mu_N(z)\right]
\\
=\frac{N}{2}\E_f\left[\int_{S_{\rm in}(\eta\log N)}\bar\partial h(z)\int_{{\rm B}(z,\eta\log N)} \frac{\bar z - \bar w}{z-w}e^{-\frac{|z-w|^2}{2\eta^2}}(\rho_{V}(w)-\rho_{V}(z))\rd m(w)\rd\tilde\mu_N(z)\right]
\\
+\frac{N}{2}\E_f\left[\int_{S_{\rm in}(\eta\log N)}\bar\partial h(z)\int_{S\setminus {\rm B}(z,\eta\log N)} \frac{\bar z - \bar w}{z-w}e^{-\frac{|z-w|^2}{2\eta^2}}\rd\mu_{V}(w)\rd\tilde\mu_N(z)\right]
\\
+\frac{N}{2}\E_f\left[\int_{S_{\rm in}(\eta\log N)^{\rm c}}|\bar\partial h(z)|\eta^2\big(\rd\hat\mu_N(z)+\rd\mu_{V}(z)\big)\right]=\OO(\frac{N\eta^3}{\tau})\ll N^{-\kappa/2}\label{eqn:E_1case3_2}
\end{multline}
which can be explained as follows. The first expectation in the right-hand side $\OO(N\eta^3\tau^2\frac{1}{\tau^2})=\OO(N^{-1/2+\kappa/5})$ using $|\frac{\tilde{\mathbf{K}}(z,z)}{N}|+|\frac{\Delta V(z)}{4\pi}|=\OO(1)$; because the integral over $w$ is $\OO(\eta^3)$ and $|\partial h|=\OO(1/\tau^2)$ with the area of integration region is the support of $h$ which is $\OO(\tau^2)$. The second term is subpolynomial due to the $e^{-\frac{|z-w|^2}{2\eta^2}}$ term. The third integral is subpolynomial when evaluated for $z\in S_{\rm out}(\delta_N)$ due to Lemma \ref{lem:1point}.(ii), and the integral on the edge layer $S_{\rm edge}(\eta\log N)$ is bounded by $N (\eta \tau)\frac{1}{\tau^2}\eta^2$ by Lemma \ref{lem:1point}.(i).

 The only term left is the integral with $\rd\hat\mu_N(w)\rd \mu_{V}(z)$, which gives
\begin{align} 
\frac{N}{2}\int_{S}\E_f&\left[\int \bar\partial h(z)\frac{\bar z - \bar w}{z-w}e^{-\frac{|z-w|^2}{2\eta^2}}\rd\hat\mu_N(w)\right]\rd \mu_{V}(z)\nonumber
\\
=&\,\frac{N}{2}\int_{S}\bar\partial h(z)\int_{S_{\rm in}(\eta)} \frac{\bar z - \bar w}{z-w}e^{-\frac{|z-w|^2}{2\eta^2}}\big(\rho_{V}(w)+\OO(N^{-2\kappa}\log N)\big)\rd m(w)\rd \mu_{V}(z)\nonumber
\\
&+\frac{N}{2}\int_{S\setminus S_{\rm in}(\eta\log N)}\bar\partial h(z)\int_{S_{\rm in}(\eta)^{\rm c}} \frac{\bar z - \bar w}{z-w}e^{-\frac{|z-w|^2}{2\eta^2}}\frac{\wt{\bf K}_N(w,w)}{N}\rd m(w)\rd \mu_{V}(z)\nonumber
\\
&+\frac{N}{2}\int_{S_{\rm in}(\eta\log N)}\bar\partial h(z)\int_{S_{\rm in}(\eta)^{\rm c}} \frac{\bar z - \bar w}{z-w}e^{-\frac{|z-w|^2}{2\eta^2}}\frac{\wt{\bf K}_N(w,w)}{N}\rd m(w)\rd \mu_{V}(z)\label{eqn:E_1case3_3}
\end{align}
where we have used $|\frac{\wt{\bf K}_N}{N}(z,z)-\rho_{V}(z)|\leq N^{-2\kappa}\log N$ in the bulk. The third integral is subpolynomial due to the exponential term and Lemma \ref{lem:1point}.(i). The second integral is bounded by $N (\tau\eta\log N)\frac{1}{\tau^2}\eta^2\ll N^{-\kappa/2}$ by Lemma \ref{lem:1point}.(i). On the other hand, in the first integral, the contribution of $\OO(N^{-2\kappa}\log N)$ will be bounded by $N \eta^2 N^{-2\kappa}\log N$ which is negligible. The only term left is $\frac{N}{2}\int_{S}\bar\partial h(z)\int_{S_{\rm in}(\eta)} \frac{\bar z - \bar w}{z-w}e^{-\frac{|z-w|^2}{2\eta^2}}\rd\mu_{V}\rd \mu_{V}$. This can be easily bounded using the rotational symmetry as done above, more explicitly it's equal to:
\begin{multline*} 
\frac{N}{2}\int_{S_{\rm in}(2\eta\log N)}\bar\partial h(z)\int_{{\rm B}(z,\eta\log N)} \frac{\bar z - \bar w}{z-w}e^{-\frac{|z-w|^2}{2\eta^2}}(\rho_{V}(w)-\rho_{V}(z))\rd m(w)\rd \mu_{V}(z)
\\
+\frac{N}{2}\int_{S_{\rm in}(2\eta\log N)}\bar\partial h(z)\int_{S_{\rm in}(\eta)\setminus {\rm B}(z,\eta\log N)} \frac{\bar z - \bar w}{z-w}e^{-\frac{|z-w|^2}{2\eta^2}}\rho_{V}(w)\rd m(w)\rd \mu_{V}(z)
\\
+\frac{N}{2}\int_{S_{\rm edge}(2\eta\log N)}\bar\partial h(z)\int_{S_{\rm in}(\eta)} \frac{\bar z - \bar w}{z-w}e^{-\frac{|z-w|^2}{2\eta^2}}\rho_{V}(w)\rd m(w)\rd \mu_{V}(z)
\\
=\OO(N\eta^3)+\oo(N^{-1})+\OO(N(\tau\eta\log N)\frac{1}{\tau^2}\eta^2)=\OO(N^{-\kappa/2}).
\end{multline*}
This completes the proof.
\end{proof}

\subsection{Concentration of linear statistics.}\label{app:stability}\
We provide a uniform exponential tail bound for the centered linear statistics of sufficiently regular functions which implies the stability of overwhelmingly likely events as will be discussed at the end of this subsection. While Gaussian-type bounds for general inverse temperature $\beta$ are available in the literature (see e.g. \cite{LebSer2018, BauBouNikYau2017, serfaty2023gaussian}), we opt for a simple, self-contained argument relying only on the kernel estimates from Section \ref{sec:kernel}.
In contrast to the cited results, our approach treats the edge without difficulty, allowing us to cover the bulk regime in our main Theorem \ref{Theorem.2singularities} up to the mesoscopic distance $N^{-1/2+\kappa}$ from the boundary.

The proposition below highlights the essential role of the isotropy lemma, Lemma \ref{lem:anis}, in exploiting angular cancellations to achieve $\oo(1)$ estimates. In fact, if we only aim for an error of order $O((\log N)^c)$, the kernel estimates alone suffice as follows.

\begin{proposition}[Uniform Exponential Tail]\label{prop:exp_tail}
Fix $C,\kappa>0$. Then there exist $N_0=N_0(V,C,\kappa)$ such that, for every $\tau_0\in[N^{-1/2+\kappa},1]$ $z_0\in\C$, and $f\in\mathscr{A}_{3,\tau_0}$ with $\supp(f)\subset{\rm B}(z_0,\tau_0)$,
\begin{align*} 
\log\E\left[e^{tX_f}\right]\leq (\log N)^7\quad\textnormal{for all}\ |t|\leq C\log N\ \textnormal{and}\ N\geq N_0.
\end{align*}
Consequently, $\mathbb{P}\left(\left|X_f\right|>t\right)\leq e^{(\log N)^7-t}$ and $\E[|X_f|]\leq (\log N)^7$ for all $t\geq 0$ and $N\geq N_0$.
\end{proposition}

\begin{proof}
If $z_0\in S_{\rm out}(2\tau_0)$, then the result follows easily by Lemma \ref{lem:1point}.(ii). If $z_0\in S_{\rm in}(2\tau_0)$ it follows easily by Theorem \ref{prop:kernel} as follows. As in the proof of Proposition \ref{prop:Laplace_loop},
\begin{align*} 
\log\E (e^{tX_f})=\int_{0}^{1}\E_{stf}[tX_f]\rd s&=\int_0^{1}\int_{\C} t f(z)\big(\wt{\bf{K}}_N(z,z)-{\wt{\bf{K}}}^{\#}_N(z,z)\big)\rd m(z)\rd s
\\
&=\OO(\log N)\int_{S_{\rm in}(\tau_0)}\max_{s\in[0,1]}|\wt{\bf{K}}_N(z,z)-{\wt{\bf{K}}}^{\#}_N(z,z)|\rd m(z)
\end{align*}
where we suppress the dependence on $s$ in the notation of the kernel corresponding to the $stf$-tilted measure $\wt{\bf{K}}_N$. By Theorem \ref{prop:kernel}, the integral is $\OO(\log N)$, hence the bound follows. Thus, the only remaining case is $z_0\in S_{\rm edge}(2\tau_0)$. For simplicity, we can assume that $z_0\in\partial S$, because there exists a $z_0'\in\partial S$ such that $\supp(f)\subset{\rm B}(z_0',3\tau_0)$ and $\|f\|_{3,3\tau_0}\leq 1$.

By Proposition \ref{prop:Laplace_loop}, it suffices to show that $\max_{t\in[0,1]}\mathscr{E}_N(tf,f)\leq (\log N)^7$. Lemmas \ref{lem1} and \ref{lem2} remain valid, as their proofs do not invoke stability. The only modification needed is in the proof of Lemma \ref{lem3}. By Lemma \ref{lem:harm_ext_reg}, it suffices to show the following for any $\tau\in[\tau_0,1]$, continuous $h$ with $\sup_{\C\setminus\partial S}|\nabla h|\leq\frac{1}{\tau^2}$ and $\supp(h)\subset{\rm B}(z_0,\tau)$,
\begin{align*}
N\,\E_{f}\Big[\iint_{z\neq w} \frac{h(z)-h(w)}{z-w} \, \tilde \mu_N(\rd z) \, \tilde \mu_N(\rd w)\Big]=\OO((\log N)^6).
\end{align*}
Let $H(z,w)=\frac{h(z)-h(w)}{z-w}$.  The left-hand side equals to $I(\C,\C)-J$, where we define
\begin{align*} 
I(A,B)&=N\int_{A}\int_{B} H(z,w)\Big(\frac{\wt{\bf{K}}_N(z,z)}{N}-\rho_V(z)\Big)\Big(\frac{\wt{\bf{K}}_N(w,w)}{N}-\rho_V(w)\Big)\rd m(w)\rd m(z),
\\
J&=N\iint H(z,w)\frac{|\wt{\bf{K}}_N(z,w)|^2}{N^2}\rd m(w)\rd m(z).
\end{align*}
The term $J$ can be bounded easily using the reproducing property of the kernel \eqref{eqn:reprod} and Lemma \ref{lem:1point}.(i), as follows:
\begin{align*} 
|J|\leq 2N\iint_{z\in{\rm B}(z_0,\tau)} |H(z,w)|\frac{|\wt{\bf{K}}_N(z,w)|^2}{N^2}\rd m(w)\rd m(z)=\OO(\frac{1}{N \tau^2})\iint_{z\in{\rm B}(z_0,\tau)} |\wt{\bf{K}}_N(z,w)|^2\rd m(w)\rd m(z)
\\
=\OO(\frac{1}{N \tau^2})\int_{z\in{\rm B}(z_0,\tau)} \wt{\bf{K}}_N(z,z)\rd m(z)=\OO(1)
\end{align*}
where in the first line we used that $H(z,w)$ is $0$ when both $z,w\notin{\rm B}(z_0,\tau)$.

We now turn to the integral $I$, which will be analysed according to the locations of $z$ and $w$. When at least one of the variables $z$ and $w$ lies in $S_{\rm out}(\delta_N)$, the integral is already negligible by Lemma \ref{lem:1point}.(ii), in particular $I(S_{\rm out}(\delta_N),\C)$ is subpolynomial. Thus, it remains to consider the three remaining configurations where $z$ and $w$ belong to $S_{\rm edge}(\delta_N)$ or $S_{\rm in}(\delta_N)$.\\

\noindent\textit{First case: If $z,w\in S_{\rm edge}(\delta_N)$.}  We rely only on the bound $\wt{\bf{K}}_N=\OO(N)$ from Lemma \ref{lem:1point}.(i):
\begin{align*} 
I(S_{\rm edge}(\delta_N),S_{\rm edge}(\delta_N))=N \iint_{S_{\rm edge}(\delta_N)^2}\frac{|h(z)|}{|z-w|} \OO(1)\rd m(z)\rd m(w)=\OO(N\delta_N\log N)\int_{S_{\rm edge}(\delta_N)\cap {\rm B}(z_0,\tau)}|h(z)|\rd m(z)
\\
=\OO((\log N)^5).
\end{align*}

\noindent\textit{Second case: If $z\in S_{\rm in}(\delta_N),w\in S_{\rm edge}(\delta_N)$.}  Using Theorem \ref{prop:kernel} for $z$, we obtain
\begin{align*} 
I(S_{\rm in}(\delta_N), S_{\rm edge}(\delta_N))=N \iint_{\substack{z\in S_{\rm in}(\delta_N),\\w\in S_{\rm edge}(\delta_N)}}|H(z,w)|\OO\Big(\frac{\log N}{N\tau_0^2}\mathds{1}_{{\rm B}(z_0,2\tau_0)}(z)+N^{-3/2}\Big)\rd m(z)\rd m(w).
\end{align*}
Note that
\begin{multline*} 
N \iint_{\substack{z\in S_{\rm in}(\delta_N)\cap {\rm B}(z_0,2\tau_0),\\ w\in S_{\rm edge}(\delta_N)}}|H(z,w)|\frac{1}{N\tau_0^2}\rd m(z)\rd m(w)=\frac{1}{\tau_0^2\tau} \iint_{\substack{z\in S_{\rm in}(\delta_N)\cap {\rm B}(z_0,2\tau_0),\\w\in S_{\rm edge}(\delta_N)}}\frac{1}{|z-w|}\rd m(z)\rd m(w)
\\
=O(\frac{\delta_N^2\log N}{\tau\tau_0})=\oo(1),
\end{multline*}
and the contribution from $N^{-3/2}$ is easily seen to be $\oo(1)$.

\noindent\textit{Third case: If $z,w\in S_{\rm in}(\delta_N)$.} By Theorem \ref{prop:kernel} applied to both $z$ and $w$, we have
\begin{multline*} 
I(S_{\rm in}(\delta_N), S_{\rm in}(\delta_N))
\\
=\OO(N)\iint_{S_{\rm in}(\delta_N)^2} |H(z,w)|\Big(\frac{\log N}{N\tau_0^2}\mathds{1}_{{\rm B}(z_0,2\tau_0)}(z)+N^{-3/2}\Big)\Big(\frac{\log N}{N\tau_0^2}\mathds{1}_{{\rm B}(z_0,2\tau_0)}(w)+N^{-3/2}\Big)\rd  m(z)\rd m(w).
\end{multline*}
We estimate the contributions from each terms as follows:
\begin{align*} 
N\iint_{(S_{\rm in}(\delta_N)\cap {\rm B}(z_0,2\tau_0))^2} |H(z,w)|\frac{1}{N\tau_0^2}\frac{1}{N\tau_0^2}\rd  m(z)\rd m(w)=\OO(\frac{1}{N\tau^2})=\oo(1)
\end{align*}
and
\begin{align*} 
\OO(N)\iint_{S_{\rm in}(\delta_N)^2} |H(z,w)|N^{-3/2}\rd  m(z)\rd m(w)=\OO(\frac{1}{N^{1/2}\tau})=\oo(1),
\end{align*}
concluding the proof.
\end{proof}

By H\"{o}lder's inequality Proposition \ref{prop:exp_tail} implies the following corollary.

\begin{corollary}\label{prop:Laplace_S_kappa} For any $C,\kappa>0$, there exist $N_0=N_0(C,\kappa)$ such that, for every $f\in\mathscr{S}_{3,C,\kappa}$,
\begin{align*} 
\log\E\left[e^{tX_f}\right]\leq (\log N)^7\quad\textnormal{for all}\ t\in[-C,C]\ \textnormal{and}\ N\geq N_0.
\end{align*}
Consequently, $\mathbb{P}\left(\left|X_f\right|>t\right)\leq e^{(\log N)^7-t}$ and $\E[|X_f|]\leq (\log N)^7$ for all $t\geq 0$ and $N\geq N_0$.
\end{corollary}

By \textit{stability} we mean that an event $A$ with overwhelming probability, $\mathbb{P}(A)\geq 1-e^{-(\log N)^{D}}$, remains overwhelmingly likely under the biased measure $\mathbb{P}_f$ given by $\mathbb{P}_f(A)=\mathbb{E}\big[\mathds{1}_{A}\frac{e^{X_f}}{\mathbb{E}e^{X_f}}\big]$, $X_f=\sum_{i}f(z_i)-N\int f\rd\mu_{V}$. Using Cauchy-Schwarz and Jensen’s inequalities, the results in this section yield stability for $\mathscr{S}_{3,C,\kappa}$ biases. We record this as a corollary for ease of reference.

\begin{corollary}[Stability]\label{cor:stability}
Fix $C,\kappa>0$. For every $f\in\mathscr{S}_{3,C,\kappa}$, if an event $A$ satisfies $\mathbb{P}(A)\geq1-e^{-(\log N)^D}$ for some $D\geq10$, then $\mathbb{P}_{f}(A)\geq 1-e^{-(\log N)^D/2}$.
\end{corollary}

\subsection{Local isotropy.}\ 
In the Lemma below, we define $(\bar z-\bar w)/(z-w)=0$ for $z=w$.

\begin{lemma}\label{lem:anis} Consider a fixed small parameter $\kappa>0$ and let $\delta=N^{-\frac{1}{2}+\kappa}$. Then, for any $D>0$, there exists an $N_0=N_0(\kappa,D)\in\N$ such that for every $\tau\in[N^{-1/2+15\kappa},1]$, $z_0\in\C$ and continuous function $g:\mathbb{C}\to\mathbb{R}$ satisfying $\supp(g)\subseteq {\rm B}(z_0,\tau)$ and $\|g\|_{\infty}\leq 1$, we have
\begin{align}\label{eqn:anisotropy}
\mathbb{P}\left(
\left|\sum_{1\leq i,j\leq N}
g(z_i)\frac{\overline{z_i}-\overline{z_j}}{z_i-z_j}e^{-\frac{|z_i-z_j|^2}{\delta^2}}
\right|>\tau^2N^{1-\kappa}
\right)\leq e^{-(\log N)^{D}}
\end{align}
for all $N\geq N_0$.
\end{lemma}

Note that if no cancellation were expected, with overwhelming probability, the number of indices $i$ with $g(z_i)\asymp 1$ would be $\asymp \tau^2 N$. For each such $i$, the sum over $j$'s for which $\frac{\overline{z_i}-\overline{z_j}}{z_i-z_j}e^{-\frac{|z_i-z_j|^2}{\delta^2}}\asymp 1$ would contribute about $N^{2\kappa}$. Consequently, one would expect the total sum to be on the order of $\tau^{2}N^{1+2\kappa}$. The lemma above shows, however, that angular cancellations occur and the apparent leading contribution is substantially reduced.
Some angular cancellations were also instrumental in  \cite{LebSer2018,BauBouNikYau2019}, in which isotropy was obtained differently from partition function asymptotics, for arbitrary inverse temperature.  
The above lemma is specific to the determinantal setting, and the proofs proceeds very differently, to cover cancellations up to the edge of the droplet,  and for general biased measures.

\begin{proof} The proof follows by four main steps.\\

\noindent\textit{First step: reducing to the bulk.}
Let $\chi$ be a smooth bump function that is equal to $1$ on $S_{\rm in}(2N^{2\kappa}\delta)$, and $0$ outside $S_{\rm in}(N^{2\kappa}\delta)$. We first consider the sum
\begin{align*} 
\sum_{1\leq i,j\leq N}
g(z_i)(1-\chi(z_i))\frac{\overline{z_i}-\overline{z_j}}{z_i-z_j}e^{-\frac{|z_i-z_j|^2}{\delta^2}}.
\end{align*}
Note that $S_{\rm edge}(2N^{2\kappa}\delta)$ can be covered by $\asymp (N^{2\kappa}\delta)^{-1}$ many balls of radius $4N^{2\kappa}\delta$. Applying Proposition \ref{prop:Laplace_S_kappa} to each ball and using a union bound we obtain that, with overwhelming probability, the number of particles in ${\rm B}(z_0,\tau)\cap S_{\rm edge}(2N^{2\kappa}\delta)$ is bounded by $\tau N^{1/2+4\kappa}$. Again with overwhelming probability, for each particle $z_i$ in this outer ring, the number of $z_j$'s inside ${\rm B}(z_i,N^{\kappa}\delta)$ is bounded by $N^{7\kappa}$; because ${\rm B}(z_i,N^{\kappa}\delta)$ is a subset of union of $\OO(1)$ many balls of radius $4N^{2\kappa}\delta$ constructed above. This leads to a rough upper bound $\tau N^{1/2+11\kappa}$ which is smaller than $N^{1-2\kappa}\tau^2$. On the other hand, the contribution from $z_i$ with $\dist(z_i,S)>\delta$ is already negligible by Lemma \ref{lem:1point}.(ii) and Markov's inequality. 

 So, it suffices to prove the equation \eqref{eqn:anisotropy} when $g$ is replaced by $g\chi$. For notational simplicity, we will henceforth write $g$ in place of $g\chi$.\\
 
\noindent\textit{Second step: grouping and centering.} Consider $\mathscr{L}=(2\delta N^{\kappa}\mathbb{Z}^2)\cap {\rm B}(z_0,\tau)\cap \{z\in S:\dist(z,\partial S)>N^{2\kappa}\delta\}$, $\ell=|\mathscr{L}|$. Let $\mathscr{S}$ be the set of squares with vertices $w\pm \delta N^\kappa\pm \ii \delta N^\kappa$, where $w\in\mathscr{L}$. For any $\alpha\in\mathscr{S}$ we define
\begin{align*}
X_\alpha=\sum_{1\leq i,j\leq N}f_{\alpha}(z_i,z_j),\ \ 
f_{\alpha}(w,z)=
g(w)\frac{\overline{w}-\overline{z}}{w-z}e^{-\frac{|w-z|^2}{\delta^2}}\mathds{1}_{w\in\alpha}\mathds{1}_{|w-z|<\delta N^{\kappa/2}}.
\end{align*}
We have the deterministic estimate
\begin{align*}
\left|\sum_{1\leq i,j\leq N}
g(z_i)\frac{\overline{z_i}-\overline{z_j}}{z_i-z_j}e^{-\frac{|z_i-z_j|^2}{\delta^2}}-\sum_{\alpha\in\mathscr{S}} X_\alpha\right|\leq N^2 e^{-N^{\kappa}}
\end{align*}
so we only need to prove that, under the same hypotheses as Lemma \ref{lem:anis} we have
\begin{equation*}
\mathbb{P}\left(
\left|\sum_{\alpha\in\mathscr{S}} X_\alpha
\right|>\tau^2N^{1-\kappa}
\right)\leq e^{-(\log N)^{D}}.
\end{equation*}
Denote the centered $X_{\alpha}$ by $Y_{\alpha}=X_{\alpha}-\E [X_{\alpha}]$. Note that for any $\alpha\in\mathscr{S}$, 
\begin{align} \label{eqn:EX_alpha}
\E [X_{\alpha}]=\iint g(w)\frac{\overline{w}-\overline{z}}{w-z}e^{-\frac{|w-z|^2}{\delta^2}}\mathds{1}_{w\in\alpha}\mathds{1}_{|w-z|<\delta N^{\kappa/2}}\rho_2^{(N)}(z,w) \rd m(z)\rd m(w)
\end{align}
where $\rho_2^{(N)}(z,w)={\bf K}_N(z,z){\bf K}_N(w,w)-|{\bf K}_N(z,w)|^2$ is the two-point correlation function. Using Theorems \ref{prop:kernel} and \ref{thm:off_diag_decay} with the Taylor expansion in \eqref{eqn:taylorforQ} we can write
\begin{align*} 
\rho_2^{(N)}(z,w)=\frac{N^2}{\pi^2}(\partial\bar{\partial}V(w))^2\big(1-e^{-N|z-w|^2\frac{\Delta V(w)}{4}}\big)+\OO(N^{3/2+2\kappa})
\end{align*}
when both $z$ and $w$ are in $\{\xi\in S:\dist(\xi,S)<N^{-1/2+2\kappa}\}$ with $|z-w|<N^{-1/2+2\kappa}$. Due to the rotational symmetry of the integrand (over the variable $z$) in \eqref{eqn:EX_alpha}, this leads to a rough bound 
\begin{align*} 
\E [X_{\alpha}]=\OO(N^{3/2+2\kappa}(\delta N^{\kappa/2})^2(\delta N^{\kappa})^2)=\OO(N^{-1/2+9\kappa}).
\end{align*}
Thus $\sum_{\alpha}\E [Y_{\alpha}]=\OO(N^{-1/2+9\kappa}|\mathscr{S}|)=\tau^2N^{1-\kappa}\OO(N^{-1/2+6\kappa})$. So, given $\kappa<1/13$, it suffices to prove
\begin{equation}\label{eqn:enough2}
\mathbb{P}\left(
\left|\sum_{\alpha\in\mathscr{S}} Y_\alpha
\right|>\tau^2N^{1-\kappa}
\right)\leq e^{-(\log N)^{D}}.
\end{equation}\\

\noindent\textit{Third step: decoupling.}
For a given exponent $p\geq 5$ (possibly $N$-dependent) we expand
\begin{equation}\label{eqn:boundelem}
\E\big[|\sum_{\alpha\in  \mathscr{S}}Y_\alpha|^{2p}\big]=\sum_{\substack{p_\alpha,q_\alpha\geq 0, \alpha\in \mathscr{S},\\ \sum_{\alpha\in\mathscr{S}}p_\alpha=\sum_{\alpha\in\mathscr{S}}q_\alpha=p}}\frac{(2p)!}{\prod_{\alpha\in\mathscr{S}}p_\alpha!q_\alpha!}\E\big[\prod_{\alpha\in  \mathscr{S}}Y_{\alpha}^{p_\alpha}\overline{Y}_{\alpha}^{q_\alpha}\big]
\leq
p^{2p}\sum_{\substack{I\subset  \mathscr{S},\\|I|\leq 2p}}\sum_{\substack{\sum_{\alpha\in I}p_\alpha=\sum_{\alpha\in I}q_\alpha=p,\\ p_\alpha+q_\alpha\geq 1,\\p_\alpha,q_\alpha\geq 0}}\left|\E\big[\prod_{\alpha\in I}Y_\alpha^{p_\alpha}\overline{Y}_{\alpha}^{q_\alpha}\big]\right|.
\end{equation}

To bound the above expectations, we now consider a given $I\subset \mathscr{S}$ and exponents $p_\alpha+q_\alpha\geq 1$. Let $G$ be the graph with vertices $I$ and edges $\{\{w_1,w_2\}\in\mathscr{S}: w_1\ {\rm and }\ w_2\ \textnormal{are distinct and adjacent horizontally, vertically,}\ \allowbreak\textnormal{or diagonally}\}$.  Denote $I=\cup_{i=1}^{{\rm c}(I)} I_i$ the decomposition of $I$ into connected components for $G$ (${\rm c}(I)$ is the number of connected components).

For any $i$, denote $m_i=|I_i|$ and $n_i=2\sum_{\alpha\in I_i} p_\alpha+q_{\alpha}$. 
Let $\Theta_r$ (resp. $\Theta_r^*$) be the set of $r$-tuples of elements (resp. distinct elements) in $\{z_1,\dots,z_N\}$.
We expand
\begin{equation}\label{eqn:interm1}
\prod_{\alpha\in I_i}Y_{\alpha}^{p_\alpha}\overline{Y}_{\alpha}^{q_\alpha}
=
\sum_{\bz\in\Theta_{n_i}}f^{(i)}(\bz)\ {\rm with}\ f^{(i)}(\bz)=\prod_{\alpha\in I_i}\left(\prod_{k=1}^{p_\alpha}f_{\alpha}(w^{(\alpha)}_k,z^{(\alpha)}_k)\prod_{k=p_\alpha+1}^{p_\alpha+q_\alpha}\overline{f_{\alpha}(w^{(\alpha)}_k,z^{(\alpha)}_k)}\right)
\end{equation}
where we have decomposed $\bz=(\bw^{(\alpha)},\bz^{(\alpha)})_{\alpha\in I_i}$ and $\bw^{(\alpha)},\bz^{(\alpha)}$  have length $p_\alpha+q_\alpha$. 
Let $d(\bw)$ be the number of distinct eigenvalues appearing in the tuple $\bw$.
Equation \eqref{eqn:interm1} implies that if we define 
\begin{align*}
f^{(i)}_r(\bz)=
\sum_{\bw\in\{z_1,\dots,z_r\}^{n_i},\ d(\bw)=r}f^{(i)}(\bw),
\end{align*}
then we have
\begin{equation}\label{eqn:decomp}
\prod_{\alpha\in I_i}Y_{\alpha}^{p_\alpha}\overline{Y}_{\alpha}^{q_\alpha}
=\sum_{r\leq n_i}\sum_{\bz\in\Theta_r^*}f^{(i)}_r(\bz).
\end{equation}
This expansion over $r$-tuples of distinct eigenvalues allows to write expectations in terms of correlation functions.
Note that there is a constant $C>0$ depending only on $p$ such that  $\|f_r^{(i)}\|_\infty\leq C$. Moreover, $f_r^{(i)}$ is supported on the $\delta N^{\kappa/2}$-neighbourhood of  $(\cup_{\alpha\in I_i} \alpha)^r$,
so that if $f_r^{(i)}(\bz)\neq 0$ and $f_s^{(j)}(\bw)\neq 0$ or $i\neq j$, then the minimal distance between points from $\bz$ and $\bw$ is at least $\delta N^{\kappa}/2$. In particular the supports of $f_r^{(i)}$ and 
$f_s^{(j)}$ are distinct so that
\begin{align*}
\prod_{i=1}^{{\rm c}(I)}\sum_{r\leq n_i}\sum_{\bz\in\Theta_r^*}f^{(i)}_r(\bz)
=
\sum_{r_i\leq n_i}\sum_{\bz\in\Theta^*_{r_1+\dots+r_{{\rm c}(I)}}}g_{r_1,\dots,r_{{\rm c}(I)}}(\bz),\ {\rm where}\ 
g_{r_1,\dots,r_{{\rm c}(I)}}(\bz)
=
\prod_{i=1}^{{\rm c}(I)}f^{(i)}_{r_i}(\bz_{r_i}).
\end{align*}
Here, we have decomposed  $\bz=(\bz_{r_1},\dots,\bz_{r_{{\rm c}(I)}})$ with $\bz_{r_i}\in\Theta_{r_i}^*$.

Moreover, by definition of the correlation functions we have
\begin{equation}\label{eqn:g}
\E\bigg[\sum_{\bz\in\Theta^*_{r_1+\dots+r_{{\rm c}(I)}}}g_{r_1,\dots,r_{{\rm c}(I)}}(\bz)
\bigg]
=
\int_{\mathbb{C}^{r_1+\dots+r_{{\rm c}(I)}}} g_{r_1,\dots,r_{{\rm c}(I)}}(\bz)
\rho_{r_1+\dots+r_{{\rm c}(I)}}^{(N)}(\bz)
\rd m(\bz).
\end{equation}
When $g_{r_1,\dots,r_{{\rm c}(I)}}(\bz)\neq 0$, as mentioned before the minimal distance between $\bz_{r_i}$ and $\bz_{r_j}$ for $i\neq j$ is at least $\delta N^\kappa/2$, so that, using Theorem \ref{thm:off_diag_decay} and assuming $p\leq (\log N)^{\log\log N}$, we obtain 
\begin{align}
\rho_{r_1+\dots+r_{{\rm c}(I)}}^{(N)}(\bz)&=\det_{\sum r_i\times\sum r_i}{\bf K}_N(z_i,z_j)\notag\\
&=
\prod_{k=1}^{{\rm c}(I)}\det_{r_k\times r_k}{\bf K}_N(z_{r_k}(i),z_{r_k}(j))
+
\OO\left(
p!\|{\bf K}_N\|_\infty^p\max_{\substack{|z-w|>N^\kappa\delta/2,\\|z|,|w|\leq 1-\delta}}|{\bf K}_N(z,w)| 
\right)\notag\\
&=
\rho_{r_1}^{(N)}(\bz_{r_1})
\cdots
\rho_{r_{{\rm c}(I)}}^{(N)}(\bz_{r_{{\rm c}(I)}})+
\OO\left(e^{-cN^{2\kappa}}\right)\label{eqn:deco}.
\end{align}
From \eqref{eqn:g} and \eqref{eqn:deco} we have
\begin{multline*}
\E\Big[\prod_{i=1}^{{\rm c}(I)}\sum_{r\leq n_i}\sum_{\bz\in\Theta_r^*}f^{(i)}_r(\bz)\Big]=\sum_{r_i\leq n_i}\left(
\prod_{i=1}^{{\rm c}(I)}\E \Big[\sum_{\bz_{r_i}\in\Theta_{r_i}^*}f^{(i)}_{r_i}(\bz_{r_i})\Big]+
\OO\left(e^{-cN^{2\kappa}}\right)\right)\\
=\prod_{i=1}^{{\rm c}(I)}\E \Big[\sum_{r\leq n_i}\sum_{\bz\in\Theta_{r}^*}f_r^{(i)}(\bz)
\Big]+
\OO\left(e^{-cN^{2\kappa}}\right).
\end{multline*}
With the definition \eqref{eqn:decomp}, this means that there exists $N_0>0$ such that for any $p\leq(\log N)^{\log\log N}$, $N\geq N_0$,
$I\subset\mathscr{S}$ and $\sum_{\alpha\in I}p_\alpha+q_\alpha=2p$ with $p_\alpha+q_\alpha\geq 1$, $p_\alpha,q_\alpha\geq 0$ we have
\begin{equation*}
\E\Big[\prod_{i=1}^{{\rm c}(I)}\prod_{\alpha\in I_i}Y_{\alpha}^{p_\alpha}\overline{Y}_{\alpha}^{q_\alpha}\Big]
=
\prod_{i=1}^{{\rm c}(I)}\E\Big[\prod_{\alpha\in I_i}Y_{\alpha}^{p_\alpha}\overline{Y}_{\alpha}^{q_\alpha}\Big]+
\OO\left(e^{-cN^{2\kappa}}\right).
\end{equation*}\\

\noindent\textit{Fourth step: counting and Markov.}
We can now estimate the right-hand side of \eqref{eqn:boundelem}. First, observe that the contribution of $I$'s with ${\rm c}(I)>p$ is zero. Indeed, for any such $I$, there exists $i\in\llbracket 1,{\rm c}(I)\rrbracket$ such that $I_i=\{\alpha\}$, and $(p_\alpha,q_\alpha)=(1,0)$ or $(p_\alpha,q_\alpha)=(0,1)$. Due to the centering, we already have $\E [Y_\alpha]=\E[\bar{Y}_{\alpha}]=0$, which forces the entire corresponding term in the sum to be zero.

 For any $I$ such that ${\rm c}(I)\leq p$, we simply bound $|X_\alpha|\leq |\{z_i:{\rm d}(z_i,\alpha)\leq \delta N^\kappa\}|^2$ and $\E [X_\alpha]\leq N^{8\kappa}$. So, by Corollary \ref{prop:Laplace_S_kappa}, we have
 $
\left|\E\big[\prod_{\alpha\in I}Y_\alpha^{p_\alpha}\overline{Y}_{\alpha}^{q_\alpha}\big]\right|\leq N^{17p\kappa}
 $.
 For any $I\subset\mathscr{S}$ satisfying $|I|\leq 2p$, we have $|\{(p_\alpha,q_\alpha)_{\alpha\in I}:p_{\alpha},q_\alpha\geq 0,\sum_{\alpha\in I}p_\alpha=\sum_{\alpha\in I}q_\alpha=p\}|\leq p^{7p}$. The contribution to \eqref{eqn:boundelem} is therefore bounded by $p^{2p}N^{17p\kappa}p^{7p}|\{I\subset\mathscr{S}:|I|\leq 2p, {\rm c}(I)\leq p\}|\leq p^{10p} N^{17p\kappa}|\mathscr{S}|^p (18p)^{2p}\leq  p^{15p}\tau^{2p}N^{p(1+13\kappa)}$ for any $p>20$.
 
 In sum, we have proved that for any $20<p\leq (\log N)^{\log\log N}$ we have
\begin{align*}
\E\left[\Big|\sum_{\alpha\in  \mathscr{S}}Y_\alpha\Big|^{2p}\right]
\leq
p^{15p}\tau^{2p}N^{p(1+13\kappa)},
\end{align*}
so by Markov's inequality
\begin{align*}
\mathbb{P}\left(
\Big|\sum_{\alpha\in\mathscr{S}} Y_\alpha
\Big|>\tau^2N^{1-\kappa}
\right)\leq 
N^{-p10\kappa }
\end{align*}
where we have used $p\leq (\log N)^{\log\log N}$. 
For fixed $\kappa,D>0$, the choice $p=(\log N)^{D}/ (10\kappa)$ concludes the proof of \eqref{eqn:enough2} and the lemma.
\end{proof}

\section{Singular potentials}\label{sec:sing_pot}

In this section, we prove Laplace transform estimates similarly to Proposition \ref{proposition.nosingularities}, but now with logarithmic singularities. The approach contains key modifications compared to the previous section,  in particular a decoupling of local singularities and a comparison with a reference model,  the electric potential for the  Ginibre ensemble at $\zeta=0$. At a high level,  this section mirrors a method initiated in \cite{BouFal2025} in the context of Fisher-Hartwig singularities in a space-time setting.  

\subsection{Regularizing and localizing log-singularities.}\
Recall the logarithmic regularization $\log_{\epsilon}=\log|\cdot|\ast\frac{\chi_{\epsilon}}{\|\chi_{\epsilon}\|_{L^1}}$ defined for $\epsilon>0$. We begin by noting that the root-type singularities in Theorem \ref{Theorem.2singularities} can be replaced by their submicroscopic regularizations. Indeed,  compared to  \cite{BouFal2025,Keles2025} the submicroscopic smoothing is considerably simpler for two dimensional Coulomb gasses. Owing to the subharmonicity and harmonicity properties of the logarithm in 2d, there is no need to invoke Hua-Pickrell type kernel estimates as in \cite[Lemmas 2.4–2.5]{BouFal2025} or to use change of variables, see \cite[Proposition 6.1]{Keles2025}.

\begin{proposition}\label{prop:submic_log_reg}
Let $C>0$, $m\in\N$, $\alpha,\kappa>0$, $\Delta=N^{-1/2-\alpha}$, $f\in\mathscr{S}_{3,C,\kappa}$, $\gamma_1,\dots,\gamma_m\in[0,C]$ and $\zeta_1,\dots,\zeta_m\in S_{\rm in}(2\delta_N)$. Then
\begin{align*} 
\mathbb{E}\Big[e^{\tr(f+\sum_{j=1}^{m}\gamma_j\log^{\zeta_j})}\Big]=\mathbb{E}\Big[e^{\tr(f+\sum_{j=1}^{m}\gamma_j\log_{\Delta}^{\zeta_j})}\Big](1+\OO(N^{-2\alpha})).
\end{align*}
\end{proposition}

\begin{proof}
Because $\chi$ is radial and $\log$ is subharmonic, $\log_{\Delta}\geq\log$, which proves one side of the proposition. For the other direction, note that $\log_{\Delta}=\log$ outside ${\rm B}(0,\Delta)$ due to harmonicity of $\log$. Defining the set $\mathcal{G}=\{\textnormal{there is no eigenvalue in }\cup_{j=1}^{m}{\rm B}(\zeta_j,\Delta)\}$, we can write
\begin{align*} 
\mathbb{E}\Big[e^{\tr(f+\sum_{j=1}^{m}\gamma_j\log_{\Delta}^{\zeta_j})}\Big]\leq \mathbb{E}\Big[e^{\tr(f+\sum_{j=1}^{m}\gamma_j\log^{\zeta_j})}\Big]+\mathbb{E}\Big[e^{\tr(f+\sum_{j=1}^{m}\gamma_j\log_{\Delta}^{\zeta_j})}\mathds{1}_{\mathcal{G}^{\rm c}}\Big].
\end{align*}
Moreover, using the estimation on the one-point correlation function in Lemma \ref{lem:1point}.(i) we obtain
\begin{align*} 
\frac{\mathbb{E}\Big[e^{\tr(f+\sum_{j=1}^{m}\gamma_j\log_{\Delta}^{\zeta_j})}\mathds{1}_{\mathcal{G}^{\rm c}}\Big]}{\mathbb{E}\Big[e^{\tr(f+\sum_{j=1}^{m}\gamma_j\log_{\Delta}^{\zeta_j})}\Big]}=\mathbb{E}_{f+\sum_{j=1}^{m}\gamma_j\log_{\Delta}^{\zeta_j}}\big[\mathds{1}_{\mathcal{G}^{\rm c}}\big]\leq \int_{\cup_{j=1}^{m}{\rm B}(\zeta_j,2\Delta)}\wt{\bf K}_N(z,z)\rd m(z)=\OO(N^{-2\alpha})
\end{align*} 
which completes the proof.
\end{proof}

 Fix $\kappa>0$, $\delta=N^{-1/2+\kappa}$, $\zeta\in\C$ and let $\chi$ be a smooth bump function such that $\chi=1$ on ${\rm B}(\zeta,1)$ and $0$ outside ${\rm B}(\zeta,2)$. For any $\zeta\in S_{\rm in}(\delta_N)$, one has $\chi\log_{\delta}^{\zeta}\in\mathscr{S}_{C,\kappa}$ for some suitable constant $C$ (depending only on $S$). This follows from the multiscale decomposition:
\begin{align*}
\log_{\delta}(z-\zeta)=\sum_{n=0}^{\infty} g_n(z),\ g_n(z)=\log_{\delta}(z-\zeta)\chi_n(N^{\frac{1}{2}-\kappa}|z-\zeta|),
\end{align*}
where $1=\sum_{n\geq 0} \chi_n$ is a partition of $\mathbb{R}_+$, $\chi_0$ is supported on $[0,2]$, $\chi_n$ is supported on $[2^{n-1},2^{n+1}]$, and $\|\chi_n^{(p)}\|_\infty\leq C_p 2^{-np}$ for $n,p\geq 0$. It's straightforward to verify that $g_n\in \mathscr{A}_{C\delta2^{n}}$ for some constant $C$. Using this decomposition for indices $n$ such that $2^n<\delta$ yields the desired result $\chi\log_{\delta}^{\zeta}\in\mathscr{S}_{C,\kappa}$.

\begin{proposition}\label{proposition.singularities}
Let $\kappa,\alpha>0$ be small, $C>0$, $m\in\N$ be some fixed constants, $\Delta=N^{-1/2-\alpha}$ and $\delta=N^{-1/2+\kappa}/3$. Then, uniformly in $f\in\mathscr{S}_{5,C,\kappa}$, $\zeta_1,\dots,\zeta_m\in S_{\rm in}(N^{-1/2+\kappa})$ satisfying $\min_{i\neq j}|\zeta_i-\zeta_j|>N^{-1/2+\kappa}$ and $\gamma_1,\dots,\gamma_m\in[0,C]$ we have
\begin{align*} 
\E\Big[e^{\sum_{i=1}^Nf(z_i)}\prod_{j=1}^{m}e^{\sum_{i=1}^N\gamma_j\log_{\Delta}(z_i-\zeta_j)}\Big]
=e^{N\int \big(f+\sum_{j}\gamma_j\log_{\delta}^{\zeta_j}\big)\rd\mu_{V}}
e^{\frac{1}{8\pi}\int_{\C}\big(|\nabla f^S|^2
+\mathds{1}_{S}\Delta f
+L^{S}\Delta f\big)\rd m
}
\\
\times\prod_{j=1}^{m}e^{\frac{\gamma_j}{2}\big(\int_{\partial S}f\rd\omega^{\infty}-f(\zeta_j)\big)} 
\prod_{j=1}^{m}e^{\frac{\gamma_j^2}{8\pi}\int_{\C}\nabla(\log_{\delta}^{\zeta_j})^{S}\cdot\nabla(2\log_{\Delta}^{\zeta_j}-\log_{\delta}^{\zeta_j})^{S}\,\rd m
+\frac{\gamma_j}{4}(1+L(\zeta_j)-L^{S}(\infty))}
\\
\times\E\Big[\prod_{j=1}^{m}e^{\sum_{i=1}^N\gamma_j(\log_{\Delta}^{\zeta_j}-\log_{\delta}^{\zeta_j})(z_i)}\Big]
\prod_{j\neq k}|\zeta_j-\zeta_k|^{-\frac{\gamma_i\gamma_j}{4}}
e^{\frac{\gamma_i\gamma_j}{4}\mathfrak{s}} \left(1+\OO(N^{-\kappa/10})\right).
\end{align*}
\end{proposition}

\begin{proof} We begin similarly to the proof of Proposition \ref{proposition.nosingularities}, with the exception that we remove only the function $f$ and the mesoscopic contributions from logarithmic singularities using Proposition \ref{prop:Laplace_loop}. More explicitly, let $\chi$ be an order $1$ bump function that is equal to $1$ on $S\cup S_{\rm edge}(1/2)$ and $0$ on $S_{\rm out}(2/3)$.  Then $F:=\sum_{j}\gamma_j(\log_{\Delta}^{\zeta_j}-\log_{\delta}^{\zeta_j})$ and $g:=f+\sum_{j}\gamma_j\log_{\delta}^{\zeta_j}$ substituted into $f$ and $g$ respectively (with $g_1=g\chi$, $g_2=g(1-\chi)$) in Proposition \ref{prop:Laplace_loop} yields
\begin{multline*} 
\log\E [e^{\sum_{i}(F+g)(z_i)}]=\log\E [e^{\sum_{i}F(z_i)}]+\Gamma(f)+\sum_{j=1}^{m}\gamma_j\,\Gamma(f,\log^{\zeta_j})+\sum_{j=1}^{m}\gamma_j\,\Gamma(\log^{\zeta_j})+\sum_{j=1}^{m}\gamma_j^2\,\Gamma(\log^{\zeta_j},\log^{\zeta_j})
\\
+\sum_{j\neq k}\gamma_j\gamma_k\,\Gamma(\log^{\zeta_j},\log^{\zeta_k})
+N\int g\,\rd\mu_{V}+\int_{0}^{1}\mathscr{E}_N(F+t g,g_1)\,\rd t
+\OO(N^{-D})
\end{multline*}
where
\begin{align*} 
\Gamma(f)&=\frac{1}{8\pi}\int_{\C}\big(|\nabla f^S|^2+\mathds{1}_{S}\Delta f+L^{S}\Delta f\big)\rd m,
\\
\Gamma(f,\log^{\zeta_j})&= \frac{1}{4\pi} \int_{C}\nabla f^{S}\cdot\nabla(\log_{\Delta}^{\zeta_j})^{S}\rd m,
\\
\Gamma(\log^{\zeta_j})&=\frac{1}{8\pi}\int_{S}\Delta \log_{\delta}^{\zeta_j}\,\rd m+\frac{1}{8\pi}\int_{S} \log_{\delta}^{\zeta_j}\Delta L\,\rd m-\frac{1}{8\pi}\int_{\partial S}\log_{\delta}^{\zeta_j}\mathcal{N}(L)\rd s,
\\
\Gamma(\log^{\zeta_j},\log^{\zeta_j})&=\frac{1}{8\pi}\int_{\C}\nabla(\log_{\delta}^{\zeta_j})^{S}\cdot\nabla(2\log_{\Delta}^{\zeta_j}-\log_{\delta}^{\zeta_j})^{S}\,\rd m,
\\
\Gamma(\log^{\zeta_j},\log^{\zeta_k})&=\frac{1}{8\pi}\int_{\C}\nabla(\log_{\delta}^{\zeta_j})^{S}\cdot\nabla(\log_{\Delta}^{\zeta_k})^{S}\,\rd m+\frac{1}{8\pi}\int_{\C}\nabla(\log_{\delta}^{\zeta_k})^{S}\cdot\nabla(\log_{\Delta}^{\zeta_j}-\log_{\delta}^{\zeta_j})^{S}\,\rd m.
\end{align*}
We now evaluate each term one-by-one. $\Gamma(f)$ is already in the desired form. Recalling the definition of the harmonic measure $\omega^{\infty}$ given in \eqref{eqn:Neumann_log} and the fact that $\Delta(\log^{\zeta}_{\Delta})=2\pi\frac{\chi^{\zeta}_{\Delta}}{\|\chi^{\zeta}_{\Delta}\|_{L^1}}$, $\Gamma(f,\log^{\zeta_j})$ can be simplified easily as follows:
\begin{align*} 
\Gamma(f,\log^{\zeta_j})&=\frac{1}{4\pi}\int_{S}\nabla f\cdot\nabla\log_{\Delta}^{\zeta_j}\rd m+\frac{1}{4\pi}\int_{S^{\rm c}}\nabla f^{S}\cdot\nabla(\log_{\Delta}^{\zeta_j})^{S}\rd m
\\
&=-\frac{1}{4\pi}\int_{S} f\Delta\log_{\Delta}^{\zeta_j}\rd m+\frac{1}{4\pi}\int_{\partial S}f\,\frac{\partial \log_{\Delta}^{\zeta_j}}{\partial n}\rd s+\frac{1}{4\pi}\int_{\partial S}f\,\frac{\partial (\log_{\Delta}^{\zeta_j})^{S}|_{S^{\rm c}}}{\partial (-n)}\rd s
\\
&=-\frac{1}{2} f(\zeta_j)+\frac{1}{2}\int_{\partial S}f\,\rd\omega^{\infty}+\OO(N^{-\kappa}).
\end{align*}
Next, for $\Gamma(\log^{\zeta_j})$ term, we begin by establishing the analogue of \eqref{eqn:green_to_neumann}, accounting for the non-compact support of the logarithm. We denote $L^{S}(\infty)$ by $L^{S}_{\infty}$ and modify $L^{S}$ to $L^{S}-L^{S}_{\infty}$ so that Green's identity can be applied on $S^{\rm c}$:
\begin{multline*} 
\int_{\C}(L^{S}-L^{S}_{\infty})\Delta \log_{\delta}^{\zeta_j} \rd m=\int_{S}(L-L^{S}_{\infty})\Delta \log_{\delta}^{\zeta_j} \rd m+\int_{S^{\rm c}}\Big((L^{S}-L^{S}_{\infty})\Delta \log_{\delta}^{\zeta_j} -\log_{\delta}^{\zeta_j}\Delta (L^{S}-L^{S}_{\infty})\Big)\rd m
\\
=\int_{S} \log_{\delta}^{\zeta_j}\Delta L\,\rd m+\int_{\partial S}\Big((L-L^{S}_{\infty})\frac{\partial \log_{\delta}^{\zeta_j}}{\partial n}-\log_{\delta}^{\zeta_j}\frac{\partial L}{\partial n}\Big)\rd s+\int_{\partial S}\Big((L-L^{S}_{\infty})\frac{\partial \log_{\delta}^{\zeta_j}}{\partial (-n)}-\log_{\delta}^{\zeta_j}\frac{\partial L^{S}|_{S^{\rm c}}}{\partial (-n)}\Big)\rd s
\\
=\int_{S} \log_{\delta}^{\zeta_j}\Delta L\rd m-\int_{\partial S}\log_{\delta}^{\zeta_j}\mathcal{N}(L)\,\rd s.
\end{multline*}
Substituting $\Delta(\log^{\zeta}_{\Delta})=2\pi\frac{\chi^{\zeta}_{\Delta}}{\|\chi^{\zeta}_{\Delta}\|_{L^1}}$ again we obtain that
\begin{align*} 
\Gamma(\log^{\zeta_j})=\frac{1}{8\pi}\int_{S}\Delta \log_{\delta}^{\zeta_j}\,\rd m+\frac{1}{8\pi}\int_{\C}(L^{S}-L^{S}_{\infty})\Delta \log_{\delta}^{\zeta_j} \rd m=\frac{1}{4}+\frac{1}{4}(L(\zeta_j)-L^{S}_{\infty})+\OO(N^{-1/2+\kappa}).
\end{align*}
We leave the $\Gamma(\log^{\zeta_j},\log^{\zeta_j})$ term as it is. So, the last term we consider is $\Gamma(\log^{\zeta_j},\log^{\zeta_k})$. The calculation of this term relies on the fact that the singularities are separated on a scale greater than the mesoscopic log-regularization. More explicitly,
\begin{align*} 
\int_{\C}\nabla(\log_{\delta}^{\zeta_k})^{S}\cdot\nabla(\log_{\Delta}^{\zeta_j}-\log_{\delta}^{\zeta_j})^{S}\,\rd m=\OO\big(\frac{1}{|\zeta_j-\zeta_k|}\big)\int_{{\rm B}(\zeta_j,\delta)}\frac{1}{|z-\zeta_j|}\,\rd m=\OO\big(\frac{\delta}{|\zeta_j-\zeta_k|}\big)=\OO(N^{-\kappa}).
\end{align*}
Thus
\begin{align*} 
\Gamma(\log^{\zeta_j},\log^{\zeta_k})&=\frac{1}{8\pi}\int_{S}\nabla\log_{\delta}^{\zeta_j}\cdot\nabla\log_{\Delta}^{\zeta_k}\,\rd m+\frac{1}{8\pi}\int_{S^{\rm c}}\nabla(\log_{\delta}^{\zeta_j})^{S}\cdot\nabla(\log_{\Delta}^{\zeta_k})^{S}\,\rd m+\OO(N^{-\kappa})
\\
&=-\frac{1}{8\pi}\int_{S}\log_{\delta}^{\zeta_j}\Delta\log_{\Delta}^{\zeta_k}\,\rd m+\frac{1}{8\pi}\int_{\partial S}\log_{\delta}^{\zeta_j}\frac{\partial \log_{\Delta}^{\zeta_k}}{\partial n}\rd s+\frac{1}{8\pi}\int_{\partial S}\log_{\delta}^{\zeta_j}\frac{\partial (\log_{\Delta}^{\zeta_k})^S}{\partial (-n)}\rd s+\OO(N^{-\kappa})
\\
&=-\frac{1}{4}\log|\zeta_j-\zeta_k|+\frac{1}{4}\int_{\partial S}\log^{\zeta_j}\rd\omega^{\infty}+\OO(N^{-\kappa})
\end{align*}
and the second integral in the right hand side is equal to $\mathfrak{s}$ by \eqref{eqn:lop_capacity}.
Therefore it remains only to show that $\mathscr{E}_N(F+tg,g_1)=\OO(N^{-\kappa/10})$. Let $h_t=F+tg$. Since the expression is linear in its second argument, it suffices to establish this bound for $\mathscr{E}_N(h_t,g)$ uniformly on $g\in\mathscr{A}_{C,\tau}$ with $\tau\in[N^{-1/2+\kappa},1]$, which is done in the next lemma.
\end{proof}

\begin{lemma}\label{lem:negligible_h_t,g}
$\mathscr{E}_N(h_t,g)$ term in the proof of Proposition \ref{proposition.singularities} is $\OO(N^{-\kappa/10})$.
\end{lemma}

\begin{proof} The proof follows the same general strategy as in Lemmas \ref{lem1}-\ref{lem3}, with only a few modifications that we now indicate. Most of the steps in the proofs of the lemmas remain valid without change: in particular, using the notation from Appendix \ref{sec:kernel}, we still have $\wt{\bf K}_N=\OO(N)$, exponential decay of $\wt{\bf K}_N(z,z)$ outside the unit disc, and $|\frac{\wt{\bf K}_N(z,z)}{N}-\rho_{V}(z)|\leq N^{-2\kappa+\epsilon}$ when $z\in S_{\rm in}(\delta_N)\setminus \cup_{j}{\rm B}(\zeta_j,\delta)$ by Lemma \ref{lem:1point} and Theorem \ref{prop:kernel}. The stability of the overwhelming probability events also continues to hold by Corollary \ref{cor:stability2}. The only points where the argument requires modification are those that invoke integration in a neighbourhood of the singularities, i.e. on $\cup_{j}{\rm B}(\zeta_j,\delta)$. In these steps, the local behaviour around the singularities must be handled separately, as we describe below.\\

\noindent\textit{Adjustments to the proof of Lemma \ref{lem1}:} The argument differs in \eqref{eqn:needchange} and \eqref{eqn:ring}, when integrals are evaluated on on $\cup_{j}{\rm B}(\zeta_j,\delta)$. For \eqref{eqn:needchange}, we split the integration domain $\cup_{j}{\rm B}(\zeta_j,\delta)$ into two regions as follows.
\begin{align}\label{eqn:changed1}
\left|\E_{h_t}\Big[\int_{S_{\rm in}(\delta_N)\cap{\rm B}(\zeta_j,3\delta_N)}\frac{\bar{\partial} p}{\Delta V}\partial f\rd\tilde\mu_N\Big]\right|\lesssim \int_{{\rm B}(\zeta_j,3\delta_N)}\frac{1}{\tau}\frac{1}{|z-\zeta_j|}\rd m(z)\lesssim \frac{\delta_N}{\tau}=\OO(N^{-\kappa+\epsilon})
\end{align}
where we used Lemma \ref{lem:1point}.(i). Moreover,
\begin{align}\label{eqn:changed2}
\left|\E_{h_t}\Big[\int_{S_{\rm in}(\delta_N)\cap{\rm B}(\zeta_j,\delta)\setminus {\rm B}(\zeta_j,3\delta_N)}\frac{\bar{\partial} p}{\Delta V}\partial f\rd\tilde\mu_N\Big]\right|\lesssim \int_{{\rm B}(\zeta_j,\delta)\setminus {\rm B}(\zeta_j,3\delta_N)}\frac{1}{\tau}\frac{1}{|z-\zeta_j|}\frac{1}{N|z-\zeta_j|^2}\rd m(z)\lesssim \frac{1}{\tau\delta_N N}=\OO(N^{-\kappa})
\end{align}
where we used Theorem \ref{prop:kernel}. 
For \eqref{eqn:ring}, the argument is simpler, relies only on ${\bf K}_N(z,z)=\OO(N)$:
\begin{equation*}
\int_{S_{\rm edge}(\delta_N)\cap{\rm B}(\zeta_j,\delta)} |\frac{\bar{\partial} p}{\Delta V}\partial f| \left|\frac{\wt {\bf K}_N}{N}-\rho_{V}\right|\rd m\lesssim \int_{S_{\rm edge}(\delta_N)\cap{\rm B}(\zeta_j,\delta)} \frac{1}{\tau |\zeta_j-z|}\rd m(z)=\OO(N^{-\kappa+\epsilon}).
\end{equation*}\\

\noindent\textit{Adjustments to the proof of Lemma \ref{lem2}:} The adjustments are identical to the previous case, and in fact simpler. We again re-evaluate the integrals over the neighbours of the singularities, while the rest of the argument remains unchanged. Moreover, the absence of a multiplicative term involving $f$ in the integrand further simplifies the analysis, since the estimations along $S_{\rm edge}(\delta_N)\cap{\rm B}(\zeta_j,\delta)$ proceed in the same manner as in the proof of Lemma \ref{lem2}, relying only on the bound $\wt{\bf K}_N = \OO(N)$. Hence, it suffices to estimate the integrals on $S_{\rm in}(\delta_N)\cap{\rm B}(\zeta_j,\delta)$, so, it suffices to bound integrations \eqref{eqn:changeinlem2} and \eqref{eqn:change2inlem2} on ${\rm B}(\zeta_j,\delta)$. 

We begin with the integral in \eqref{eqn:changeinlem2}. This can be done similarly to \eqref{eqn:changed1}-\eqref{eqn:changed2}:
\begin{align*} 
E_{h_t}\Big[\int_{{\rm B}(\zeta_j,3\delta_N)} \frac{\Delta g}{\Delta V} \,\rd\tilde{\mu}_{N}\Big]=\OO(\frac{\delta_N^2}{\tau^2}), \quad E_{h_t}\Big[\int_{S_{\rm in}(\delta_N)\cap{\rm B}(\zeta_j,\delta)\setminus {\rm B}(\zeta_j,3\delta_N)} \frac{\Delta g}{\Delta V} \,\rd\tilde{\mu}_{N}\Big]=\OO(\frac{\log N}{N\tau^2}).
\end{align*}
The bound on \eqref{eqn:change2inlem2} follows easier, directly by $\wt{\bf K}_N = \OO(N)$:
\begin{align*} 
E_{h_t}\Big[\int_{{\rm B}(\zeta_j,\delta)} \frac{\bar{\partial}p\partial L}{\Delta V/4} \,\rd\tilde{\mu}_{N}\Big]=\OO(\delta^2/\tau)
\end{align*}
and completes the proof.\\

\noindent\textit{Adjustments to the proof of Lemma \ref{lem3}:} For the first half of the proof concerning the estimation 
\begin{align}\label{eqn:r_evaluation}
\E_{h_t}\Big[N\iint_{z\neq w}\frac{r(z)-r(w)}{z-w}\tilde\mu_N(\rd z)\tilde\mu_N(\rd w)\Big]=\OO(N^{-\kappa+\epsilon}),
\end{align}
the required changes are, in essence, the same as those made in \eqref{eqn:changed1}-\eqref{eqn:changed2}. The bound for $I(r_k;S_{\rm edge}(\delta_N),\allowbreak S_{\rm edge}(\delta_N))$ remains the same as it relies only on $\wt{\bf K}_N(z,z)=\OO(N)$. The evaluation of $I(r_k;S_{\rm in}(\delta_N),S_{\rm edge}(\delta_N))$ requires a minor adjustment as follows. We split the integral $I(r_k;S_{\rm in}(\delta_N),S_{\rm edge}(\delta_N))$ into two regions, according to whether $z$ is close the the singularities or not. When $z\in\cup_j {\rm B}(\zeta_j,3\delta_N)$ we simply use $\wt{\bf K}_N(z,z)=\OO(N)$. Otherwise we apply Theorem \ref{prop:kernel} again. Hence the integral is bounded by,
\begin{multline*} 
\OO(N) \iint_{\substack{z\in S_{\rm in}(\delta_N)\setminus \cup_j {\rm B}(\zeta_j,3\delta_N),\\w\in S_{\rm edge}(\delta_N)}}\frac{|r_k(w)|}{|z-w|}\OO\left(\frac{\sup_{{\rm B}(z,2\delta_N)}|\nabla^2 f|}{N}+N^{-3/2}+\frac{|\wt{\bf{K}}_N(z,w)|^2}{N^2}\right)\rd m(z)\rd m(w)
\\
+\OO(N) \iint_{\substack{z\in \cup_j {\rm B}(\zeta_j,3\delta_N),\\w\in S_{\rm edge}(\delta_N)}}\frac{|r_k(w)|}{|z-w|}\rd m(z)\rd m(w).
\end{multline*}
The first double integral is estimated exactly as before, cf. \eqref{eqn:another_change}. For the second double integral we use the bounds $|r_k(w)|\leq\delta_N/\tau_k^2$ together with $|z-w|\geq N^{-1/2+\kappa}$ on the domain of integration. Thus,
\begin{align*} 
\OO(N) \iint_{\substack{z\in \cup_j {\rm B}(\zeta_j,3\delta_N),\\w\in S_{\rm edge}(\delta_N)}}\frac{|r_k(w)|}{|z-w|}\rd m(z)\rd m(w)=\OO(N\frac{\delta_N}{\tau_k^2}\frac{1}{N^{-1/2+\kappa}}) \iint_{\substack{z\in \cup_j {\rm B}(\zeta_j,3\delta_N),\\w\in S_{\rm edge}(\delta_N)\cap\supp(r_k)}}\rd m(z)\rd m(w)
\\
=\OO\big(N\frac{\delta_N}{\tau_k^2}\frac{1}{N^{-1/2+\kappa}}(\delta_N\tau_k)\delta_N^2\big)=\OO(N^{-\kappa}).
\end{align*}
This completes the necessary modifications for \eqref{eqn:r_evaluation}.

What remains in the proof are just the adjustments required to evaluate the estimates for $h_1$, $h_2$, and $(h_3-r)$; these are detailed below. The treatment of $E_2(\eta,h)$ in \eqref{eqn:E_1} is unaffected: the bounds for centered linear statistics (Proposition \ref{prop:exp_tail}) and transfer of the overwhelming probability events to the biased measure remain valid (Corollary \ref{cor:stability2}).

For $E_1(\eta,h)$, the three contributions from $\partial h(z)$, $\bar\partial h(z) \frac{\bar z - \bar w}{z-w}$, and $\OO(|z-w|)$ behave as follows. (i) The $\OO(|z-w|)$ term: \eqref{eqn:E_1case1} relies only on $\wt{\bf K}_N(z,z)=\OO(N)$ and exponential decay, hence unchanged. (ii) The $\partial h(z)$ term: \eqref{eqn:E_1case2} depends solely on the estimates for centered linear statistics, also unchanged by the stability. (iii) The $\bar\partial h(z)\frac{\bar z-\bar w}{z-w}$ term: the first two cases \eqref{eqn:anis_app}, \eqref{eqn:E_1case3_2} remain valid by Lemma \ref{lem:anis}, stability and Lemma \ref{lem:1point}. The only adjustment is required in the third case \eqref{eqn:E_1case3_3}. As the calculations are more intricate here, we re-evaluate the full integral, not just near the singularities, for the reader’s convenience despite some repetition. With $\theta=N^{-1/2+\kappa/5}$, we decompose the left-hand side of \eqref{eqn:E_1case3_3} as,
\begin{multline*}
\frac{N}{2}\int_{S\setminus\cup_{j}{\rm B}(\zeta_j,\theta)}\bar\partial h(z)\int_{S_{\rm in}(\eta)\setminus \cup_{j}{\rm B}(\zeta_j,\theta/2)} \frac{\bar z - \bar w}{z-w}e^{-\frac{|z-w|^2}{2\eta^2}}\big(\rho_{V}(w)+\OO(\frac{1}{\theta^2N}\log N)\big)\rd m(w)\rd \mu_{V}(z)
\\
\frac{N}{2}\int_{S\setminus\cup_{j}{\rm B}(\zeta_j,\theta)}\bar\partial h(z)\int_{\cup_{j}{\rm B}(\zeta_j,\theta/2)} \frac{\bar z - \bar w}{z-w}e^{-\frac{|z-w|^2}{2\eta^2}}\frac{\wt{\bf K}_N(w,w)}{N}\rd m(w)\rd \mu_{V}(z)
\\
+\frac{N}{2}\int_{\cup_{j}{\rm B}(\zeta_j,\theta)}\bar\partial h(z)\int_{S_{\rm in}(\eta)} \frac{\bar z - \bar w}{z-w}e^{-\frac{|z-w|^2}{2\eta^2}}\frac{\wt{\bf K}_N(w,w)}{N}\rd m(w)\rd \mu_{V}(z)
\\
+\frac{N}{2}\int_{S}\bar\partial h(z)\int_{S_{\rm in}(\eta)^{\rm c}} \frac{\bar z - \bar w}{z-w}e^{-\frac{|z-w|^2}{2\eta^2}}\frac{\wt{\bf K}_N(w,w)}{N}\rd m(w)\rd \mu_{V}(z).
\end{multline*}
(i) In the first term, the contribution from $\OO(\frac{1}{\theta^2N}\log N)$ can bounded similarly by $N \tau^2\frac{1}{\tau^2}\eta^2 \frac{\log N}{\theta^2N}=\oo(N^{-\kappa/5})$ where $\tau^2$ is the support size of $h$, $\frac{1}{\tau^2}$ is the maximum of $\bar{\partial}h$, and $\eta^2$ is the integral of $e^{\frac{-|z-w|^2}{2\eta^2}}$ on $w$. (ii) The second term is subpolynomial due to the exponential term, as $|z-w|>\theta/2\gg \eta$. (iii) The third term is controlled by $N\theta^2\frac{1}{\tau^2}\eta^2=\OO(N^{-\kappa})$ using $\wt{\bf K}_N(w,w)/N=\OO(1)$, where $\theta^2$ stands for the integration area over $z$, $\frac{1}{\tau^2}$ for maximum of $\bar\partial h$ and $\eta^2$ for the integral over $w$. (iv) Lastly, the fourth term term is bounded exactly as in \eqref{eqn:E_1case3_3}, using $\wt{\bf K}_N(w,w)/N=\OO(1)$. Thus, only the following term remains:
\begin{align*} 
\frac{N}{2}\int_{S\setminus\cup_{j}{\rm B}(\zeta_j,\theta)}&\bar\partial h(z)\int_{S_{\rm in}(\eta)} \frac{\bar z - \bar w}{z-w}e^{-\frac{|z-w|^2}{2\eta^2}}\rho_{V}(w)\rd m(w)\rd \mu_{V}(z)
\\
=&\,\frac{N}{2}\int_{S_{\rm in}(2\eta\log N)\setminus\cup_{j}{\rm B}(\zeta_j,\theta)}\bar\partial h(z)\int_{{\rm B}(z,\eta\log N)} \frac{\bar z - \bar w}{z-w}e^{-\frac{|z-w|^2}{2\eta^2}}\big(\rho_{V}(w)-\rho_{V}(z)\big)\rd m(w)\rd \mu_{V}(z)
\\
&+\frac{N}{2}\int_{S_{\rm in}(2\eta\log N)\setminus\cup_{j}{\rm B}(\zeta_j,\theta)}\bar\partial h(z)\int_{S_{\rm in}(\eta)\setminus{\rm B}(z,\eta\log N)} \frac{\bar z - \bar w}{z-w}e^{-\frac{|z-w|^2}{2\eta^2}}\rho_{V}(w)\rd m(w)\rd \mu_{V}(z)
\\
&+\frac{N}{2}\int_{S_{\rm edge}(2\eta\log N)\setminus\cup_{j}{\rm B}(\zeta_j,\theta)}\bar\partial h(z)\int_{S_{\rm in}(\eta)} \frac{\bar z - \bar w}{z-w}e^{-\frac{|z-w|^2}{2\eta^2}}\rho_{V}(w)\rd m(w)\rd \mu_{V}(z).
\end{align*} 
Analogously to the previous calculations, the first integral on the right-hand side is $\OO(N\eta^3\log N)$, the second is subpolynomial, and the third is $\OO(N(\tau\eta\log N)\frac{1}{\tau^2}\eta^2)$; each of these contributions is negligible.
\end{proof}

\subsection{Decoupling and comparison.}\label{subsec:local_log}\ We discuss two important results concerning submicroscopically regularized local log-singularities, both of which rely on Fredholm determinant theory: the decoupling and the Ginibre comparison.

The first proposition is an analogue of \cite[Proposition 3.4]{BouFal2025}. Thanks to the exponential decay of the determinantal kernel, Theorem \ref{thm:off_diag_decay}, the decoupling holds with much more relaxed conditions for 2d Coulomb gases. Moreover, because we have a single-time problem, the correlation kernel is already self-adjoint,  which makes the proof considerably easier compared to \cite{BouFal2025}. Here we give the details for the sake of completeness.

\begin{lemma}\label{lem:decoupling}
Assume that $C>0$, $m\in\N$, $\alpha>0$, $\kappa\in[0,1/2)$, $\Delta=N^{-1/2-\alpha}$, $\delta=N^{-1/2+\kappa}/3$, $\gamma_1,\dots,\gamma_m\in[0,C]$ and $\zeta_1,\dots,\zeta_m\in S_{\rm in}(N^{-1/2+\kappa})$ satisfy the separation condition $\min_{1\leq i\neq j\leq m}|\zeta_i-\zeta_j|>N^{-1+\kappa}$. Then
\begin{align*} 
\mathbb{E}\Big[\prod_{j=1}^{m} e^{\sum_{i=1}^{N}\gamma_j(\log^{\zeta_j}_{\Delta}-\log_{\delta}^{\zeta_j})(z_i)}\Big]=\prod_{j=1}^{m}\mathbb{E}\Big[e^{\sum_{i=1}^{N}\gamma_j(\log^{\zeta_j}_{\Delta}-\log_{\delta}^{\zeta_j})(z_i)}\Big](1+\OO(e^{-N^{\kappa/5}})).
\end{align*}
\end{lemma}

\begin{proof}
Define
\begin{align*} 
h_j=e^{\gamma_j(\log^{\zeta_j}_{\Delta}-\log_{\delta}^{\zeta_j})},\ g_j=\sqrt{1-h_j},\ h=\prod_{j=1}^{m}h_j,\ g=\sqrt{1-h}
\end{align*}
and let $\epsilon=\min_{j\in[\![m]\!]}\min_{z\in\C}h_j(z)$. Consider the operators $\mathcal{K},\mathcal{K}_1,\dots,\mathcal{K}_m$ acting on ${\rm L}^2(\mathbb{C})$,
\begin{align*}
(\mathcal{K}f)(z)=\int_{\mathbb{C}} g(z) {\bf K}_N(z,w)g(w) f(w)\rd m(w), 
\quad (\mathcal{K}_jf)(z)=\int_{\mathbb{C}} g_j(z) {\bf K}_N(z,w)g_j(w) f(w)\rd m(w)
\end{align*}
where ${\bf K}_N$ is the determinantal kernel, recall \eqref{eqn:Gin_kernel_defn}. Note that these are self-adjoint operators because ${\bf K}_N(z,w)=\overline{{\bf K}_N(w,z)}$ and $g$ and $g_j$'s are real valued:
\begin{align*}
\langle \phi,\mathcal{K}\psi\rangle=\iint \bar\phi(z)g(z)K(z,w)g(w)\psi(w)\rd m(z)\rd m(w)=\iint \bar g(w)\overline{K(w,z)}\bar g(z)\bar\phi(z)\psi(w)\rd m(z)\rd m(w)=\langle \mathcal{K}\phi,\psi\rangle.
\end{align*}
As ${\bf K}_N$ is a projection-type sum, we see that these are also positive semi-definite operators. From \cite[Theorem 2.12]{Sim2005}, we conclude that these are trace-class operators.  Application of the Cauchy-Binet formula (see e.g. \cite[Proposition 2.11]{Joh2006}),
\begin{align*} 
\mathbb{E}\big[\prod_{i=1}^{N}h(z_i)\big]=\det(\Id-\mathcal{K}), \quad \mathbb{E}\big[\prod_{i=1}^{N}h_j(z_i)\big]=\det(\Id-\mathcal{K}_j)
\end{align*}
where the right-hand sides are Fredholm determinants. Using the property $\det(\Id+A)\det(\Id+B)=\det(\Id+A+B+AB)$ for the Fredholm determinants we obtain
\begin{align*} 
\prod_{j=1}^{m}\mathbb{E}\big[\prod_{i=1}^{N}h_j(z_i)\big]=\prod_{j=1}^{m}\det(\Id-\mathcal{K}_j)=\det(\Id-\sum_{j=1}^{m}\mathcal{K}_j)
\end{align*}
where we have used the fact that each $g_j$ is supported on a disjoint set which yields $\mathcal{K}_i\mathcal{K}_j$ to be the zero operator for every $i\neq j$. For convenience we define, $\wt{\mathcal{K}}=\sum_{j=1}^{m}\mathcal{K}_j$.

From \cite[Theorem 3.7]{Sim2005}, we know that for every constant $\xi\in\R$, the Fredholm determinants $\det(\Id-\xi\mathcal{K})$ and $\det(\Id-\xi\wt{\mathcal{K}})$ can be expressed in terms of the eigenvalues $\{\lambda_k\}_{k=1}^{\infty}$, $\{\wt{\lambda_k}\}_{k=1}^{\infty}$ of the operators $\mathcal{K}$ and $\wt{\mathcal{K}}$ respectively, as follows:
\begin{gather*} 
\mathbb{E}\left[\prod_{n=1}^{N}\left(1-\xi (1-h(z_n))\right)\right]=\det(\Id-\xi\mathcal{K})=\prod_{k=1}^{\infty}(1-\xi\lambda_k), 
\\
 \prod_{j=1}^{m}\mathbb{E}\left[\prod_{n=1}^{N}\left(1-\xi (1-h_j(z_n))\right)\right]=\det(\Id-\xi\wt{\mathcal{K}})=\prod_{k=1}^{\infty}(1-\xi\wt{\lambda_k}).
\end{gather*}
Note that $(1-h)$ and $(1-h_j)$ lives inside the interval $[0,1-\epsilon]$. Therefore, the left-hand sides are strictly positive when $\xi<\frac{1}{1-\epsilon}$. This leads to the restriction on the spectra $\sigma(\mathcal{K}),\sigma(\wt{\mathcal{K}})\subset[0,1-\epsilon]$, when the right-hand side expressions are considered.

Following the same steps with the third step of \cite[Proof of Proposition 3.4]{BouFal2025}, for any $n\geq 2$, we can write
\begin{align}\label{eqn:fredholm_diff} 
\big|\log \det(\Id-\mathcal{K})-\log\det(\Id-\wt{\mathcal{K}})\big|\leq |\tr\mathcal{K}-\tr\wt{\mathcal{K}}|+\sum_{k=2}^{n}\big|\frac{\tr(\mathcal{K}^k)-\tr(\wt{\mathcal{K}}^k)}{k}\big|+\frac{1}{n\epsilon}(\|\mathcal{K}\|^2_{\rm{HS}}+\|\wt{\mathcal{K}}\|^2_{\rm{HS}}).
\end{align}
By \cite[Theorem 2.12]{Sim2005}, the traces can be expressed as the integrals along the diagonal, i.e.
\begin{align*} 
\tr\mathcal{K}=\int_{\C}g(z){\bf K}_N(z,z)g(z)\rd m(z)=\sum_{j=1}^{m}\int_{\C}g_j(z){\bf K}_N(z,z)g_j(z)\rd m(z)=\tr\wt{\mathcal{K}}
\end{align*}
since the supports of $g_j$'s are disjoint. On the other hand, for every $k\geq 2$
\begin{multline*} 
\left|\frac{\tr(\mathcal{K}^k)-\tr(\wt{\mathcal{K}}^k)}{k}\right|=\frac{1}{k}\left|\int_0^1\frac{\rd}{\rd t}{\rm Tr}(({t\mathcal{K}+(1-t)\wt{\mathcal{K}}})^k)\rd t \right|= \left|\int_0^1{\rm Tr}((\mathcal{K}-\wt{\mathcal{K}})({t\mathcal{K}+(1-t)\wt{\mathcal{K}}})^{k-1})\rd t \right|
\\
\leq \|\mathcal{K}-\wt{\mathcal{K}}\|_{\rm HS}\max_{t\in[0,1]}\|({t\mathcal{K}+(1-t)\wt{\mathcal{K}}})^{k-1}\|_{\rm HS}\leq \|\mathcal{K}-\wt{\mathcal{K}}\|_{\rm HS}\max_{t\in[0,1]}\left(\|{t\mathcal{K}+(1-t)\wt{\mathcal{K}}}\|_{\rm HS}\|{t\mathcal{K}+(1-t)\wt{\mathcal{K}}}\|^{k-2}\right)
\end{multline*}
where we have used $\tr(AB)\leq \|A\|_{\rm HS}\|B\|_{\rm HS}$ and $\|AB\|_{\rm HS}\leq \|A\|_{\rm HS}\|B\|$ for Hilbert-Schmidt operators ($\|\cdot\|$ stands for the operator norm).  As $\|{t\mathcal{K}+(1-t)\wt{\mathcal{K}}}\|\leq t\|{\mathcal{K}\|+(1-t)\|\wt{\mathcal{K}}}\|\leq 1-\epsilon$,  we have
\begin{align*} 
\left|\frac{\tr(\mathcal{K}^k)-\tr(\wt{\mathcal{K}}^k)}{k}\right|\leq \|\mathcal{K}-\wt{\mathcal{K}}\|_{\rm HS}\max(\|\mathcal{K}\|_{\rm HS},\|\wt{\mathcal{K}}\|_{\rm HS}).
\end{align*}
Finally, the following Hilbert-Schmidt norm bounds are straightforward, given the uniform estimations on the kernel, Theorem \ref{thm:off_diag_decay} and Lemma \ref{lem:1point}.(i):
\begin{align*} 
\|\mathcal{K}-\wt{\mathcal{K}}\|_{\rm HS}=\OO( e^{-N^{\kappa/2}}),\  \|\mathcal{K}\|_{\rm HS}=\OO(N^{\kappa}),\ \|\wt{\mathcal{K}}\|_{\rm HS}=\OO(N^{\kappa}). 
\end{align*}
Substituting these, we get 
\begin{align*} 
\big|\log \det(\Id-\mathcal{K})-\log\det(\Id-\wt{\mathcal{K}})\big|\lesssim n e^{-N^{\kappa/2}}N^{3\kappa/2}+\frac{1}{n\epsilon}N^{2\kappa}.
\end{align*}
Choosing $n=e^{N^{\kappa/4}}$ concludes the proof.
\end{proof}

The next observation follows from calculations closely analogous to those in \cite[Proposition 2.3]{Keles2025}, with one key modification. We do not estimate the Hilbert-Schmidt norm of the operator difference, since the angular term prevents the determinantal kernel for a general potential from being approximated by the Ginibre kernel. Instead, we use that these angular contributions vanish in the trace calculations, which is all that we require.

We denote the expectation for quadratic potential $V(z)=|z|^2$ case, i.e.  the Ginibre case, by $\E^{\rm Gin}$.

\begin{lemma}\label{lem:V_Gin_comp} Given $C>2$. Then, uniformly in $\gamma\in[0,C]$, $\alpha,\kappa>0$ satisfying $\alpha+\kappa<\frac{1}{20C}$, and $\zeta\in S_{\rm in}(N^{-1/2+\kappa})$, denoting $\delta=N^{-1/2+\kappa}/3$, $\Delta=N^{-1/2-\alpha}$ we have
\begin{align}\label{eqn:shift_kernel} 
\mathbb{E}\big[e^{\gamma\tr(\log^{\zeta}_{\Delta}-\log^{\zeta}_{\delta})}\big]=\mathbb{E}^{\rm Gin}\big[e^{\gamma\tr(\log_{\breve\Delta}-\log_{\breve\delta})}\big](1+\OO(N^{-1/8}))
\end{align}
where $\breve\Delta=\Delta\sqrt{\frac{\Delta V(\zeta)}{4}}$ and $\breve\delta=\delta\sqrt{\frac{\Delta V(\zeta)}{4}}$.
\end{lemma}

\begin{proof} Let $\theta=\sqrt{\frac{4}{\Delta V(\zeta)}}$ and $\breve V(z)=V(\theta z+\zeta)$. By change of variables, using the fact that $\log_{\epsilon}(cz)=\log_{\epsilon/c}z+\log c$, we have 
\begin{align*} 
\mathbb{E}^{V}\big[e^{\gamma\sum_{i=1}^{N}(\log^{\zeta}_{\Delta}-\log^{\zeta}_{\delta})(z_i)}\big]=\mathbb{E}^{\breve V}\big[e^{\gamma\sum_{i=1}^{N}(\log_{\Delta}-\log_{\delta})(\theta z_i)}\big]=\mathbb{E}^{\breve V}\big[e^{\gamma\sum_{i=1}^{N}(\log_{\breve\Delta}-\log_{\breve\delta})(z_i)}\big]
\end{align*}
and this shift applies directly to the droplet $S$ and the determinantal kernels. More explicitly, the droplet for the potential $\breve{V}$ is $\breve{S}=\{z\in C: \theta z+\zeta\in S\}$ and the determinantal kernel for $\breve{V}$, denoted by $\breve{\bf{K}}_N(z,w)$, is
\begin{align*} 
\breve{\mathbf{K}}_N(z,w)=\mathbf{K}_N(\theta z+\zeta,\theta w+\zeta)\,\theta^2.
\end{align*}
Using the Taylor expansion we obtain that for any $|z|,|w|<\delta$:
\begin{gather*} 
V(\theta z+\zeta)-2V(\theta z+\zeta,\overline{\theta w+\zeta})+V(\theta w+\zeta)=|z|^2-2z\bar{w}+|w|^2-2\ii F(z,w)+\OO(\delta^3),
\\
\textnormal{where}\ F(z,w)=\im\big(\partial V(\zeta)\theta(z-w)+\frac{\partial^2 V(\zeta)}{2}\theta^2(z^2-w^2)\big).
\end{gather*}
 Theorems \ref{prop:kernel} and \ref{thm:off_diag_decay} imply that for any $z,w\in\{z\in S:\dist(z,\partial S)>\delta\}$ with $|z-w|<2\breve\delta$:
\begin{align}\label{eqn:local_kern_est}
\breve{\mathbf{K}}_N(z,w)=\frac{N}{\pi}e^{-\frac{N}{2}(|z|^2-2z\bar{w}+|w|^2)+\ii NF(z,w)}+\OO(N^2\delta^3).
\end{align} 
Let $h=e^{\gamma(\log_{\breve\Delta}-\log_{\breve\delta})},\ g=\sqrt{1-h}$. Similar to the proof of Proposition \ref{lem:decoupling}, we define integral operators $\mathcal{K}_1(z,w)=g(z)\breve{\mathbf{K}}_{N}(z,w)g(w)$ and $\mathcal{K}_2(z,w)=g(z)\mathbf{K}^{\rm Gin}_{N}(z,w)g(w)$ where the Ginibre kernel is $\mathbf{K}^{\rm Gin}_{N}(z,w)=\frac{N}{\pi}e^{-\frac{N}{2}|z|^2-2z\bar{w}+|w|^2}$. In order to prove \eqref{eqn:shift_kernel}, it suffices to achieve an $\OO(N^{-1/8})$ bound to the following expression, as in \eqref{eqn:fredholm_diff}:
\begin{align}\label{eqn:fredholm_diff2}
\sum_{k=1}^{n}\big|\frac{\tr(\mathcal{K}_1^k)-\tr(\mathcal{K}_2^k)}{k}\big|+\frac{1}{n\epsilon}(\|\mathcal{K}_1\|^2_{\rm{HS}}+\|\mathcal{K}_2\|^2_{\rm{HS}})
\end{align}
for some $n$ to be chosen later, and $\epsilon=e^{\gamma(\log_{\Delta}0-\log_{\delta}0)}\asymp N^{-\gamma(\alpha+\kappa)}$. The rough bounds $|\breve{{\bf K}}_N(z,w)|=\OO(N)$ and $|{\bf K}^{\rm Gin}_N(z,w)|=\OO(N)$ give $\|\mathcal{K}_i\|_{\rm HS}=\OO(N^{\kappa})$ for $i=1,2$. On the other hand, note that 
\begin{align*} 
\tr(\mathcal{K}_i^k)=\int_{\C^k}\mathcal{K}_i(z_1,z_2)\mathcal{K}_i(z_2,z_3)\cdots\mathcal{K}_i(z_{k},z_1)\rd m(z_1)\dots\rd m(z_k)
\end{align*}
for $i=1,2$. For convenience, let $z_{k+1}=z_1$. Plugging the estimation \eqref{eqn:local_kern_est}, due to $\sum_{\ell=1}^{k}F(z_\ell,z_{\ell+1})=0$, we get
\begin{align*} 
\tr(\mathcal{K}_1^k)&=\int_{{\rm B}(0,\breve\delta)}\prod_{\ell=1}^{k}\big(\mathcal{K}_2(z_\ell,z_{\ell+1})+\OO(N^2\delta^3)\big)\rd m(z_1)\dots\rd m(z_k)
\\
&=\tr(\mathcal{K}_2^k)+\OO\big(\sum_{\ell=1}^{k}{k\choose\ell}(N^2\delta^5)^\ell\big)=\tr(\mathcal{K}_2^k)+k\OO(N^2\delta^5)
\end{align*}
uniformly in $k\leq (N^{2}\delta^5)^{-1}$, where we used the fact that $\int_{\C}Ne^{-\frac{N}{2}|z-w|^2}\rd m(z)=2\pi$. Substituting this into \eqref{eqn:fredholm_diff2} yields the following bound:
\begin{align*} 
\Big|\log\mathbb{E}\big[e^{\gamma\tr(\log^{\zeta}_{\Delta}-\log^{\zeta}_{\delta})}\big]-\log\mathbb{E}\big[e^{\gamma\tr(\log_{\Delta}-\log_{\delta})}\big]\Big|\leq n\OO(N^2\delta^5)+\frac{1}{n}N^{\gamma(\alpha+\kappa)+2\kappa}
\end{align*}
for $n\leq(N^{2}\delta^5)^{-1}$. Choosing $n=N^{1/4}$ completes the proof.
\end{proof}

\subsection{Determinant of the Ginibre matrix.}\ In this section, we derive the asymptotic formula for a single logarithmic singularity at the origin in the Ginibre case. 

\begin{lemma}\label{lemma:stab_reg_log}
Fix $C>0$. For all $\gamma\in[0,C]$, $\alpha>0$ and $\Delta=N^{-1/2-\alpha}$, we have
\begin{align*}
\mathbb{E}^{\rm Gin}\big[e^{\gamma \tr(\log_{\Delta})}\big]= e^{-N\gamma/2}N^{\frac{\gamma^2}{8}}\frac{(2\pi)^{\frac{\gamma}{4}}}{{\rm G}\left(1+\frac{\gamma}{2}\right)}\big(1+\OO(N^{-\min(\alpha,1)})\big)
\end{align*}
for all sufficiently large $N$.
\end{lemma}

\begin{proof}
By Kostlan's theorem \cite{Kos1992}, $\{|z_1|^2,\dots,|z_N|^2\}\stackrel{(d)}{=}\{\Gamma_1/N,\dots,\Gamma_N/N\}$ where $z_i$ are eigenvalues of Ginibre ensemble and $\Gamma_i$ are independent, with $\Gamma_i\sim\mathrm{Gamma}(i,1)$. Defining $f(x)=e^{\gamma\log_{\Delta}\sqrt{|x|}}$, this leads
\begin{align*} 
\mathbb{E}^{\rm Gin}\big[e^{\gamma\tr(\log_{\Delta})}\big]=\prod_{j=1}^{N}\mathbb{E}\big[f\big(\Gamma_j/N\big)\big].
\end{align*}
Writing the probability density of the gamma random variable we get
\begin{align*} 
\mathbb{E}^{\rm Gin}\big[f\big(\Gamma_j/N\big)\big]=\frac{1}{(j-1)!}\int_{0}^{\infty} \big(x/N\big)^{\gamma/2}x^{j-1}e^{-x}\rd x-\frac{1}{(j-1)!}\int_{0}^{\Delta^2} \Big(\big(x/N\big)^{\gamma/2}-f(x/N)\Big)x^{j-1}e^{-x}\rd x
\\
=\frac{1}{N^{\gamma/2}}\frac{\Gamma(j+\gamma/2)}{\Gamma(j)}+\OO\Big(\frac{\Delta^{2j+2\gamma}}{j!}\Big)=\frac{1}{N^{\gamma/2}}\frac{\Gamma(j+\gamma/2)}{\Gamma(j)}\Big(1+\OO(N^{-\alpha(j+1)})\Big).
\end{align*}
Using Barnes ${\rm G}$-function asymptotics (e.g. see \cite[(A.6)]{Vor1987}):
\begin{align*}
\log {\rm G}(z+1)=z^2\left(\frac{\log z}{2}-\frac{3}{4}\right)+\frac{z}{2}\log 2\pi-\frac{\log z}{12}+B+\OO_{|z|\to\infty}\left(\frac{1}{z}\right)
\end{align*}
for some constant $B$, we obtain
\begin{align*} 
\mathbb{E}^{\rm Gin}\big[e^{\gamma\tr(\log_{\Delta})}\big]=(1+\OO(N^{-\alpha}))N^{-N\gamma/2}\prod_{j=1}^{N}\frac{\Gamma(j+\gamma/2)}{\Gamma(j)}=(1+\OO(N^{-\alpha}))N^{-N\gamma/2}\frac{{\rm G}(N+1+\frac{\gamma}{2})}{{\rm G}(N+1){\rm G}(1+\frac{\gamma}{2})}
\\
=e^{-N\gamma/2}N^{\frac{\gamma^2}{8}}\frac{(2\pi)^{\frac{\gamma}{4}}}{{\rm G}\left(1+\frac{\gamma}{2}\right)}\big(1+\OO(N^{-\min(\alpha,1)})\big)
\end{align*}
concluding the proof.
\end{proof}

For later reference, we also record the formula non-regularized case which improves the error term in \cite[Theorem 1]{WebWong2019}. The proof follows identically.

\begin{lemma}\label{lem:centeredsing}
Let $C>0$. Uniformly in $\gamma\in[0,C]$, we have
\begin{align*}
\E^{\rm Gin}\Big[\prod_{k=1}^N|z_k|^\gamma\Big]
=N^{-N\frac{\gamma}{2}}\frac{{\rm G}(N+1+\frac{\gamma}{2})}{{\rm G}(N+1){\rm G}(1+\frac{\gamma}{2})}
=e^{-N\frac{\gamma}{2}}N^{\frac{\gamma^2}{8}}\frac{(2\pi)^{\frac{\gamma}{4}}}{{\rm G}\left(1+\frac{\gamma}{2}\right)}\left(1+\OO\left(\frac{1}{N}\right)\right).
\end{align*} 
\end{lemma}

\subsection{Stability for log-biases.}\label{subsec:stability2}\ 
Here, we derive some estimates on linear statistics of log-singular test functions regularized at submicroscopic scales under general external potential and state the related stability result, analogous to the Subsection \ref{app:stability}.

 Let $\alpha,\kappa>0$ and define $\delta=N^{-1/2+\kappa}/3$, $\Delta=N^{-1/2-\alpha}$ and $\zeta\in S_{\rm in}(N^{-1/2+\kappa})$. Proposition \ref{prop:Laplace_S_kappa} yields $\E[|X_{\log_{\delta}^{\zeta}}|]\leq (\log N)^{7}$. Moreover, using the kernel estimate in Theorem \ref{prop:kernel}, we get
\begin{align*} 
|X_{\log_{\Delta}^{\zeta}}-X_{\log_{\delta}^{\zeta}}|=\left|\int (\log_{\Delta}^{\zeta}z-\log_{\delta}^{\zeta}z)({\bf{K}}_N(z,z)-N\rho_{V}(z))\rd m(z) \right|=\oo(N^{-1}).
\end{align*}
Hence $\E[|X_{\log_{\Delta}^{\zeta}}|]\leq (\log N)^{7}$ for all sufficiently large $N$. For the Laplace transform bound, we leverage the Ginibre comparison for the local contribution. Combining Proposition \ref{prop:Laplace_S_kappa}, Lemma \ref{lemma:stab_reg_log}, and Lemma \ref{lem:V_Gin_comp}, we obtain the uniform Laplace bound for sub-microscopically regularized log-singularities in the bulk. More explicitly, a Cauchy-Schwarz inequality gives the following corollary.
\begin{corollary}
For every $C>1$, $\kappa\in(0,\frac{1}{100C})$, there exist $N_0$ such that for every $0<\alpha<\frac{1}{100C}$, $\Delta=N^{-1/2-\alpha}$ and $\zeta\in\{z\in S:\dist(z,\partial S)>N^{-1/2+\kappa}\}$, defining, we have
\begin{align*} 
\log\E\big[e^{t\,X_{\log_{\Delta}^{\zeta}}}\big]\leq (\log N)^7
\end{align*}
for all $t\in[0,C]$ and $N\geq N_0$. Moreover, $\E\big[|X_{\log_{\Delta}^{\zeta}}|\big]\leq (\log N)^7$.
\end{corollary}

An analogous statement to Corollary \ref{cor:stability} for biases involving submicroscopically regularized logarithms is as follows.

\begin{corollary}[Stability for log-biases]\label{cor:stability2}
Fix $C>1$, $\kappa,\alpha\in(0,\frac{1}{100C})$ and a non-negative integer $m$. Let $\gamma_1,\dots,\gamma_m\in[0,C]$ and $\zeta_1,\dots,\zeta_m\in S_{\rm in}(N^{-1/2+\kappa})$. Set
\begin{align*} 
f=h+\sum_{j=1}^{m}\gamma_j\log_{\Delta}^{\zeta_j},\quad h\in\mathscr{S}_{3,C,\kappa},\ \Delta=N^{-1/2-\alpha}.
\end{align*}
If an event $A$ satisfies $\mathbb{P}(A)\geq1-e^{-(\log N)^D}$ for some $D\geq 10$, then $\mathbb{P}_{f}(A)\geq 1-e^{-(\log N)^D/2}$.
\end{corollary}

\subsection{Proof of Theorem \ref{Theorem.2singularities}.}\ We now prove Theorem \ref{Theorem.2singularities}, by combining Proposition \ref{proposition.singularities}, Lemmas \ref{lem:decoupling}, \ref{lem:V_Gin_comp}, \ref{lem:centeredsing}.

\begin{proof}[Proof of Theorem \ref{Theorem.2singularities}] By Proposition \ref{prop:submic_log_reg}, it suffices to prove the statement of the Theorem for exponentials of submicroscopically regularized logarithmic singularities, instead of root-type singularities. Let $\alpha=\kappa$ and take the submicroscopic scale $\Delta=N^{-1/2-\alpha}$ in the regularizations.

Combining the asymptotics in Proposition \ref{proposition.singularities} with the results from Lemmas \ref{lem:decoupling} and \ref{lem:V_Gin_comp} we obtain that
\begin{multline*} 
\E\Big[e^{\sum_{i=1}^Nf(z_i)}\prod_{j=1}^{m}e^{\sum_{i=1}^N\gamma_j\log_{\Delta}(z_i-\zeta_j)}\Big]
=e^{N\int f+\sum_{j}\gamma_j\log_{\delta}^{\zeta_j}\rd\mu_{V}}
e^{\frac{1}{8\pi}\big(\int_{\C}|\nabla f^S|^2\rd m
+\mathds{1}_{S}\Delta f
+L^{S}\Delta f\rd m\big)
}
\\
\times\prod_{j=1}^{m}e^{\frac{\gamma_j}{2}\big(\int_{\partial S}f\rd\omega^{\infty}-f(\zeta_j)\big)} 
\prod_{j=1}^{m}e^{\frac{\gamma_j^2}{8\pi}\int_{\C}\nabla(\log_{\delta}^{\zeta_j})^{S}\cdot\nabla(2\log_{\Delta}^{\zeta_j}-\log_{\delta}^{\zeta_j})^{S}\,\rd m
+\frac{\gamma_j}{4}(1+L(\zeta_j)-L^{S}(\infty))}
\\
\times\prod_{j=1}^{m}\E^{\rm Gin}\left[e^{\sum_{i=1}^N\gamma_j(\log_{\Delta_{j}}-\log_{\delta_{j}})(z_i)}\right]
\prod_{j\neq k}|\zeta_j-\zeta_k|^{-\frac{\gamma_i\gamma_j}{4}}
e^{\frac{\gamma_i\gamma_j}{4}\mathfrak{s}} \left(1+\OO(N^{-\kappa/10})\right)
\end{multline*}
where $\Delta_j=\Delta\sqrt{\frac{\Delta V(\zeta_j)}{4}}$ and $\delta_j=\delta\sqrt{\frac{\Delta V(\zeta_j)}{4}}$. On the other hand, applying Proposition \ref{proposition.singularities} for quadratic potential with a single singularity, by the same way as above, we obtain that for any $j=1,\dots,m$
\begin{align*} 
\E^{\rm Gin}\Big[\prod_{i=1}^N|z_i|^{\gamma_j}\Big]=&\,e^{\frac{N}{\pi}\gamma_j\int_{\D} \log_{\delta_j}\rd m} e^{\gamma_j^2\int_{\C}\nabla(\log_{\delta_j})^{\D}\cdot\nabla(2\log_{\Delta_j}-\log_{\delta_j})^{\D}\,\rd m
+\frac{\gamma_j}{4}} 
\\
&\times\E^{\rm Gin}\left[e^{\sum_{i=1}^N\gamma_j(\log_{\Delta_{j}}-\log_{\delta_{j}})(z_i)}\right] \left(1+\OO(N^{-\kappa/10})\right).
\end{align*}
Combining the last two asymptotics, we can use the following simplifications:
\begin{multline*} 
\int \log_{\delta}^{\zeta_j}\rd\mu_{V}-\frac{N}{\pi}\int_{\D}\log_{\delta_j}\rd m=\int \log^{\zeta_j}\rd\mu_{V}-\frac{N}{\pi}\int_{\D}\log\rd m+\int (\log_{\delta}^{\zeta_j}-\log^{\zeta_j})(\frac{\Delta V}{4\pi}-\frac{\Delta V(\zeta_j)}{4\pi})\rd m
\\
=\int \log^{\zeta_j}\rd\mu_{V}-\frac{N}{\pi}\int_{\D}\log\rd m+\OO(\delta^3)
\end{multline*}
where $\log$ in the integrands stands for $\log|\cdot|$ and
\begin{align*} 
\int_{\C}\nabla(\log_{\delta}^{\zeta_j})^{S}&\cdot\nabla(2\log_{\Delta}^{\zeta_j}-\log_{\delta}^{\zeta_j})^{S}\,\rd m-\int_{\C}\nabla(\log_{\delta_j})^{\D}\cdot\nabla(2\log_{\Delta_j}-\log_{\delta_j})^{\D}\,\rd m
\\
&=\int_{\C}\nabla(\log_{\delta}^{\zeta_j})^{S}\cdot\nabla(\log_{\Delta}^{\zeta_j})^{S}\,\rd m-\int_{\C}\nabla(\log_{\delta_j})^{\D}\cdot\nabla(\log_{\Delta_j})^{\D}\,\rd m
\\
&=\int_{S}\nabla\log_{\delta}^{\zeta_j}\cdot\nabla\log_{\Delta}^{\zeta_j}\,\rd m+\int_{S^{\rm c}}\nabla(\log_{\delta}^{\zeta_j})^{S}\cdot\nabla(\log_{\Delta}^{\zeta_j})^{S}\,\rd m-\int_{\D}\nabla\log_{\delta_j}\cdot\nabla\log_{\Delta_j}\,\rd m
\\
&=-\int_{S}\log_{\delta}^{\zeta_j}\Delta\log_{\Delta}^{\zeta_j}\rd m+\int_{\partial S} \log^{\zeta_j}\big(\frac{\partial\log^{\zeta_j}}{\partial n}-\frac{\partial(\log^{\zeta_j})^S|_{S^{\rm c}}}{\partial n}\big)\,\rd s+\int_{\D}\log_{\delta_j}\Delta\log_{\Delta_j}\,\rd m
\\
&=-2\pi\log_{\delta}(0)+2\pi\log_{\delta_j}(0)+\int_{\partial S} \log^{\zeta_j}\mathcal{N}(\log^{\zeta_j})\,\rd s +\OO(\Delta/\delta)
\\
&=\pi L(\zeta_j)+2\pi\mathfrak{s} +\OO(N^{-\kappa-\alpha}).
\end{align*}
Together with Lemma \ref{lem:centeredsing}, these complete the proof.
\end{proof}

\section{Convergence to the Gaussian multiplicative chaos}\label{sec:GMC}\ 
A collection of sufficient conditions for convergence to a GMC measure throughout the entire $L^1$-phase has been provided in \cite{CFLW2021}. We first state this result  and then verify the conditions in our setting using the main Theorem \ref{Theorem.2singularities}. This leads to the proof of Theorem \ref{GMCforGinibre}.

\subsection{GMC convergence in the $L^1$-phase.}\  
Let $U$ be a simply connected, open, bounded subset of $\R^{d}$ for some $d\geq1$. Let $X$ be a log-correlated Gaussian field on $U$ with symmetric positive semi-definite covariance kernel 
\begin{align*} 
K(z,w)=\log\frac{1}{|z-w|}+g(z,w)
\end{align*}
for a continuous function $g\in L^{2}(U\times U)$ that is bounded from above. Let $(X_N)_{N\geq 1}$ be a sequence of random functions, defined on probability spaces $\Omega_N$,
\begin{align*} 
X_N:\Omega_N\to\big\{f:U\to\R\ \big|\ f\in L^1(U),\ \sup_{z\in U}f(z)<\infty,\ f\ \textnormal{is upper semi-continuous}\big\}
\end{align*}
satisfying the following properties: for each $\beta>0$, $z\in U$ and $\varphi\in\mathscr{C}_{\rm c}^{\infty}(U)$,
\begin{align*} 
\E [e^{\beta X_N(z)}]<\infty,\quad \E [e^{\int_{U} X_N(z)\varphi(z)\rd m(z)}]<\infty, \quad \int_{U} X_N(z)\varphi(z)\rd m(z)\xrightarrow[N\to\infty]{\rm (d)}\int_{U} X(z)\varphi(z)\rd m(z).
\end{align*}
Denote the mollifications of $X_N$ and $X$ by 
\begin{align*} 
X_N^{(\epsilon)}=X_N\ast\frac{\chi_{\epsilon}}{\|\chi_{\epsilon}\|_{L^1}},\quad X^{(\epsilon)}=X\ast\frac{\chi_{\epsilon}}{\|\chi_{\epsilon}\|_{L^1}}
\end{align*}
for every $\epsilon>0$ and let $X_N^{(0)}=X_N$, $X^{(0)}=X$ for convenience. We also denote the multiplicative chaos measure generated by the Gaussian field $X$ on $U$ by $\mu^{\beta}$ which can be formally written as
\begin{align*} 
\rd\mu^{\beta}(z)=\frac{e^{\beta X(z)}}{\E[ e^{\beta X(z)}]}.
\end{align*}

\begin{theorem}[{\hspace{-.01mm}\cite[Proposition 2.13]{CFLW2021}}]\label{thm:gmc_conds} Let $X_N$ and $X$ be as defined above. Assume that for a $\beta\in(0,\sqrt{2d})$, there exists a sequence $\epsilon_N=\epsilon_N(\beta)$ converging to $0$ as $N\to\infty$ satisfying the following properties:
\begin{enumerate}[(i)]
    \item For any fixed $\epsilon,\epsilon'\geq 0$,
\begin{align} \label{eqn:cond1}
\lim_{N\to\infty}\frac{\E[e^{\beta X_N^{(\epsilon)}(z) + \beta X_N^{(\epsilon')}(w)}]}{\E[e^{\beta X_N^{(\epsilon)}(z)}]\E[e^{\beta X_N^{(\epsilon')}(w)}]}=e^{\beta^2\E[X^{(\epsilon)}(z)X^{(\epsilon')}(w)]}
\end{align}
for all $z\neq w$ in $U$ and the convergence is uniform for all $(z,w)$ in any fixed compact subset of $\{(u,v)\in U^2:u\neq v\}$.
    \item For any fixed $\epsilon\geq 0$, $\epsilon'>0$, and compact set $K\subset U$, there exists a positive constant $C=C(\beta,\epsilon,\epsilon',K)$ such that
\begin{align} \label{eqn:cond2}
\sup_{N\in\N}\sup_{z,w\in K}\frac{\E[e^{\beta X_N^{(\epsilon)}(z) + \beta X_N^{(\epsilon')}(w)}]}{\E[e^{\beta X_N^{(\epsilon)}(z)}]\E[e^{\beta X_N^{(\epsilon')}(w)}]}\leq C.
\end{align} 
	\item For any fixed $\lambda\in\R$, and compact set $K\subset U$, there exists a positive constant $C=C(\beta,\lambda,K)$ such that
\begin{align} \label{eqn:cond3}
\frac{\E[e^{\beta X_N(z)+\lambda X_N^{(\epsilon)}(z)}]}{\E[e^{\beta X_N(z)}]}\leq C\epsilon^{-\lambda\beta-\frac{\lambda^2}{2}}
\end{align}
for all $z\in K$, $\epsilon\geq\epsilon_N$, and $N\in\N$.
    \item For any fixed compact set $K\subset U$, there exists a positive constant $C=C(\beta,K)$ such that
\begin{align} \label{eqn:cond4}
\frac{\E[e^{\beta X_N(z) + \beta X_N(w)}]}{\E[e^{\beta X_N(z)}]\E[e^{\beta X_N(w)}]}\leq C|z-w|^{-\beta^2}
\end{align}
for all $(z,w)\in K^2\cap\{(u,v)\in U^2:\, |u-w|\geq\epsilon_N\}$ and $N\in\N$.
    \item For any fixed $\rho>0$, $\epsilon,\epsilon'\geq0$, $n\in\N$, $\boldsymbol{\lambda}\in\R^n$, and compact set $K\subset U$,
\begin{align} 
\frac{\E[e^{\beta X_{N}^{(\epsilon)}(z)+\beta X_N^{(\epsilon')}(w)}e^{\sum_{k=1}^{n}\lambda_k X_{N}^{(\eta_k)}(u_k)}]}{\E[e^{\beta X_{N}^{(\epsilon)}(z)+\beta X_N^{(\epsilon')}(w)}]}=&\,(1+\oo_{N\to\infty}(1))\,e^{\frac{1}{2}\E[(\sum_{k=1}^{n}\lambda_kX^{(\eta_k)}(u_k))^2]}\nonumber
\\
&\times e^{\sum_{k=1}^{n}\beta \lambda_k\E[X^{(\epsilon)}(z)X^{(\eta_k)}(u_k)+X^{(\epsilon')}(z)X^{(\eta_k)}(u_k)]} \label{eqn:cond5}
\end{align}    
uniformly for all $(z,w)\in K^2\cap\{(u,v)\in U^2:|u-v|\geq\rho\}$, $\mathbf{u}\in K^{n}$, and $\boldsymbol{\eta}\in(\epsilon_N,1]^{n}$. The implicit constant may depend on $\beta,\rho,\epsilon,\epsilon',n$, $\boldsymbol{\lambda}$, and $K$.
	\item For any fixed $\lambda\in\R$ and compact set $K\subset U$, there exists a positive constant $C=C(\beta,\lambda,K)$ such that
\begin{align}\label{eqn:cond6}
\limsup_{N\to\infty}\frac{\E[e^{\beta X_N^{(\epsilon)}(z)+\beta X_N^{(\epsilon')}(w)+\lambda X_N^{(\eta)}(w)}]}{\E[e^{\beta X_N^{(\epsilon)}(z)+\beta X_N^{(\epsilon')}(w)}]}\leq C e^{\frac{\lambda^2}{2}\E[X^{(\eta)}(w)^2]+\lambda\beta\,\E[X^{(\epsilon)}(z)X^{(\eta)}(w)+X^{(\epsilon')}(w)X^{(\eta)}(w)]}
\end{align}
for all $\eta\geq\epsilon'\geq\epsilon$ with $\epsilon'>0$ and $x,y\in K$.
	\item For any fixed $\lambda\in\R$ and compact set $K\in U$, there exists a positive constant $C=C(\beta,\lambda,K)$ such that
\begin{align} \label{eqn:cond7}
\frac{\E[e^{\beta X_N(z)+\beta X_N(w)+\lambda X_{N}^{(\eta)}(w)}]}{\E[e^{\beta X_N(z)+\beta X_N(w)}]}\leq Ce^{\frac{\lambda^2}{2}\E[X^{(\eta)}(w)^2]+\lambda\beta\,\E[X(z)X^{\eta}(w)+X(w)X^{\eta}(w)]}
\end{align}
for all $\eta\geq\epsilon_N$ and $(z,w)\in K^2\cap\{(u,v):\,|u-v|\geq\epsilon_N\}$.
	\item There exists a small parameter $0<\theta<\min(\frac{(d-\frac{\beta^2}{2})^2}{4d},\frac{\beta^2}{8})$ such that for any compact set $K\subset U$ there is a positive constant $C=C(\beta,K)$ satisfying
\begin{align} \label{eqn:cond8}
\iint_{\substack{z,w\in K,\\|z-w|\leq\epsilon_N}}\frac{\E[e^{\beta X_N(z)+\beta X_N(w)}]}{\E[e^{\beta X_N(z)}]\E[e^{\beta X_N(w)}]}\rd m(z)\rd m(w)\leq C\epsilon_N^{-\max(\beta^2-d,0)-\theta}.
\end{align}	
\end{enumerate}
Then
\begin{align*} 
\frac{e^{\beta X_N(z)}}{\E[e^{\beta X_N(z)}]}\rd m(z)\xrightarrow[N\to\infty]{}\rd\mu^{\beta}(z)
\end{align*}
where the convergence holds in distribution with respect to the weak topology of measures, i.e., for every bounded continuous function $\varphi:U\to\R$,
\begin{align*} 
\int_{U} \varphi(z)\frac{e^{\beta X_N(z)}}{\E[e^{\beta X_N(z)}]}\rd m(z)\xrightarrow[N\to\infty]{\rm (d)}\int_{U} \varphi(z)\,\rd\mu^{\beta}(z).
\end{align*}
\end{theorem}

In \cite{CFLW2021}, condition \eqref{eqn:cond7} is stated uniformly for all $(z,w)\in K^2$. However, tracing through their proof shows that the requirement is only necessary when $z$ and $w$ are separated at scale $\epsilon_N$. Consequently, it is sufficient to assume \eqref{eqn:cond7} holds uniformly on the restricted set as given above. In addition, we correct a typo in the denominator of this condition: in \cite{CFLW2021}, it was written as $\E[e^{\beta X(z)+\beta X(w)}]$.

Moreover, for any condition in the theorem, once a compact set $K \subset U$ is fixed, it is enough to verify the equation only for mollification parameters $\epsilon$, $\epsilon'$, $\eta$, $\eta_k$ smaller than $\dist(K,\partial U)$.

\subsection{Proof of Theorem \ref{GMCforGinibre}.}\ To obtain a log-correlated limit without any additional scaling, we modify \eqref{eqn:defn_X_N} and for convenience in the calculations below, define
\begin{align*} 
X_N(z)=\sqrt{2}\cdot\Big(\sum_{i=1}^{N}\log|z-z_i|-\E\Big[\sum_{i=1}^{N}\log|z-z_i|\Big]\Big)
\end{align*}
which yields
\begin{align*} 
X_N^{(\epsilon)}(z)=\sqrt{2}\cdot\Big(\sum_{i=1}^{N}\log_{\epsilon}(z-z_i)-\E\Big[\sum_{i=1}^{N}\log_{\epsilon}(z-z_i)\Big]\Big).
\end{align*}
Denote the limit field by $X$. It follows easily by Theorem \ref{Theorem.2singularities} that $X$ is a centered log-correlated field with covariance structure
\begin{align*} 
\E \big[X(z)X(w)\big]=\log\frac{1}{|z-w|}+\mathfrak{s}=\frac{1}{2\pi}\int_{\C}\nabla(\log^{z})^S\cdot\nabla(\log^{w})^S\rd m,\quad z,w\in S,
\end{align*}
where $\mathfrak{s}$ is as in \eqref{eqn:capacity_defn}. Then
\begin{align*} 
\E \big[X^{(\epsilon)}(z)X^{(\epsilon')}(w)\big]=-\int_{\C}\log_{\epsilon}^{z-w}\frac{\chi_{\epsilon'}}{\|\chi_{\epsilon'}\|_{L^{1}}}\rd m+\mathfrak{s},\quad z,w\in S
\end{align*}
for all $\epsilon,\epsilon'\geq0$. Recall that
\begin{align*} 
-\int_{\C}\log_{\epsilon}^{z-w}\frac{\chi_{\epsilon'}}{\|\chi_{\epsilon'}\|_{L^{1}}}\rd m+\mathfrak{s}=\frac{1}{2\pi}\int_{\C}\nabla(\log_{\epsilon}^{z})^S\cdot\nabla(\log_{\epsilon'}^{w})^S\rd m
\end{align*}
if $\epsilon<\dist(z,\partial S)$ and $\epsilon'<\dist(w,\partial S)$.

Moreover, the following formulation of Theorem \ref{Theorem.2singularities} will be useful in the calculations:
\begin{align*}
\E\left[e^{\tr f+\sum_{j=1}^m \gamma_j\tr\log^{\zeta_j}}\right]=
\E\big[e^{\tr f}\big]\prod_{j=1}^m \E[e^{\sum_{i=1}^{N}\gamma_j\log^{\zeta_j}z_i}]
\prod_{j=1}^{m}e^{\frac{1}{4\pi}\int_{\C}\nabla f^S\cdot\nabla(\gamma_j\log^{\zeta_j})^S\rd m} 
\\
\times\prod_{j< k}e^{\frac{1}{4\pi}\int_{\C}\nabla (\gamma_j\log^{\zeta_j})^S\cdot\nabla(\gamma_k\log^{\zeta_k})^S\rd m} 
\big(1+\OO(N^{-\kappa/10})\big)
\end{align*}
with
\begin{align*} 
\E\big[e^{\tr f-\E[\tr f]}\big]=e^{\frac{1}{8\pi}\int_{\C}|\nabla f^S|^2\rd m}(1+\OO(N^{-\kappa/10})\big).
\end{align*}

\begin{proof}[Proof of Theorem \ref{GMCforGinibre}] Fix a $\beta\in(0,2)$. We use Theorem \ref{thm:gmc_conds}, choosing $\epsilon_N=N^{-1/2+\kappa}$ for fixed $\kappa=\min\big(\frac{1}{1000},\frac{1-\frac{\beta^2}{2}}{5}\big)$. We now verify every condition in Theorem \ref{thm:gmc_conds} one-by-one, using the asymptotic formula in Theorem \ref{Theorem.2singularities}.\\

\noindent\textit{Condition \eqref{eqn:cond1}.} Fix an arbitrary compact subset of $\{(u,v)\in U^2:u\neq v\}$. Uniformly for all $(z,w)$ in the compact subset,
\begin{align*} 
\frac{\E[e^{\beta X_N^{(\epsilon)}(z) + \beta X_N^{(\epsilon')}(w)}]}{\E[e^{\beta X_N^{(\epsilon)}(z)}]\E[e^{\beta X_N^{(\epsilon')}(w)}]}=e^{\beta^2\frac{1}{2\pi}\int_{\C}\nabla(\log_{\epsilon}^{z})^S\cdot\nabla(\log_{\epsilon'}^{w})^S\rd m}(1+\OO(N^{-\kappa/10}))=e^{\beta^2\E[X^{(\epsilon)}(z)X^{(\epsilon')}(w)]}(1+\OO(N^{-\kappa/10})).
\end{align*}

\noindent\textit{Condition \eqref{eqn:cond2}.} Given $\epsilon'>0$ and compact set $K$ in $U$, uniformly for all $z,w\in K$,
\begin{align*} 
\frac{\E[e^{\beta X_N^{(\epsilon)}(z) + \beta X_N^{(\epsilon')}(w)}]}{\E[e^{\beta X_N^{(\epsilon)}(z)}]\E[e^{\beta X_N^{(\epsilon')}(w)}]}=e^{-\beta^2\int\log_{\epsilon'}^{w-z}\frac{\chi_{\epsilon}}{\|\chi_{\epsilon}\|_{L^{1}}}\rd m+\beta^2\mathfrak{s}}\OO(1)=e^{\beta^2\OO(1+\log\epsilon')}.
\end{align*}

\noindent\textit{Condition \eqref{eqn:cond3}.} Fix a compact $K\subset U$. Uniformly for all $z\in K$,
\begin{multline*} 
\frac{\E[e^{\beta X_N(z)+\lambda X_N^{(\epsilon)}(z)}]}{\E[e^{\beta X_N(z)}]}=e^{-\lambda\beta\frac{1}{2\pi}\int_{\C}\nabla(\log^{z})^S\cdot\nabla(\log_{\epsilon}^{z})^S\rd m}e^{\frac{\lambda^2}{2}\frac{1}{2\pi}\int_{\C}|\nabla(\log_{\epsilon}^{z})^S|^2\rd m}\OO(1)
\\
=e^{-\lambda\beta\log_{\epsilon}0+\lambda\beta \mathfrak{s}}e^{-\frac{\lambda^2}{2}\int\log_{\epsilon}\frac{\chi_{\epsilon}}{\|\chi_{\epsilon}\|_{L^{1}}}\rd m+\frac{\lambda^2}{2}\mathfrak{s}}\OO(1)=\epsilon^{-\lambda\beta-\frac{\lambda^2}{2}}\OO_{\lambda,\beta}(1).
\end{multline*}

\noindent\textit{Condition \eqref{eqn:cond4}.} Uniformly for $z$ and $w$ in the given set,
\begin{align*} 
\frac{\E[e^{\beta X_N(z) + \beta X_N(w)}]}{\E[e^{\beta X_N(z)}]\E[e^{\beta X_N(w)}]}=|z-w|^{-\beta^2}e^{\beta^2\mathfrak{s}}\OO(1).
\end{align*}

\noindent\textit{Condition \eqref{eqn:cond5}.} Uniformly for all $z$, $w$, $\mathbf{u}$, and $\boldsymbol{\eta}$ satisfying the given conditions:
\begin{align*} 
\frac{\E[e^{\beta X_{N}^{(\epsilon)}(z)+\beta X_N^{(\epsilon')}(w)}e^{\sum_{k=1}^{n}\lambda_k X_{N}^{(\eta_k)}(u_k)}]}{\E[e^{\beta X_{N}^{(\epsilon)}(z)+\beta X_N^{(\epsilon')}(w)}]}=&\,(1+\OO(N^{-\kappa/10}))e^{\frac{1}{4\pi}\int_{\C}\nabla|(\sum_{k=1}^{n}\lambda_k\log^{u_k}_{\eta_k})^S|^2\rd m} 
\\
&\times  e^{\sum_{k=1}^{n}\frac{1}{2\pi}\int_{\C}\nabla(\beta\log^{z}_{\epsilon})^S\cdot\nabla(\lambda_k\log^{u_k}_{\eta_k})^S\rd m}e^{\sum_{k=1}^{n}\frac{1}{2\pi}\int_{\C}\nabla(\beta\log^{w}_{\epsilon'})^S\cdot\nabla(\lambda_k\log^{u_k}_{\eta_k})^S\rd m}
\\
=&\,(1+\OO(N^{-\kappa/10}))e^{\frac{1}{2}\E[(\sum_{k=1}^{n}\lambda_kX^{(\eta_k)}(u_k))^2]}\nonumber
\\
&\times e^{\sum_{k=1}^{n}\beta \lambda_k\E[X^{(\epsilon)}(z)X^{(\eta_k)}(u_k)+X^{(\epsilon')}(z)X^{(\eta_k)}(u_k)]}.
\end{align*}

\noindent\textit{Condition \eqref{eqn:cond6}.} Fix $\eta\geq\epsilon'\geq\epsilon$ with $\epsilon'>0$ and $z,w\in K$. Then
\begin{align*} 
\frac{\E[e^{\beta X_N^{(\epsilon)}(z)+\beta X_N^{(\epsilon')}(w)+\lambda X_N^{(\eta)}(w)}]}{\E[e^{\beta X_N^{(\epsilon)}(z)+\beta X_N^{(\epsilon')}(w)}]}=\OO(1)\,e^{\frac{1}{4\pi}\int_{\C}\nabla|(\lambda\log^{w}_{\eta})^S|^2\rd m}e^{\frac{1}{2\pi}\int_{\C}\nabla(\beta\log^{z}_{\epsilon})^S\cdot\nabla(\lambda\log^{w}_{\eta})^S\rd m}e^{\frac{1}{2\pi}\int_{\C}\nabla(\beta\log^{w}_{\epsilon'})^S\cdot\nabla(\lambda\log^{w}_{\eta})^S\rd m}
\\
=\OO(1)\,e^{\frac{\lambda^2}{2}\E[X^{(\eta)}(w)^2]+\lambda\beta\,\E[X^{(\epsilon)}(z)X^{(\eta)}(w)+X^{(\epsilon')}(w)X^{(\eta)}(w)]}.
\end{align*}

\noindent\textit{Condition \eqref{eqn:cond7}}. The calculation is identical to that of condition \eqref{eqn:cond6}.

\vspace{3pt}\noindent\textit{Condition \eqref{eqn:cond8}}. Applying Cauchy-Schwarz to the numerator we obtain that the expression in \eqref{eqn:cond8} is bounded by
\begin{align*} 
\iint_{\substack{z,w\in K,\\|z-w|\leq\epsilon_N}}\frac{\E[e^{2\beta X_N(z)}]^{1/2}\E[e^{2\beta X_N(w)}]^{1/2}}{\E[e^{\beta X_N(z)}]\E[e^{\beta X_N(w)}]}\rd m(z)\rd m(w)=\iint_{\substack{z,w\in K,\\|z-w|\leq\epsilon_N}} \OO(N^{\frac{\beta^2}{2}})\rd m(z)\rd m(w)=\OO(N^{-1+\frac{\beta^2}{2}+2\kappa}).
\end{align*}	
By the choice of $\kappa\leq\frac{1-\frac{\beta^2}{2}}{5}$ the desired result follows easily.
\end{proof}

\appendix

\setcounter{equation}{0}
\setcounter{theorem}{0}
\renewcommand{\theequation}{\thesection.\arabic{equation}}
\renewcommand{\thetheorem}{\thesection.\arabic{theorem}}
\section{Kernel estimates}

\label{sec:kernel}

In this appendix,  we obtain kernel estimates by following  \cite[Section 3]{AmeHedMak2010} and \cite[Appendix]{AmeHedMak2015}, with certain modifications, in particular to cover some logarithmic singularities. We provide the details for completeness.

\subsection{Setup.}\ Let $P_0(z),P_1(z),\dots$ be the monic analytic polynomials of degrees $0,1,\dots$ respectively that are orthogonal with respect to the measure $e^{-NV(z)}\rd m(z)$, i.e., the inner product is given by $\langle f_1,f_2\rangle=\int f_1 \bar f_2 e^{-N V}$. Let $K_N$ be the reproducing kernel of the space of analytic polynomials of degree at most $N-1$ with norm induced from the inner product above; more explicitly,
\begin{align*} 
K_N(z,w)=\sum_{n=0}^{N-1}\frac{P_n(z)\overline{P_n(w)}}{\|P_n\|^2_{L^2(e^{-NV})}}
\end{align*}
which satisfies the reproducing property
\begin{align}\label{eqn:reprod}
\int_{\C} K_N(z,\xi)K_N(\xi,w) e^{-NV(\xi)} \rd m(\xi)=K_N(z,w).
\end{align}
When viewed as an integral operator, this is a projection from $L^2(e^{-NV})$ to $\mathscr{P}_N$, the space of analytic polynomials of degree at most $N-1$:
\begin{align*} 
\Pi_N f(z)=\int_{\C} K_N(z,w)f(w)e^{-NV(w)}\rd m(w)=\sum_{n=0}^{N-1} \langle f,P_n\rangle P_n(z)
\end{align*}
Recall from Subsection \ref{subsec:NotConv} that we assume $V$ is real-analytic in a neighborhood of the droplet. Hence we may extend $V(z,\bar{z})=V(z)$ to a function $V(\cdot,\cdot)$ that is complex-analytic in two variables in a neighbourhood of the diagonal $\{(z,\bar{z}):z\in\C\}\cap S$; see, for example, \cite[Section 2]{AmeHedMak2011}. For this extension we have
\begin{align*} 
V(\bar{z},\bar{w})=\overline{V(w,z)},\quad \partial_{1}^{n}\partial_{2}^{k}V(z,\bar{z})=\partial^{n}\bar\partial^{k}V(z)
\end{align*}
where $\partial_{1}$ and $\partial_{2}$ are partial derivatives with respect to the first and second coordinates. The first order approximation of $K_N(z,w)$ inside the droplet $S$ is given by
\begin{align*} 
K_N^{\#}(z,w)=\frac{N}{\pi}\partial_{1}\partial_{2}V(z,\bar{w})e^{NV(z,\bar{w})}.
\end{align*}
Embedding the measure into the kernels, we obtain the determinantal kernel ${\bf K}_N$ and we define its approximation ${\bf K}_N^{\#}$ similarly
\begin{align} \label{eqn:Gin_kernel_defn}
{\bf K}_N(z,w)=K_N(z,w)e^{-\frac{N}{2}(V(z)+V(w))}, \quad {\bf K}_N^{\#}(z,w)=K_N^{\#}(z,w)e^{-\frac{N}{2}(V(z)+V(w))}.
\end{align}

The same definitions applies to $\wt K_N$, $\wt {\bf K}_N$ where $V$ is replaced by $\wt V=V-\frac{f}{N}$ where $f$ is a smooth function to be determined later. We define the projection $\wt \Pi_N g(z)=\langle g,\wt K_{N}(\cdot,z)\rangle_{L^2(e^{-N\wt V})}$ similarly. The main result of this appendix is that for sufficiently regular functions $f$,
\begin{align*} 
\wt K_N^{\#}(z,w)=\frac{N}{\pi}\partial_{1}\partial_{2}V(z,\bar{w})e^{NV(z,\bar{w})-f_{w}(z)}, \quad  \wt{\bf  K}_N^{\#}(z,w)=\wt K^{\#}(z,w)e^{-\frac{N}{2}(\wt V(z)+\wt V(w))}
\end{align*}
approximates $\wt K_N$ and $\wt {\bf K}_N$ where $f_w(z)=f(w)+(z-w)\partial f(w)$. We also define the projection $\wt \Pi_N^{\#}$ for $\wt K_{N}^{\#}$ similarly to $\wt \Pi_N$, i.e. using the weight $e^{-N\wt V}$.

\subsection{Results.}\ 
In the remainder of this section, we fix an arbitrarily large constant $C>0$, small parameters $\kappa,\alpha>0$ and a non-negative integer $m$. Let $\gamma_1,\dots,\gamma_m\in[0,C]$ and $\zeta_1,\dots,\zeta_m\in S_{\rm in}(2\delta_N)$ satisfy the separation condition $\min_{1\leq i\neq j\leq m}|\zeta_i-\zeta_j|>2\delta_N$, $h\in\mathscr{S}_{3,C,\kappa}$ and $\Delta=N^{-1/2-\alpha}$. We then set 
\begin{align*}
f=h+\sum_{j=1}^{m}\gamma_j\log_{\Delta}^{\zeta_j}.
\end{align*}

\begin{theorem}\label{prop:kernel}
Let $f$ be defined as above and $\theta>0$ be a fixed constant. Then, uniformly in $z\in  S_{\rm in}(2\delta_N)$ satisfying $\inf_{{\rm B}(z,2\delta_N)}\Delta V>\theta$ with $\min_{j\in[\![m]\!]}|z-\zeta_j|>3\delta_N$ and any $w\in {\rm B}(z,\delta_N/2)$, we have
\begin{align*}
|\wt {\bf K}_N(z,w)-\wt{\bf  K}_N^{\#}(z,w)|=\OO\left(\sup_{{\rm B}(z,2 \delta_N)}|\nabla^2 f|+N^{-1/2}\right).
\end{align*}
\end{theorem}

Theorem \ref{prop:kernel} is an immediate consequence of the following two propositions. Here, in the rest of this section, for a generic kernel $K$ we often write $K_w(z)=K_N(z,w)$, omitting the $N$-dependence. We also fix a cutoff function $\chi_z$ such that $\chi_z=1$ in ${\rm B}(z,\delta_N)$, $\chi_z=0$ outside ${\rm B}(z,3\delta_N/2)$, and $\|\nabla \chi_z\|_{\infty}\lesssim \delta_N^{-1}$. 

\begin{proposition}\label{prop:calculus}
Assume $f$ satisfies the assumptions of Theorem \ref{prop:kernel}. Let $\theta>0$ be a fixed constant. Then, uniformly in $z\in S_{\rm in}(2\delta_N)$ satisfying $\inf_{{\rm B}(z,2\delta_N)}\Delta V>\theta$ with $\min_{j\in[\![m]\!]}|z-\zeta_j|>3\delta_N$ and any $w\in {\rm B}(z,\delta_N/2)$, we have
\begin{align*}
\wt K_N(z,w)-\wt\Pi_N\big(\wt K_w^{\#}\chi_{z}\big)(z)=\OO\left(\sup_{{\rm B}(z,2 \delta_N)}|\nabla^2 f|+N^{-1/2}\right) e^{\frac{N}{2}(\wt V(z)+\wt V(w))}.
\end{align*}
\end{proposition}

\begin{proposition}\label{prop:Hormander}
Assume $f$ satisfies the assumptions of Theorem \ref{prop:kernel}. Let $\theta>0$ be a fixed constant. Then, uniformly in $z\in S_{\rm in}(2\delta_N)$ satisfying $\inf_{{\rm B}(z,2\delta_N)}\Delta V>\theta$ with $\min_{j\in[\![m]\!]}|z-\zeta_j|>3\delta_N$ and any $w\in {\rm B}(z,\delta_N/2)$, we have
\begin{align*}
\wt K_N^{\#}(z,w)-\wt\Pi_N\big(\wt K_w^{\#}\chi_{z}\big)(z)=\OO\left(e^{-c (\log N)^4}\right)e^{\frac{N}{2}(\wt V(z)+\wt V(w))}.
\end{align*}
\end{proposition}

We also have the following lemma which is particularly useful when evaluating the kernel outside the droplet $S$.

\begin{lemma}\label{lem:1point}
Assume $f$ satisfies the assumptions of Theorem \ref{prop:kernel}. Then:
\begin{itemize}
\item[\textnormal{(i)}] For all $z\in\C$, $|\wt{\bf{K}}_N(z,z)|=\OO(N)$.
\item[\textnormal{(ii)}] For all $z\in\C$, $|\wt{\bf{K}}_N(z,z)|=N^{\OO(1)}e^{-N(\wt V(z)-\check V(z))}$.
\end{itemize}
\end{lemma}
Recall from the assumptions on $V$ (Subsection \ref{subsec:potential}) that there exists $\epsilon,c>0$ such that $V(z)-\check{V}(z)\geq c\dist(z,\partial S)^2$ on $S_{\rm edge}(\epsilon)\setminus S$; where the constant $c$ may change line to line.   Moreover, combining the growth condition on $V$ with $\check{V}(z)=2\log |z|+\OO_{|z|\to\infty}(1)$, we obtain a constant $c>0$ such that $V(z)-\check{V}(z)\geq c\log |z|$ for all $z\in{\rm B}(0,C)^{\rm c}$. Since the coincidence set $\{z\in\C:V(z)=\check{V}(z)\}$ coincides with $\bar{S}$, the difference $V-\check{V}$ is strictly positive on $S_{\rm out}(\epsilon)$. By continuity, there exists $c>0$ such that $V(z)-\check{V}(z)\geq c$ on ${\rm B}(0,C)\cap S_{\rm out}(\epsilon)$. Putting these estimates together, we conclude that there exist constants $c,C>0$ such that
\begin{align*} 
V(z)-\check{V}(z)\geq
\begin{cases}
c\dist(z,\partial S)^2, &z\in {\rm B}(0,C)\setminus S\\
c\log|z|, &z\in{\rm B}(0,C)^{\rm c}
\end{cases}.
\end{align*}

In addition to these results, whose proofs are provided in the following section, we quote without proof a theorem concerning the off-diagonal decay of the kernel from \cite{AmeHedMak2010}. The same estimate was independently established in a more general setting in \cite[Theorem 5.7]{Ber2016}.

\begin{theorem}[{\hspace{-.01mm}\cite[Corollary 8.2]{AmeHedMak2010}}]\label{thm:off_diag_decay} There exists $C>0$ and $c>0$ depending only on $V$ such that for any $z\in S$, denoting $r={\rm dist}(z,\partial S)/2$, we have
\begin{align*} 
|{\bf K}_N(z,w)|\leq CNe^{-c(\inf_{{\rm B}(z,r)}\Delta V)^{1/2}\sqrt{N}\min(r,|w-z|)}e^{-N(V(w)-\check{V}(w))}
\end{align*}
for all $w\in\C$.
\end{theorem}

\subsection{Proofs.}\ We begin with the proof of Lemma \ref{lem:1point}, followed by those of Propositions \ref{prop:calculus} and \ref{prop:Hormander}.

\begin{proof}[Proof of Lemma \ref{lem:1point}] Let $\delta>0$ and $u:\C\to\C$ be a function (which are possibly $N$-dependent). Define, for any $z\in\C$
\begin{align*} 
F(\xi)=|u(z+\delta\xi)|^2e^{-N\wt V(z+\delta\xi)+A(z)|\delta\xi|^2}
\end{align*}
where $A(z)=\max(0,\sup_{{\rm B}(z,\delta)}N\Delta \wt V)$. Easy to see that, if $u$ is analytic on ${\rm B}(z,\delta)$, then $\Delta(\log F)\geq 0$ for all $\xi\in \D$, i.e. $\log F$ is subharmonic, so is $F$. By the mean value inequality for the subharmonic functions we obtain
\begin{align}\label{eqn:F_subharm} 
|u(z)|^2e^{-N\wt V(z)}\leq \int_{\D}|u(z+\delta\xi)|^2e^{-N\wt V(z+\delta\xi)+A(z)|\delta\xi|^2}\rd m(\xi)\leq \frac{e^{A(z)\delta^2}}{\delta^2}\int_{{\rm B}(z,\delta)}|u(\xi)|^2e^{-N\wt V(\xi)}\rd m(\xi)
\end{align}
for every $z\in \C$. We fix an arbitrary $w\in\C$ and take $u(z)=\frac{\wt K_N(z,w)}{\sqrt{\wt K_N(w,w)}}$. Notice that by reproducing property of the kernel we have $\int_{\C}|u|^2e^{-N\wt V}=1$. Using the above inequality, we get
\begin{align*} 
\frac{|\wt K_N(z,w)|^2}{\wt K_N(w,w)}e^{-N\wt V(z)}\leq \frac{e^{A(z)\delta^2}}{\delta^2}.
\end{align*}
Notice that $A(z)=\OO(\max(0,N \sup_{{\rm B}(z,\delta)}\Delta V))$ as $\Delta(\log^{\zeta}_{\Delta})=2\pi\frac{\chi^{\zeta}_{\Delta}}{\|\chi^{\zeta}_{\Delta}\|_{L^1}}\geq 0$. Thus, taking $w=z$ and $\delta=N^{-1/2}$ if $\sup_{{\rm B}(z,1)}\Delta V\leq 0$ and $\delta=N^{-1/2}\min(1,(\sup_{{\rm B}(z,1)}\Delta V)^{-1/2})$ otherwise,  we have
\begin{align} \label{eqn:almost_O(N)}
|\wt{\bf{K}}_N(z,z)|=\OO(N)\max(1,\sup_{{\rm B}(z,1)}\Delta V)
\end{align}
which proves the estimation (i) uniformly in compact subsets of $\C$.

On the other hand, note that by the choice of $f$, there is a constant $c$ such that for all $z\in\C$, $e^{-N\wt V(z)}\geq e^{-NV(z)-c\log N}$ ($c$ may change line-by-line). Hence, taking $\delta=N^{-1/2-\alpha}$ we obtain
\begin{align*} 
\Big|\frac{\wt K_N(z,w)}{N^c\sqrt{\wt K_N(w,w)}}\Big|^2 e^{-NV(z)}\leq 1, \quad \textnormal{for all $z\in S$.}
\end{align*}
Thus, by \cite[Theorem III.2.1]{SaffTotik2024LogPot}, we get
\begin{align*} 
\Big|\frac{\wt K_N(z,w)}{N^c\sqrt{\wt K_N(w,w)}}\Big|^2 e^{-N\check{V}(z)}\leq 1, \quad \textnormal{for all $z\in\C$}
\end{align*}
(cf. \cite[Lemma 3.4]{AmeHedMak2010}). Taking $z=w$, this gives 
\begin{align} \label{eqn:exterior_bound}
|\wt K_N(z,z)|e^{-N\wt V(z)}\leq N^{c}e^{-N(\wt V(z)-\check{V}(z))}
\end{align}
completing the proof of (ii) in the lemma. Moreover, combining \eqref{eqn:almost_O(N)} and \eqref{eqn:exterior_bound} with the growth condition on $V$, it follows immediately that (i) holds uniformly in $z\in\C$.
\end{proof}
\vspace{1mm}

Before moving onto the proof of Proposition \ref{prop:calculus}, we start with the following quantitative strict analogue of \cite[Lemma A.2]{AmeHedMak2015}.
\begin{lemma}\label{analytic} Let $z\in S_{\rm in}(2\delta_N)$ with $\min_{j\in[\![m]\!]}|z-\zeta_j|>3\delta_N$ and $g$ be an analytic function on ${\rm B}(z,2\delta_N)$. Then, uniformly in $w\in {\rm B}(z,\delta_N/2)$,
\begin{align*} 
g(w)-\wt\Pi^{\#}_N\big(g\chi_{z}\big)(w)=\OO\left(\int
\Big(|(w-\xi)g\chi_{z}|+\frac{|g\bar\partial\chi_{z}|}{\delta_N}+|g\chi_{z} |\sup_{{\rm B}(z,2\delta_N)}|\nabla^2f| \Big) e^{N\re(V(w,\bar{\xi})-V(\xi))}\rd m(\xi)\right)
\end{align*}
where the implicit constant in the error term does not depend on $z$ and $g$.
\end{lemma}

\begin{proof}
Substituting the definitions, $\wt\Pi^{\#}_N\big(g\chi_{z}\big)(w)$ is equal to
\begin{align*}
\int_{\C} \wt K^{\#}_N(w,\xi)g(\xi)\chi_{z}(\xi)e^{-N\wt{V}(\xi)}\rd m(\xi)=\frac{N}{\pi}\int g(\xi)\partial_1\partial_2V(w,\bar{\xi}) e^{N(V(w,\bar{\xi})-V(\xi))}\chi_{z}(\xi) e^{-(w-\xi)\partial f(\xi)}\rd m(\xi).
\end{align*}
So, defining $a(\xi)=g(\xi)e^{N(V(w,\bar{\xi})-V(\xi))}$ and $b(\xi)=\chi_{z}(\xi) e^{-(w-\xi)\partial f(\xi)}$, the above expression can be written as
\begin{align*} 
\frac{ 1}{\pi}\int\frac{b(\xi)\bar\partial a(\xi)}{w-\xi}\rd m(\xi)+\OO\left(\int\Big|(w-\xi) g(\xi)e^{N(V(w,\bar{\xi})-V(\xi))}\chi_{z}(\xi)e^{-(w-\xi)\partial f(\xi)} \Big|\rd m(\xi)\right).
\end{align*}
By Green's formula applied to $a\cdot b$ (here we use $\bar\partial g=0$), the first term is equal to 
\begin{align*}
g(w)-\frac{ 1}{\pi}\int\frac{a(\xi)\bar\partial b(\xi)}{w-\xi}\rd m(\xi).
\end{align*}

To bound the errors, note $\sup_{{\rm B}(z,2\delta_N)}|\nabla f|\lesssim \delta_N^{-1}$, so that $e^{-(w-\xi)\partial f(\xi)}=\OO(1)$ when $b$ is non-zero, and $\bar \partial b(\xi)=\OO(|\bar\partial \chi_{z}(\xi)|+\chi_{z}(\xi)|w-\xi|\sup_{{\rm B}(z,2\delta_N)}|\nabla^2 f|)$.
Noting that $|w-\xi|\asymp\delta_N$ when $\bar{\partial}\chi_{z}$ is non-zero the result follows easily.
\end{proof}
\vspace{1mm}

\begin{proof}[Proof of Proposition \ref{prop:calculus}]
We apply the above lemma to $g(w)=\wt K_N(w,z)$ for an arbitrary $z\in S$ with $\min_{j\in[\![J]\!]}|z-\zeta_j|>3\delta_N$, and obtain that uniformly in $w\in \rm B(z,\delta_N/2)$,
\begin{align*} 
\wt K_N(w,z)-\wt\Pi^{\#}_N\big(\wt K_z\chi_{z}\big)(w)
=\OO\left(\int
\Big(|(w-\xi)g\chi_{z}|+\frac{|g\bar\partial\chi_{z}|}{\delta_N}+|g\chi_{z} |\sup_{{\rm B}(z,2\delta_N)}|\nabla^2f| \Big) e^{N\re(V(w,\bar{\xi})-V(\xi))}\rd m(\xi)\right).
\end{align*}
First, by the Cauchy-Schwarz inequality we have 
\begin{align*}
|\wt K_N(\xi,z)|\leq |\wt K_N(\xi,\xi)\wt K_N(z,z)|^{1/2}\leq C N e^{\frac{N}{2}(\wt V(\xi)+\wt V(z))}.
\end{align*}
The first inequality means the two-point function is non-negative, and the second is a consequence of Lemma \ref{lem:1point}.(i). Second, on the support of $\chi_z$, we have $f(\xi)=f(w)+\OO(1)$. Third, by Taylor expansion,
\begin{align}\label{eqn:taylorforQ} 
\re\left(-V(w,\bar{w})+2V(w,\bar{\xi})-V(\xi,\bar{\xi}))\right)=-|w-\xi|^2\frac{\Delta V(w)}{4}+\OO(|w-\xi|^3).
\end{align}
Hence, substituting these, we have proved that 
\begin{align*}
\wt K_N(w,z)=\wt\Pi^{\#}_N\big(\wt K_z\chi_{z}\big)(w)+\OO
\left(N^{-1/2}+\sup_{{\rm B}(z,2\delta_N)}|\nabla^2f|\right)e^{\frac{N}{2}(\wt V(z)+\wt V(w))}.
\end{align*}
Moreover, noting that $\chi_{z}$ and $\wt{V}$ are real-valued, taking conjugates of both sides completes the proof of Proposition \ref{prop:calculus}.
\end{proof}
\vspace{1mm}

\begin{proof}[Proof of Proposition \ref{prop:Hormander}]
First, viewing $\wt\Pi_N$ as a projection into the set of analytic polynomials with degree up to $N-1$ (denoted by $\mathscr{P}_N$), the function $\wt u=\wt K_w^{\#}\chi_{z}-\wt\Pi_N\big(\wt K_w^{\#}\chi_{z}\big)$ can be described as the $L^2(e^{-N\wt V})$-minimum solution of the following system:
\begin{align*} 
\begin{cases}
\bar{\partial}u=\bar{\partial}\big(\wt K_w^{\#}\chi_z\big),\\
u-\wt K_w^{\#}\chi_{z}\in\mathscr{P}_N.
\end{cases}
\end{align*}
Similarly, we also denote $u=\wt K_w^{\#}\chi_{z}-\Pi_N\big(\wt K_w^{\#}\chi_{z}\big)$ which is the $L^2(e^{-N V})$-minimum solution of the same system. Our goal is to bound $\wt u(z)$. Because $|z-w|<\delta_N$, $\wt u$ is analytic on $B(z,N^{-1/2})$. So, by Equation \eqref{eqn:F_subharm} we have
\begin{align}\label{eqn:u(z)bound} 
|\wt u(z)|^2e^{-N\wt V(z)}\leq N\|\wt u\|^2_{L^2(e^{-N\wt V})}\leq N\| u\|^2_{L^2(e^{-N\wt V})}
\end{align}
where in the second inequality we used the fact that $\wt u$ is the $L^2(e^{-N\wt V})$-minimum solution. Note that, due to the choice of function $f$ in $\wt V=V-\frac{f}{N}$, we have $f=\OO(\log N)$ on ${\rm B}(0,N)$. This gives
\begin{align*} 
\| u\cdot \mathds{1}_{{\rm B}(0,N)}\|^2_{L^2(e^{-N\wt V})}\leq N^c\| u\|^2_{L^2(e^{-N V})}
\end{align*}
for some constant $c>0$. For the remainder of the proof, $c$ will be used to denote arbitrary constants, possibly taking different values in different occurrences. On the other hand, on $|\xi|\geq N$,
\begin{align*} 
u(\xi)\leq\int \big|K_N(\xi,y)\wt K_N^{\#}(y,w)\chi_{z}(y)e^{-NV(y)}\big|\rd m(y)\leq e^{cN}|K_N(\xi,\xi)|^{1/2}
\end{align*}
where in the second inequality we have used Cauchy-Schwarz and Lemma \ref{lem:1point}.(i) to get $|K_N(\xi,y)|\leq |K_N(\xi,\xi)|^{1/2}|K_N(y,y)|^{1/2}\leq |K_N(\xi,\xi)|^{1/2}e^{cN}$ and a rough bound $|K_N^{\#}(y,w)|\leq e^{cN}$ in the integration region. Applying Lemma \ref{lem:1point}.(ii) gives $|K_N(\xi,\xi)|e^{-N\wt V(\xi)}\leq e^{-N \epsilon \log|\xi|}$ for a fixed positive constant $\epsilon$ as $|\xi|\gg 1$. Substituting this bound, we obtain
\begin{align*} 
\| u\cdot \mathds{1}_{{\rm B}(0,N)^{\rm c}}\|^2_{L^2(e^{-N\wt V})}\leq e^{cN}\int_{{\rm B}(0,N)^{\rm c}}|K_N(\xi,\xi)|e^{-N\wt V(\xi)}\rd m(\xi)\leq e^{-cN\log N}.
\end{align*}
Thus, we obtain
\begin{align}\label{eqn:uL2norm} 
\| u\|^2_{L^2(e^{-N\wt V})}\leq N^c\| u\|^2_{L^2(e^{-N V})}+e^{-cN\log N}.
\end{align}

Next, we define $\varphi(\xi)=\check{V}(\xi)+\frac{1}{N}\log(1+|\xi|^2)$ noting that it is strictly subharmonic outside $S$; and for any entire function, $\|g\|_{L^2(e^{-N\varphi})}<+\infty$ implies $g\in\mathscr{P}_N$. Indeed, the latter follows by using polar coordinates for $\int |g|^2e^{-N\varphi}$ and substituting Cauchy integral formula $\frac{\rd^n(|g|^2)}{\rd z^n}(0)\leq \frac{n!}{2\pi r^n}\int_{0}^{2\pi}|g(re^{\ii\theta})|^2\rd \theta$. We denote the $L^2(e^{-N\varphi})$-minimum solution of 
\begin{align*} 
\bar{\partial}v=\bar{\partial}\big(\wt K_w^{\#}\chi_z\big)
\end{align*}
by $v$. Due to a.e. subharmonicity of $\varphi$, applying H\"ormander's estimate \cite[(4.2.6)]{Hormander94} (also see \cite[Section 4.2]{AmeHedMak2010}) we obtain
\begin{align}\label{eqn:hormander} 
\|v\|^2_{L^2(e^{-NV})}\leq \|v\|^2_{L^2(e^{-N\varphi})}\lesssim \int |\bar{\partial}(\wt K_w^{\#}\chi_z)|^2\frac{e^{-N\varphi}}{N\Delta\varphi}\lesssim N^{-1}\|\bar{\partial}(\wt K_w^{\#}\chi_z)\|_{L^2(e^{-NV})}^2
\end{align}
where in the first inequality we used $\varphi\leq V+c$ and in the last inequality we used $\Delta\varphi\asymp 1$ and $V\leq\varphi$ on the support of $\chi_{z}$. 
On the other hand, by definition, $v-\wt K_w^{\#}\chi_z$ is an entire function and since it has a finite $L^2(e^{-N\varphi})$-norm, it must be an analytic polynomial in $\mathscr{P}_N$. This implies, by the $L^2(e^{-N V})$-minimality of $u$, 
\begin{align*} 
\|u\|^2_{L^2(e^{-NV})}\leq \|v\|^2_{L^2(e^{-NV})}.
\end{align*}
Combining this with \eqref{eqn:u(z)bound}, \eqref{eqn:uL2norm} and \eqref{eqn:hormander} we obtain
\begin{align} \label{eqn:tildeubound}
|\wt u(z)|e^{-\frac{N}{2}\wt V(z)}\leq N^c\|\bar{\partial}(\wt K_w^{\#}\chi_z)\|_{L^2(e^{-NV})}+e^{-cN\log N}
\end{align}

Moreover, when $|\xi-w|\asymp\delta_N$, by Taylor expansion \eqref{eqn:taylorforQ} we have
\begin{align*} 
e^{N\re V(\xi,\bar{w})}\leq e^{\frac{N}{2} V(\xi)+\frac{N}{2}V(w)-c(\log N)^4}.
\end{align*}
Substituting this estimate into the definition of $\wt{K}^{\#}_N(\xi,w)$, we obtain
\begin{align*} 
\|\bar{\partial}(\wt K_w^{\#}\chi_z)\|_{L^2(e^{-NV})}\lesssim \Big(\int|\bar{\partial}\chi_{z}(\xi)|^2 e^{-c(\log N)^4}\rd m(\xi)\Big)^{1/2}e^{\frac{N}{2}V(w)}\leq e^{-c(\log N)^4}e^{\frac{N}{2}V(w)}.
\end{align*}
Combined with \eqref{eqn:tildeubound}, this completes the proof.
\end{proof}

\setcounter{equation}{0}
\setcounter{theorem}{0}
\section{The harmonic measure and capacity}\label{app:neumann}

In this appendix, we justify the definitions of the harmonic measure seen from infinity \eqref{eqn:Neumann_log} and the logarithmic capacity \eqref{eqn:lop_capacity}, and prove several basic properties mentioned in the introduction.
\\

\noindent\textit{Well-definedness of \eqref{eqn:Neumann_log}.} We show that the function $\mathcal{N}(\log^{\zeta})$ does not depend on the choice of $\zeta\in S$.

Let $\zeta_1\in S$, $r=\dist(\zeta_1,\partial S)$ and choose an arbitrary $\zeta_2\in{\rm B}(\zeta_1,r/2)$. Then
\begin{align*} 
\log^{\zeta_2}z=\log^{\zeta_1}z+\re\log\big(1+\frac{\zeta_1-\zeta_2}{z-\zeta_1}\big)
\end{align*}
where the function $\re\log\big(1+\frac{\zeta_1-\zeta_2}{z-\zeta_1}\big)$ goes to zero as $|z|\to\infty$ and is harmonic on the set $\{w\in S: \dist(w,\partial S)<r/2\}\cup S^{\rm c}$, which can be easily seen from $|(\zeta_1-\zeta_2)/(z-\zeta_1)|<1$. Thus, the Neumann jump of this function is zero, which gives
\begin{align*} 
\mathcal{N}(\log^{\zeta_1})=\mathcal{N}(\log^{\zeta_2}).
\end{align*}
By iterating this argument, we conclude that for any $\zeta_1,\zeta_2\in S$, the Neumann jumps of $\log^{\zeta_1}$ and $\log^{\zeta_2}$ coincide.\\

\noindent\textit{$\omega^{\infty}$ is non-negative.} Fix an arbitrary point $\zeta\in S$. Let $f=\log^{\zeta}-(\log^{\zeta})^S$ on $S^{\rm c}$. Note that $f$ is harmonic on $S^{\rm c}$, $f=0$ on $\partial S$ and $f(z)$ goes to $\infty$ as $|z|\to\infty$. By the maximum principle, $f$ attains its minimum value on $\partial S$. Thus, $f(z)\geq 0$ for all $z\in S^{\rm c}$. Thus, $\frac{\partial f}{\partial n}\geq 0$ on $\partial S$.\\

\noindent\textit{$\omega^{\infty}$ is a probability measure.} What remains is to verify that $\int_{\partial S}\rd\omega^{\infty}$ is $1$. This follows from the Green's theorem,
\begin{gather*} 
\int_{\partial S}\frac{\partial \log^{\zeta}}{\partial n}\rd s=\int_{S}1\Delta \log^{\zeta}\rd m+\int_{S}\nabla1\cdot\nabla\log^{\zeta}\rd m=2\pi, 
\\
\int_{\partial S}\frac{\partial (\log^{\zeta})^S}{\partial (-n)}\rd s=\lim_{R\to\infty}\Big(\int_{{\rm B}(0,R)\setminus S}1\Delta(\log^{\zeta})^S\rd m+\int_{{\rm B}(0,R)\setminus S}\nabla1\cdot\nabla(\log^{\zeta})^S\rd m-\int_{\partial {\rm B}(0,R)}\frac{\partial (\log^{\zeta})^S}{\partial n}\rd s\Big)=0.
\end{gather*}
\\

\noindent\textit{Well-definedness of \eqref{eqn:capacity_defn}.} By the second equality in \eqref{eqn:lop_capacity}, it suffices to show that $\int_{\partial S}\log^{\zeta_1}\rd \omega^{\infty}=\int_{\partial S}\log^{\zeta_2}\rd \omega^{\infty}$ for any $\zeta_1,\zeta_2$ distinct in the interior of $S$. We begin with a symmetric expression, and apply Green's theorem again,
\begin{multline*} 
\int_{\C}\nabla(\log^{\zeta_1})^{S}\cdot\nabla(\log^{\zeta_2})^{S}\rd m=\int_{S}\nabla(\log^{\zeta_1})^{S}\cdot\nabla(\log^{\zeta_2})^{S}\rd m+\int_{S^{\rm c}}\nabla(\log^{\zeta_1})^{S}\cdot\nabla(\log^{\zeta_2})^{S}\rd m
\\
=\int_{\partial S}\log^{\zeta_1}\frac{\partial \log^{\zeta_2}}{\partial n}\rd s-\int_{S}\log^{\zeta_1}\Delta \log^{\zeta_2}\rd m+\int_{\partial S}\log^{\zeta_1}\frac{\partial (\log^{\zeta_2})^S}{\partial (-n)}\rd s=-2\pi\log|\zeta_1-\zeta_2|+2\pi\int_{\partial S}\log^{\zeta_1}\rd \omega^{\infty}.
\end{multline*}
Note that the right-hand side is not symmetric with respect to $\zeta_1$ and $\zeta_2$, which implies
\begin{align*} 
\int_{\partial S}\log^{\zeta_1}\rd \omega^{\infty}=\int_{\partial S}\log^{\zeta_2}\rd \omega^{\infty}.
\end{align*}

\begin{remark}
The properties of the harmonic measure stated in the introduction and proved above can easily be extended to the case of arbitrary connected droplet $S$ with smooth boundary, without the simple connectivity assumption. 
In this general case,  independence of $\omega^{\infty}$ from the choice of $\zeta$ is proved in the same manner, $\omega^{\infty}$ is a measure with mass $1$ which is non-negative and supported on the outer boundary of $S$,  and $\int_{\partial S}\log^{\zeta}\rd\omega^{\infty}$ does not depend on $\zeta$, again.
\end{remark}

\tableofcontents

\end{document}